\documentclass{amsart}
\pdfoutput=1

\usepackage{amsmath,amsthm,amsfonts,amssymb,verbatim,dsfont}
\usepackage{hyperref, cite, fullpage, graphicx, subfigure,ifpdf, breakurl}

\newcommand{\abs}[1]{\lvert#1\rvert}
\newcommand{\norm}[1]{\lVert#1\rVert}
\newcommand{\paren}[1]{\left(#1\right)}
\newcommand{\bracket}[1]{\left[#1\right]}
\newcommand{\set}[1]{\left\{#1\right\}}

\newcommand{\loc}{\mathrm{loc}}
\newcommand{\Real}{\mathrm{Re}}
\newcommand{\Imag}{\mathrm{Im}}

\newcommand{\bigo}{\mathrm{O}}
\newcommand{\littleo}{\mathrm{o}}

\newcommand{\eps}{\epsilon}

\newcommand{\ddr}{\frac{d}{dr}}

\newcommand{\ind}{{\mathrm{ind}}}
\newcommand{\rad}{{\mathrm{rad}}}

\newcommand{\rmax}{r_{\max}}

\newcommand{\spn}{\mathrm{span}}

\newcommand{\inner}[2]{\left\langle #1,#2\right\rangle}

\newcommand{\R}{\mathbb{R}}

\newcommand{\calK}{\mathcal{K}}

\newcommand{\calL}{\mathcal{L}}
\newcommand{\calB}{\mathcal{B}}
\newcommand{\calV}{\mathcal{V}}
\newcommand{\calH}{\mathcal{H}}

\newcommand{\sech}{\mathrm{sech}}

\newtheorem{prop}{Proposition}[section]

\newtheorem{defn}[prop]{Definition}
\newtheorem{definition}[prop]{Definition}

\newtheorem{theorem}{Theorem}[section]
\newtheorem{lemma}[theorem]{Lemma}
\newtheorem{corollary}[theorem]{Corollary}
\theoremstyle{remark}
\newtheorem{remark}[theorem]{Remark}

\newcommand{\grad}{\nabla}
\newcommand{\laplacian}{\Delta}

\newcommand{\real}{\mathbb{R}}
\newcommand{\integers}{\mathbb{Z}}
\newcommand{\x}{\overline{x}}
\newcommand{\J}{{\mathcal J}}
\newcommand{\epsNorm}{\int{\abs{\grad_y\epsilon}^2\,dy} + \int_{\abs{y}\leq\frac{10}{b}}{\abs{\epsilon}^2e^{-\abs{y}}\,dy}}
\newcommand{\epsTilde}{\widetilde{\epsilon}}
\newcommand{\epsTildeNorm}{\int{\abs{\grad_y\widetilde{\epsilon}}^2\,dy} + \int_{\abs{y}\leq\frac{10}{b}}{\abs{\widetilde{\epsilon}}^2e^{-\abs{y}}\,dy}}

\newcommand{\m}{ {(m)} }
\newcommand{\Hm}{H^1_\m}
\newcommand{\dotHm}{\dot{H}^1_\m}
\newcommand{\calFm}{{\mathcal F}_\m}
\newcommand{\Qm}{Q^\m}
\newcommand{\Qmbar}{\overline{Q}^\m}
\newcommand{\Rm}{R^\m}
\newcommand{\Qmb}{Q^\m_b}
\newcommand{\Pmb}{P^\m_b}
\newcommand{\Pm}{P^\m} 
\newcommand{\PmbT}{\widetilde{P}^\m_b}
\newcommand{\PmT}{\widetilde{P}^\m} 
\newcommand{\QmbT}{\widetilde{Q}^\m_b}
\newcommand{\QmT}{\widetilde{Q}^\m} 
\newcommand{\Zmb}{\zeta^\m_b}
\newcommand{\ZmbT}{\widetilde{\zeta}^\m_b}
\newcommand{\ZmbTbar}{\overline{\widetilde{\zeta}}^\m_b}
\newcommand{\Lm}{L^\m} 
\newcommand{\calLm}{\calL^\m} 
\newcommand{\calPm}{{\mathcal P}^\m}
\newcommand{\calHm}{\calH^\m}
\newcommand{\hyp}{\mathrm{hyp}}

\begin{document}
\title[Vortex Collapse]{Vortex Collapse for the $L^2$-Critical Nonlinear Schr\"odinger Equation}
\author{G. Simpson \& I. Zwiers}
\date{\today}

\begin{abstract}
  The focusing cubic nonlinear Schr\"odinger equation in two dimensions
  admits vortex solitons, standing wave solutions with spatial
  structure, $\Qm(r,\theta) = e^{im\theta}\Rm(r)$. In the case of spin
  $m=1$, we prove there exists a class of data that collapse with the
  vortex soliton profile at the log-log rate. This extends the work of
  Merle and Rapha\"el, (the case $m=0$,) and suggests that the $L^2$
  mass that may be concentrated at a point during generic collapse may
  be unbounded.  Difficulties with $m\geq 2$ or when breaking the spin
  symmetry are discussed.
\end{abstract}

\maketitle

\tableofcontents

\section{Introduction}
We consider the $L^2$-critical nonlinear Schr\"odinger equation in two
dimensions,
\begin{equation}\label{Eqn-NLS}
  \left\{\begin{aligned}
      &iu_t+\laplacian u + u\abs{u}^2 = 0\\
      &u(0,x) = u_0 \in H^1(\R^2).
    \end{aligned}\right.
\end{equation}

Equation (\ref{Eqn-NLS}) is locally wellposed for data $u_0\in H^1$,
\cite{GinibreVelo-NLSCauchyProb-1979,Kato-NLS-87}. That is, there
exists a solution $u \in C\left([0,T_{\max}),H^1\right)$ and some
fixed negative power so that $T_{\max} \geq T_{lwp}
=\norm{u_0}_{H^1}^{-C}$. Therefore, we have the classic blowup
alternative,
\[\begin{aligned}
  T_{\max} = +\infty &&\text{ or, }&& \lim_{t\to
    T_{\max}}\norm{u(t)}_{H^1} = +\infty.
\end{aligned}\] Evolution of $u_0$ by equation (\ref{Eqn-NLS})
preserves the following quantities.
\begin{align}
  \label{ConserveMass}
  M[u_0] = M[u(t)] &= \int_{\R^{2}}{\abs{u(t,x)}^2\,dx}, &&
  \text{(mass)}
  \\
  \label{ConserveEnergy}
  E[u_0] = E[u(t)] &= \int{\abs{\grad_xu(t,x)}^2\,dx} -
  \frac{1}{2}\int{\abs{u(t,x)}^4\,dx}, && \text{(energy)}
  \\
  \label{ConserveMoment}
  P[u_0] = P[u(t)] &= \Imag\left(\int{\overline{u}(t,x)\grad
      u(t,x)\,dx}\right).  && \text{(momentum)}
\end{align}
The associated symmetries of the equation are phase, time translation,
and spatial translation. There is a Galilean symmetry,
\[\begin{aligned}
  u_{\beta_0}(t,x) = u(t,x-\beta_0
  t)e^{i\frac{\beta_0}{2}\cdot\left(x-\frac{\beta_0}{2}t\right)},
  &&\text{for any fixed }\beta_0\in\R^2,
\end{aligned}\] and a scaling symmetry,
\[\begin{aligned}
  u_{\lambda_0}(t,x) = \lambda_0 u(\lambda_0^2t,\lambda_0 x)
  &&\text{for any fixed }\lambda_0 > 0.
\end{aligned}\] The effect of scaling on Sobolev norms is,
$\norm{u_{\lambda_0}}_{\dot{H}^s} =
\lambda_0^{-s}\norm{u}_{\dot{H}^s}$, for any reasonable $s$. Note that
only the critical norm is left invariant. By choosing $\lambda_0 =
\norm{u(t)}_{H^1}$ at a fixed time, and using the minimum local
wellposedness time for unit data in $H^1$, we have the scaling lower
bound for the blowup speed,
\[\begin{aligned}
  u(t)\in C\left([0,T_{\max}), H^1\right), \text{ with }T_{\max}\text{
    maximal, then} &&\norm{u(t)}_{H^1} \gtrsim
  \frac{1}{\sqrt{T_{\max}-t}}.
\end{aligned}\]
Alternatively, the scaling lower bound can be established through
energy conservation, \cite{CW-CauchyProblemHs-90}.

Peculiar to the $L^2$-critical case, there is also the
pseudo-conformal (or lens) symmetry,
\begin{equation}\label{Defn-Eqn-PseudoConformal}
  v(t,x) = \frac{1}{T-t}u\left(\frac{1}{(T-t)^2},\frac{x}{T-t}\right)e^{-i\frac{\abs{x}^2}{4(T-t)}},
\end{equation}
which acts on the virial space, $\left\{ f\in H^1 \right\}\cap\left\{
  \abs{x}^2f \in L^2\right\}$. In particular, the pseudo-conformal
symmetry transforms standing wave solutions into blowup solutions with
$H^1$ norm growth $\frac{1}{T-t}$.

\subsection{Blowup with Soliton Profile}
To find standing wave solutions of equation (\ref{Eqn-NLS}), introduce
the usual ansatz, $u(t,x) = e^{it}Q(x)$, to derive the profile
equation,
\begin{equation}\label{Eqn-Defn-Q}
  \begin{aligned}
    \laplacian Q - Q + Q\abs{Q}^2 = 0.
  \end{aligned}
\end{equation}
There is a unique real-valued positive radial solution $Q$ to equation
(\ref{Eqn-Defn-Q}), as proved by McLeod and Serrin
\cite{McLeodSerrin-Uniqueness-87}\footnotemark.
\footnotetext{Following earlier work by Coffman
  \cite{Coffman-Uniqueness3DCubicGroundState-72} in 3D. Kwong
  \cite{Kwong-Uniqueness-89} extended the result to all
  $H^1$-subcritical nonlinearities.}  This solution we call the
soliton, or the ground-state since $E(Q) = 0$.
In this paper we will focus on other solutions of equation
(\ref{Eqn-Defn-Q}), as we discuss in the next
section. 
Weinstein \cite{Weinstein-NLSSharpInterpolation-82} identified the
soliton as the unique minimizer of $J[f] = \frac{\abs{\grad
    f}_{L^2}^2\abs{f}_{L^2}^2}{\abs{f}_{L^4}^4}$ among $H^1$
functions, thereby showing the optimal constant of the
Gagliardo-Nirenberg inequality,
\[
\norm{v}_{L^4}^4 \leq
\frac{2}{\norm{Q}_{L^2}^2}\norm{v}_{\dot{H}^1}^2\norm{v}_{L^2}^2.
\]
Note that if $M[u_0] < M[Q]$, the Gagliardo-Nirenberg inequality gives
apriori control of the $H^1$ norm from the conservation of
energy. That is, there is global wellposedness for data with $M[u_0] <
M[Q]$.

The pseudo-conformal transformation (\ref{Defn-Eqn-PseudoConformal})
applied to the standing wave solution $e^{it}Q(x)$ gives an explicit
blowup solution with $M[u_0] = M[Q]$. We denote this explicit solution
$S(t)$; Merle \cite{Merle-BlowupSolnMinimalMass-93} showed that, up to
symmetries, it is the only blowup solution with the mass of
$Q$. Bourgain and Wang \cite{BourgainWang} proved that $S(t)$ is
stable with respect to perturbations that are exceptionally flat near
the central profile.

More generally, negative energy data in the virial space leads to
blowup, as shown by Glassey's virial identity
\cite{Glassey-BlowupUpSolnsCauchyProbNLS-77},
\[
\frac{d^2}{dt^2}\int{\abs{x}^2\abs{u(t)}^2} =
4\frac{d}{dt}{Im}\int{x\cdot\grad u\overline{u}} = 16E[u_0]
\]
Ogawa and Tsutsumi
\cite{OgawaTsutsumi-NegEnerBlowupForEnergySubcrit-91} later extended
the argument to negative energy radial data.

Let us consider ${\mathcal B}_{\alpha} = \{u_0\in H^1 : M[Q] < M[u_0]
< M[Q]+\alpha\}$, where $\alpha>0$ is some small constant.  Merle and
Rapha\"el \cite{MR-UniversalityBlowupL2Critical-04} proved that there
is no solution in ${\mathcal B}_\alpha$ that blows up as predicted by
Glassey's virial identity\footnotemark.  \footnotetext{There is no
  solution in ${\mathcal B}_\alpha$ for which $\lim_{t\to
    T_{\max}}\int{\abs{x}^2\abs{u(t)}^2} = 0$, in constrast to the
  explicit solution $S(t)$.}  They also showed
\cite{MR-SharpUpperL2Critical-03,MR-SharpLowerL2Critical-06} that
there is an open subset ${\mathcal O}\subset {\mathcal B}_\alpha$,
including all the negative energy data, that lead to blowup in finite
time with the log-log rate,
\[
\norm{u(t)}_{H^1} \approx \sqrt{\frac{\log\abs{\log(T-t)}}{T-t}}.
\]
Rapha\"el \cite{R-StabilityOfLogLog-05} proved that all solutions in
${\mathcal B}_\alpha$ that lead to blowup either belong to ${\mathcal
  O}$, or blowup with at least the $H^1$ growth rate of
$S(t)$. Finally, Merle and Rapha\"el \cite{MR-ProfilesQuantization-05}
showed that all solutions in ${\mathcal B}_\alpha$ that blowup
concentrate exactly the profile $Q$ at a point, in the sense that
there are parameters $\lambda(t) > 0$, $\gamma(t)\in\R$ and
$\x(t)\in\R^2$ such that,
\[
u(t,x) -
\frac{1}{\lambda(t)}Q\left(\frac{x-\x(t)}{\lambda(t)}\right)e^{-i\gamma(t)}
\longrightarrow u^*(x),
\]
where the convergence is in $L^2$ as $t\to T_{\max}$. Moreover, the
residual profile $u^*$ identifies the blowup regime, with $u^*\not\in
H^1$ if and only if the solution belonged to ${\mathcal O}$ and
followed the log-log rate.

\subsection{Vortex Solitons}\label{SubSec-VortexSolitons}
Vortex solitons are solutions to equation (\ref{Eqn-Defn-Q}) of the
form $\Qm(r,\theta) = e^{im\theta}\Rm(r)$, where $\Rm$ is real-valued
and positive. That is, we seek a function $\Rm$ that satisfies,
\begin{equation}\label{Eqn-VortexSoliton-R}
  \left\{\begin{aligned}
      &\laplacian \Rm - \left(1+\frac{m^2}{r^2}\right)\Rm + \left(\Rm\right)^3 = 0,\\
      &\begin{aligned}
	\partial_r\Rm|_{r=0} = 0, && \Rm(\abs{x})>0, && \Rm \in
        H^1(\R^2)\cap\left\{{\abs{x}}^{-1}f(x)\in L^2\right\}
      \end{aligned}
    \end{aligned}\right.
\end{equation}
For all $m\in\integers$, Iaia and Warchall
\cite{IaiaWarchall-NonradialSolns-95} showed there exists a solution
to (\ref{Eqn-VortexSoliton-R}) and, analogous to the result of Kwong
\cite{Kwong-Uniqueness-89} in the case $m=0$, Mizumachi
\cite{Mizumachi-VortexSolitons-05} has shown it is unique.
Fibich and Gavish \cite[Lemma 12]{FG-TheorySingularVortex-08} have
remarked that the resulting profile $Q_m$ is the unique minimizer of
$J[f] = \frac{\abs{\grad f}_{L^2}^2\abs{f}_{L^2}^2}{\abs{f}_{L^4}^4}$
among $H^1$ functions with spin $m$. We denote this space by $\Hm$.
Some vortex solutions are pictured in Figure \ref{f:vortex_surfaces},
and their radial profiles appear in Figure \ref{f:vortex_radial}.

\ifpdf
\begin{figure}
  \centering
  \subfigure[]{\includegraphics[width=2.35in]{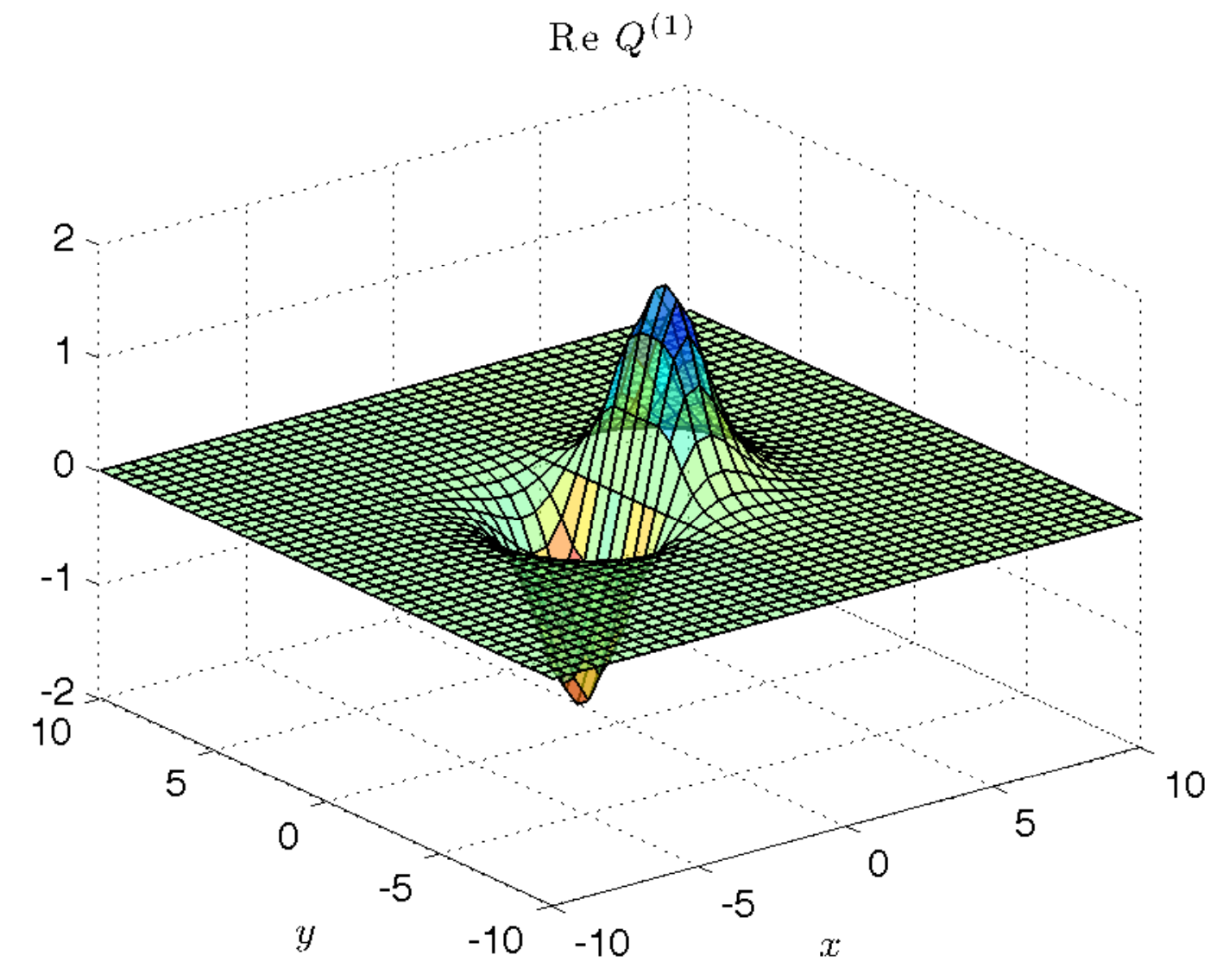}}
  \subfigure[]{\includegraphics[width=2.35in]{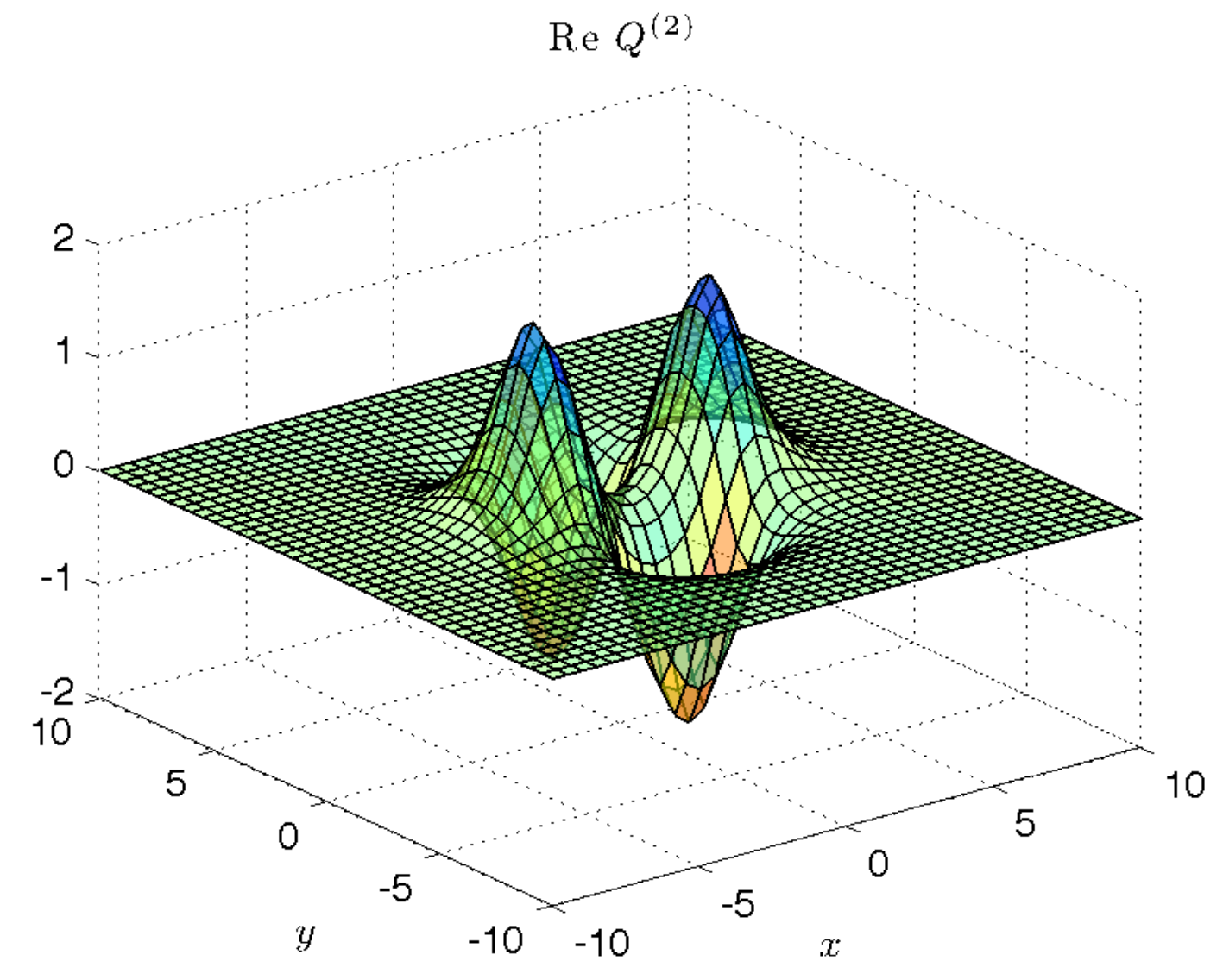}}
	
  \subfigure[]{\includegraphics[width=2.35in]{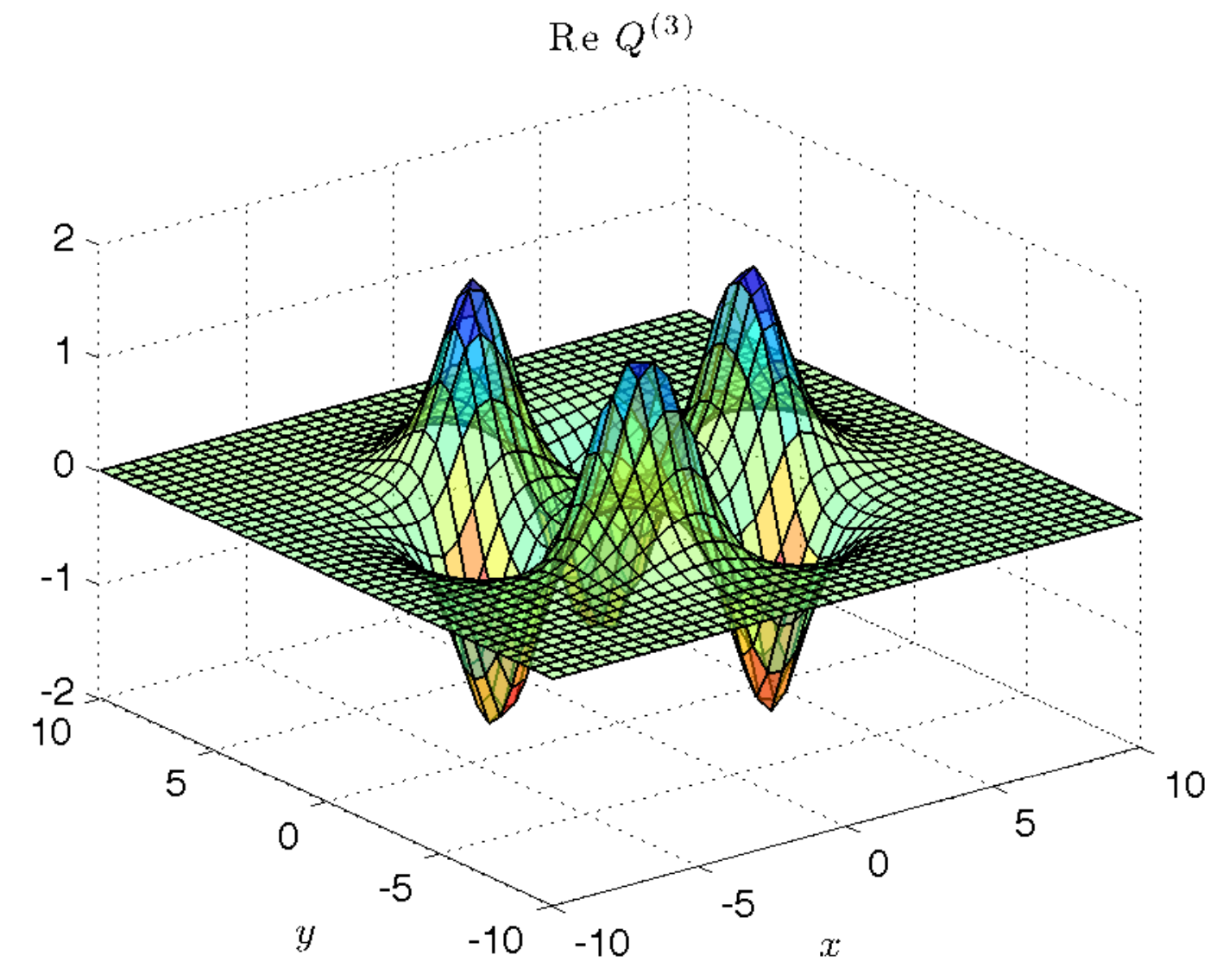}}
  \caption{The real component of some vortex solutions.}
  \label{f:vortex_surfaces}
\end{figure}
	
	\begin{figure}
          \centering
          {\includegraphics[width=3in]{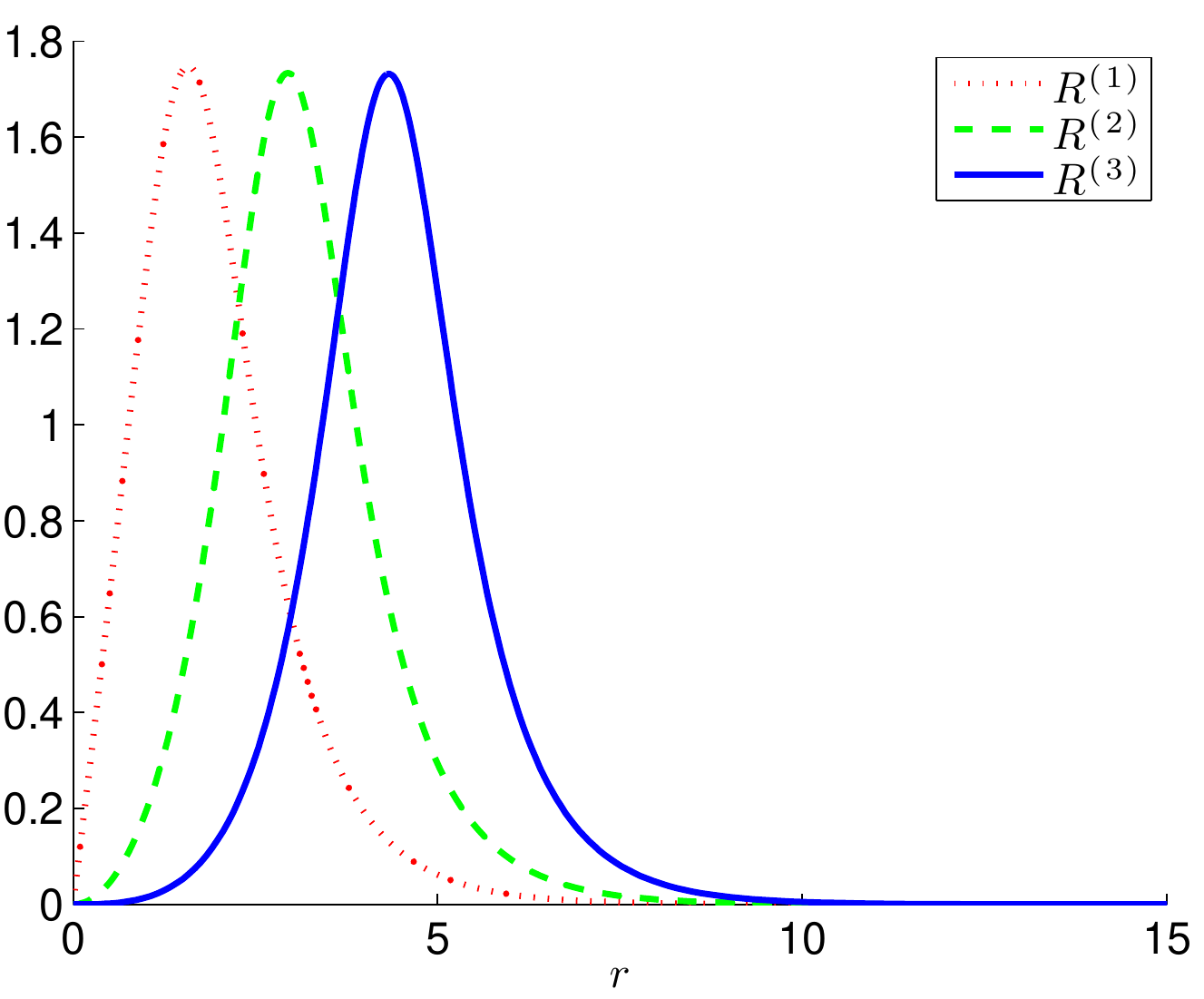}}
          \caption{The radial component of some vortex solutions.}
          \label{f:vortex_radial}
	\end{figure}
        \else (No figures in DVI mode) \fi

        This variational characterization gives an optimal
        Gagliardo-Nirenberg inequality
        for functions in $\Hm$. As a consequence, for data $u_0 \in
        \Hm$ and $L^2$ norm less than $\norm{\Qm}_{L^2}$ there is
        global wellposedness. As a second consequence, Fibich and
        Gavish \cite[Corollary 16]{FG-TheorySingularVortex-08} remark
        that $\norm{\Qm}_{L^2}^2$ is a strictly increasing sequence in
        $m$. Indeed, Pego and Warchall
        \cite{PegoWarchall-VorticesNLS-02} showed the asymptotic form,
        \[
        \Rm(r) \approx \left(1 +
          \frac{m^2}{r_{\max}^2}\right)^\frac{1}{2}\sqrt{2}\sech\left(\left(1+\frac{m^2}{r_{\max}^2}\right)^\frac{1}{2}(r-r_{\max})\right),
        \]
        where $r_{\max} \approx \sqrt{2}m$ for $m \gg 0$.
        Therefore, $\norm{\Qm}_{L^2}^2 \approx 4\sqrt{3} m$ for large
        $m$, which Fibich and Gavish found to be a good
        approximation\footnotemark even for small $m$.
        \footnotetext{Error less than 3\% for $m=2$, less than 0.4\%
          for $m\geq 5$.}

        The linearization of equation (\ref{Eqn-NLS}) near $\Qm$ is,
        \begin{equation}\label{Eqn-Linear1}\begin{aligned}
            \partial_tv = -i\Lm[v], &&\text{ where,}&& \Lm[v] \equiv
            \left(-\laplacian + 1 - \abs{\Qm}^2\right)v -
            2\left(\Qm\right)^2\overline{v}.
          \end{aligned}\end{equation}
        Written as a harmonic series,
        $v=\sum_{j\in\integers}{e^{i(m+j)\theta}f_j(r)}$,
        \begin{equation}\label{Eqn-Linear2}\begin{aligned}
            \Lm[v] = \sum_{j}\left(-\laplacian + 1 -
              \abs{\Qm}^2\right)e^{i(m+j)\theta}f_j -
            2\abs{\Qm}^2e^{i(m-j)\theta}\overline{f_j},
          \end{aligned}\end{equation}
        so that it is clear the linear system excites harmonics in
        pairs. In the case involving only $j=0$, that is, $v =
        e^{im\theta}\left(v_1 + iv_2\right)$, we may write $-i\Lm[v]$
        in matrix form as,
        \begin{equation}\label{Eqn-Defn-L}\begin{aligned}
            \left[\begin{matrix}0&\Lm_-\\-\Lm_+&0\end{matrix}\right]\left[\begin{matrix}v_1\\v_2\end{matrix}\right]
            &&\text{ where, }\;
            \begin{aligned}\Lm_+ &= -\laplacian + 1 - 3\abs{\Qm}^2,\\
              \Lm_- &= -\laplacian + 1 - \abs{\Qm}^2.\end{aligned}
          \end{aligned}
        \end{equation}
        Comparing equations (\ref{Eqn-Linear2}) and (\ref{Eqn-Defn-L})
        we see that $\Lm$ takes on the form of (\ref{Eqn-Defn-L}) on
        all of $H^1$ in the case of spin $m=0$.  In this important
        case, Weinstein \cite{Weinstein85} showed that the generalized
        nullspace of $L$ has dimension $8$ and is generated by the
        symmetries.

        In the cases $m=1$ and $m=2$, the generalized nullspace of
        $\Lm$ is generated in the same way. However, in these cases
        Pego and Warchall \cite{PegoWarchall-VorticesNLS-02} found
        unstable eigenvalues and additional eigenvalues in the
        spectral gap (all for modes with $\abs{j} \neq 0$).
        That is, there exists a function $\rho$ with spin $m=1$ such
        that,
        \[\begin{aligned}
          &\Lm_+\left(\rho\right) = -\abs{y}^2Q_m,
          &&\Lm_-\left(\abs{y}^2Q_m\right) = -4\Lambda Q_m, \\
          &\Lm_+\left(\Lambda Q_m\right) = -2Q_m,
          &&\Lm_-\left(Q_m\right) = 0,
        \end{aligned}\] where $\Lambda = 1 + y\cdot\grad$ denotes the
        scaling operator. The remaining Jordan chains, generated by
        $\grad \Qm$, consist of functions with $\abs{j} = 1$.
        Instability of vortex profiles is not restricted to the cubic
        nonlinearity. Mizumachi \cite{Mizumachi-InstabilityVortex-07}
        has shown that there are unstable vortex profiles for any
        power-type nonlinearity strictly stronger than linear.

        \subsection{Blowup with Vortex Profiles}

        Any vortex soliton becomes a blowup solution through the
        pseudo-conformal transformation. Study of the asymptotic
        profile during vortex blowup was initiated by Fibich \& Gavish
        \cite{FG-TheorySingularVortex-08}, including the variational
        structure referenced above. Their work includes numerical
        simulations where they found data with mass slightly larger
        than $\Qm$ that blowup at exactly the scaling lower bound and
        with profiles different from the vortex
        soliton.
          \footnote{In particular, they present results using
          $u_0 = 1.02 \Qm(r,\theta)$ and spin $m=2$. 
          The profiles identified, denoted $G_{m}$, are truncated solutions of equation (\ref{Eqn-Pmb}) with Cauchy boundary conditions and an
          implied value of $b$, in this case $b\approx 0.1092$. 
Our own truncated solutions of \eqref{Eqn-Pmb} are very similar.
          Fibich \& Gavish have conveyed by personal communication corresponding discoveries for spin $m=1$ and data as small as $1.00001\Qm$.

          }

        Our main result is that there is a class of solutions with
        spin $m=1$ that blowup with exactly the vortex soliton profile
        and log-log behaviour similar to that established in the case
        $m=0$.
        \begin{theorem}[Log-log Blowup with Vortex
          Profile]\label{Thm-MainResult}
          Assume the Spectral Property\footnotemark is true for spin
          $m$.  Then there exists a class of data $\calPm$, open as a
          subset of $\Hm$, such that for $u_0\in \calPm$ the evolution
          $u(t)$ by (\ref{Eqn-NLS}) blows up at finite time $T_{\max}$
          with the $\Qm$ profile and log-log rate.  That is, for
          $t\in[0,T_{\max})$ there exist continuously variable
          parameters $\lambda(t) > 0$ and $\gamma(t)\in\R$ with the
          following properties: \footnotetext{See Proposition
            \ref{Prop-SpectralProperty}, below.}
          \begin{description}
          \item[\it Log-log Blowup Rate:]
            \begin{equation}\label{Thm-MainResult-BlowupRate}
              \lim_{t\to T_{\max}} \norm{u(t)}_{\dot{H}^1}\sqrt{\frac{T_{\max}-t}{\log\abs{\log T_{\max}-t}}} = C
            \end{equation}
          \item[\it Description of the Singularity:]
            \begin{equation}\label{Thm-MainResult-L2}
              \lim_{t\to T_{\max}} u(t,x) - \frac{1}{\lambda(t)}\Qm\left(\frac{x}{\lambda(t)}\right)e^{-i\gamma(t)} = u^*(x) \in L^2(\R^2).
            \end{equation}
          \end{description}
        \end{theorem}

        We will now discuss the consistency of the self-similar regime
        discovered by Fibich \& Gavish and Theorem
        \ref{Thm-MainResult}.  Consider,
        \[
        B_{\alpha,m} = \left\{u_0\in \Hm : M[\Qm] < M[u_0] < M[\Qm] +
          \alpha\right\}.
        \]
        Then, due to the variational characterization of $\Qm$:
        \begin{theorem}[``Orbital
          Stability'']\label{Thm-OrbitalStability}
          For $\alpha>0$ sufficiently small, let $v\in B_{\alpha,m}$
          with, $E[v] \leq \alpha\norm{v}_{H^1}^2$. Then there exists
          $\lambda_v>0, \gamma_v\in\real$ such that,
          \[
          \norm{\lambda_v\,v\left(\lambda_vy\right)e^{i\gamma_v}-\Qm(y)}_{H^1}
          \leq \delta(\alpha),
          \]
          where $\delta(\alpha)\to 0$ as $\alpha\to 0$.
        \end{theorem}
        The proof is by means of concentration compactness and the
        Gagliardo Nirenberg inequality in $\Hm$, and is not
        constructive. See \cite[Theorem 6]{R-Zurich} for a clear
        exposition. The class of data $\calPm$ from Theorem 
        \ref{Thm-MainResult} belongs to $B_{\alpha,m}$, and we note that
        the orbital stability of Theorem \ref{Thm-OrbitalStability} applies 
        to all data $v_0\in B_{\alpha,m}$ that blowup in finite time.
        Indeed, we conjecture\footnotemark that finite time blowup
        solutions from the class $B_{\alpha,m}$ either obey the
        log-log blowup rate (\ref{Thm-MainResult-BlowupRate}) or the
        lower bound, $\norm{v(t)}_{H^1}\gtrsim (T_{max}-t)^{-1}$.
        \footnotetext{We expect the analysis of
          \cite{R-StabilityOfLogLog-05} to apply, and that the proof
          of Theorem \ref{Thm-MainResult} may be reformulated to apply
          to all $v_0\in B_{\alpha,m}$ with $\norm{v(t)}_{H^1}\in
          L^1\left(t\in[0,T_{max})\right)$, as in
          \cite{MR-SharpLowerL2Critical-06}.}

        Collapse at the square-root rate has also been observed numerically in
        the case with no spin, \cite{FGW-NewSingularSolns-05}.
        These examples are an important area of continuing study.
        It is possible that the threshold $\alpha$ of Theorem \ref{Thm-OrbitalStability} (and hence the applicability of Theorem \ref{Thm-MainResult}) is exceedingly small.

        \subsection{Spectral Propety}
        \label{s:intro_specprop}
        In order to demonstrate the dynamic claimed in Theorem
        \ref{Thm-MainResult}, we will attempt to parameterize the
        solution in terms of the symmetries and a suitable deformation
        of the profile $\Qm$. In order for the finite-dimensional
        system of parameters to capture the essential dynamics of the
        solution we require two things. First, that the parameter
        dynamics can be reliably predicted from a finite system of
        differential inequalities. Second, that after removing the
        central profile from the solution the error $\epsilon$ can be
        estimated in terms of those parameter dynamics.

        That the parameter dynamic are stable is an essential feature
        of the log-log regime. Indeed, Rapha\"el showed
        \cite{R-StabilityOfLogLog-05} that the relationship between a
        particular ratio of parameters\footnotemark and a fixed
        constant evolves according to a Riccati equation, with the
        log-log dynamic corresponding to the stable branch.
        \footnotetext{Namely the sign of $f_- =
          \frac{b}{\lambda}-d_0\sqrt{E_0}$. Parameter $b$ will be
          introduced in Section \ref{SubSec-DecompModulate}.}  To
        control the error $\epsilon$ in terms of the dynamics, we will
        consider the following operator, derived from the linearized
        energy,
        \begin{equation}\label{Eqn-Defn-H}\begin{aligned}
            \calHm(\epsilon,\epsilon) =
            \left\langle\calLm_1\epsilon_1, \epsilon_1\right\rangle +
            \left\langle\calLm_2\epsilon_2, \epsilon_2\right\rangle,
          \end{aligned}\end{equation}
        where
        \begin{align}
          \calLm_1 = -\laplacian + 3\Qm y\cdot\grad \Qmbar, &\quad
          \calLm_2 = -\laplacian + \Qm y\cdot\grad \Qmbar,\\
          \epsilon_1 = e^{im\theta}\,Re\left( e^{-im\theta}\epsilon
          \right),&\quad \epsilon_2 = e^{im\theta}\,Im\left(
            e^{-im\theta}\epsilon \right).
        \end{align}

        This decomposition, $\epsilon = \epsilon_1 + i\epsilon_2$, is
        powerful, as it reduces the algebraic structure of the problem
        in $H^1_m$ to that of the radially symmetric problem in $H^1$.
        For further discussion, see \eqref{Eqn-Defn-Notation}, below.

        We will prove the following for $m=1$,
        \begin{prop}[Spectral Property]
          \label{Prop-SpectralProperty}
          Let $\epsilon \in H^{1}_{(m)}$ Then there exists a universal
          constant $\delta_m$ such that
          \begin{equation}\label{Eqn-SpectralProperty}\begin{aligned}
              \calHm(\epsilon,\epsilon) \geq
              &\delta_m\left(\int{\abs{\grad_y\epsilon}^2} + \int{\abs{\epsilon}^2e^{-\abs{y}}}\right)\\
              &\quad - \frac{1}{\delta_m} \left( \left\langle
                  \epsilon_1,\Qm\right\rangle^2 +\left\langle
                  \epsilon_1,\Lambda \Qm\right\rangle^2\right. \\
              &\qquad \left.+\left \langle \epsilon_2,\Lambda
                  \Qm\right\rangle^2 +\left\langle
                  \epsilon_2,\Lambda^2 \Qm\right\rangle^2 \right).
            \end{aligned}\end{equation}
        \end{prop}
        In the case of $L^2$-critical nonlinearity, no spin, and
        dimension $N=1$, Merle and Rapha\"el \cite[Appendix
        A]{MR-BlowupDynamic-05} gave an explicit proof of the Spectral
        Property. In the case of $L^2$-critical nonlinearity, no spin,
        and dimensions $N=2,3,4,5$, including equation \eqref{Eqn-NLS}
        in the case $m=0$, Fibich, Merle and Rapha\"el
        \cite{FMR-ProofOfSpectralProperty-06} have given a numerical
        proof that inspires our own proof of Proposition
        \ref{Prop-SpectralProperty} in Section
        \ref{Section-SpectralProperty}.  Details of our numerical
        methods are provided in Appendix \ref{s:numerics}.  Code to
        reproduce our computations is available at
        \url{http://www.math.toronto.edu/simpson/files/vortex_dist.tgz}.
        As stated, the spectral property is false for $m=2, 3$.

\section{Proof of Log-log Blowup}
In this section, we prove Theorem \ref{Thm-MainResult} assuming
Proposition \ref{Prop-SpectralProperty}. Before decomposing the
solution, we introduce almost self-similar deformations of the vortex
profiles that simulate the effect of symmetries that do not belong to
$H^1$.  The standard self-similar ansatz is, $u(t,x) =
\frac{1}{\sqrt{2b(T-t)}}\Qmb\left(\frac{x}{\sqrt{2b(T-t)}}\right)e^{i\omega(t)}$,
which gives the following equation for the spatial profile,
\begin{equation}\label{Eqn-Qmb}
  \laplacian \Qmb - \Qmb +ib\Lambda \Qmb+ \Qmb\abs{\Qmb}^2=0.
\end{equation}
We seek solutions with spin $m$. Remove a quadratic phase,
$e^{im\theta}\Pmb(r) = \Qmb\,e^{ib\frac{r^2}{4}}$, and assume the
radial profile $\Pmb(r)$ is real valued. We seek solutions to,
\begin{equation}\label{Eqn-Pmb}
  \left\lbrace\begin{aligned}
      & \laplacian \Pmb - \left(1+\frac{m^2}{r^2}-\frac{b^2}{4}r^2\right)\Pmb + \left(\Pmb\right)^3 = 0,\\
      &\begin{aligned}
	\lim_{r\to 0^+}r^{-m}\Pmb(r) \neq 0,
	&& \lim_{r\to 0^+}\partial_r\left(r^{-m}\Pmb(r)\right) = 0.
      \end{aligned}
    \end{aligned}\right.
\end{equation}
As pointed out by Fibich and Gavish \cite[Lemma
8]{FG-TheorySingularVortex-08}, equation (\ref{Eqn-Pmb}) does not
admit solutions in either $L^2$ or $\dot{H}^1$, due to oscillations of
amplitude $r^{-1}$ outside the domain of uniform ellipticity of the
linear part. The argument is due to Johnson and Pan \cite{RX}.  We
truncate a solution of (\ref{Eqn-Pmb}) at an arbitrary point, chosen
to allow close approximation to the vortex profile $\Qm$. Define,
\begin{equation}\label{Eqn-Defn-Rb}\begin{aligned}
    R_b = \frac{1}{\abs{b}}\sqrt{2+2\sqrt{1+b^2m^2}} &\geq
    \frac{2}{\abs{b}}.
  \end{aligned}\end{equation}
\begin{prop}[Localized Self-Similar Profiles]\label{Prop-Qb}

  Let $a > C\eta>0$ where $C>0$ is a fixed constant and $a,\eta$ are
  sufficiently small parameters. Then for $\abs{b}>0$ sufficiently
  small, there exists $\QmbT \in H^1(\R^2)$, supported on $\abs{y} <
  (1-\eta)R_b$, with the following properties.
  \begin{itemize}
  \item {\it Simple Profile}:
    \begin{equation}\label{Prop-Qb-RealProfile}\begin{aligned}
        \QmbT = e^{im\theta}e^{-ib\frac{\abs{y}^2}{4}}\PmbT(\abs{y}),
        &&\text{for }\PmbT\text{ real-valued, non-negative.}
      \end{aligned}\end{equation}

  \item {\it Algebraic Proximity to $\Qm$}:
    \begin{equation}\label{Prop-Qb-Eqn}
      \laplacian\QmbT - \QmbT + ib\Lambda\QmbT + \QmbT\abs{\QmbT}^2 = -\Psi_b,
    \end{equation}
    for an error term $\Psi_b$, supported on
    $(1-\eta)^2R_b<\abs{y}<(1-\eta)R_b$, that satisfies the estimate,
    $\norm{P(y)\grad^k\Psi_b}_{L^\infty} \leq
    e^{-\frac{C(P)}{\abs{b}}},$
    for $k=0,1$ and any polynomial $P$.

  \item {\it Uniform Proximity to $\Qm$}:
    \begin{equation}\label{Prop-Qb-closeToQ}\begin{aligned}
        \left.\norm{e^{C\abs{y}}\left(\QmbT-\Qm\right)}_{C^{3}}
          +\norm{e^{C\abs{y}} \left(\frac{\partial}{\partial b}\QmbT +
              i\frac{\abs{y}^2}{4}\Qm\right)}_{C^2}
        \right.\longrightarrow 0 && \text{ as } && b \rightarrow 0.
      \end{aligned}\end{equation}

  \item {\it Supercritical Mass and Degenerate Energy}:
    \begin{equation}\label{Prop-Qb-Mass}\begin{aligned}
        \left. \frac{\partial}{\partial(b^2)}\norm{\QmbT}^2_{L^2}\right|_{b^2=0}
        = \frac{1}{4}\int{\abs{x}^2\abs{\Qm}^2}, && \text{ denoted }
        && d_m, && \text{ and,}
      \end{aligned}\end{equation}
    \begin{equation}\label{Prop-Qb-EnerMomentum}\begin{aligned}
        \abs{E\left[\QmbT\right]} \leq
        e^{-(1+C\eta)(1-a)\frac{\pi}{\abs{b}}}.
      \end{aligned}\end{equation}
  \end{itemize}
\end{prop}
The proof of Proposition \ref{Prop-Qb} is similar to that given by
Merle and Rapha\"el \cite{MR-SharpUpperL2Critical-03,
  MR-UniversalityBlowupL2Critical-04, MR-SharpLowerL2Critical-06} in
the case of $m=0$. An overview of the proof, and description of the
particular adaptations for $m\neq 0$, is given in Appendix
\ref{Appendix-ProofOfAlmostSelfSimilar}.

Later in the argument, Section \ref{SubSec-Lyapounov}, we will
introduce the linear radiation induced by the truncation error
$\Psi_b$. A quantity $\Gamma_b$, related to the decay of this
radiation, will be an important dynamical quantity, measuring the rate
of mass ejection from the singular region. At the time we formally
define $\Gamma_b$, Proposition \ref{Prop-Zb}, we will also prove the
following estimate,
\begin{equation}\label{Eqn-GammaBEstimate}\begin{aligned}
    e^{-(1+C\eta)\frac{\pi}{b}} \lesssim \Gamma_b \lesssim
    e^{-(1-C\eta)\frac{\pi}{b}}.
  \end{aligned}\end{equation}

\subsection{Decomposition \& Modulation}\label{SubSec-DecompModulate}

\begin{lemma}[Modulation Near $\Qm$]\label{Lemma-Modulation}
  Suppose that $v\in \Hm$ is close to $\Qm$, up to symmetries:
  \begin{equation}\label{Lemma-Mod-GeoDecomp}
    v(x) = 
    \frac{1}{\lambda_v}\left(\QmT_{b_v}+\epsilon_v\right)
    \left(\frac{x}{\lambda}\right)e^{-i\gamma_v},
  \end{equation}
  for some symmetry parameters $\lambda_v > 0$, $b_v > 0$ and
  $\gamma_v \in \R$ such that the error is comparably small,
  \begin{equation}\label{Lemma-Mod-EpsilonAssumption}
    \int{\abs{\grad_y\epsilon_v(y)}^2\,dy}
    +\int_{\abs{y}\leq\frac{10}{b_v}}{\abs{\epsilon_v}^2e^{-\abs{y}}\,dy}
    < \Gamma_{b_v}^\frac{1}{2},
  \end{equation}
  where $y$ denotes $\frac{x}{\lambda_v}$, and such that the deformed
  profile is sufficiently close to $\Qm$,
  \begin{equation}\label{Lemma-Mod-ParamAssumption}\begin{aligned}
      \lambda_v < \frac{1}{10} b_v &&\text{ and, }&& b_v < \alpha^*.
    \end{aligned}\end{equation}
  Then there are parameters $\lambda_0 > 0$, $b_0 > 0$ and $\gamma_0
  \in \R$, nearby in the sense,
  \begin{equation}\label{Lemma-Mod-NewParamAreClose}
    \abs{b_0-b_v}
    +\abs{\frac{\lambda_0}{\lambda_v} - 1} 
    \leq \Gamma_{b_0}^\frac{1}{5},
  \end{equation}
  and such that the error $\epsilon_0$ corresponding to these
  parameters,
  \begin{equation}\label{Lemma-Mod-DefnNewEpsilon}
    \epsilon_0(y) = \lambda_0\,v\left(\lambda_0 y\right)\,e^{i\gamma_0} - \QmT_{b_0},
  \end{equation}
  satisfies the following orthogonality conditions\footnotemark:
  \begin{equation}\label{Lemma-Mod-OrthogConditions}
    \Real\left\langle\epsilon_0,\abs{y}^2\QmT_{b_0}\right\rangle = 
    \Imag\left\langle\epsilon_0,\Lambda^2\QmT_{b_0}\right\rangle =
    \Imag\left\langle\epsilon_0,\Lambda\QmT_{b_0}\right\rangle = 0.
  \end{equation}
\end{lemma}
\footnotetext{ These orthogonality conditions were introduced
  \cite[Lemma 6]{MR-UniversalityBlowupL2Critical-04}, and lead to a
  better estimate on the phase parameter than achieved in
  \cite{MR-SharpUpperL2Critical-03}.  }
Let us reiterate and extend the notation alluded to by equation
(\ref{Eqn-Defn-H}),
\begin{equation}\label{Eqn-Defn-Notation}
  \begin{aligned}
    &\left.\begin{aligned}
        &\epsilon_1 = e^{im\theta}\text{Re}\left(e^{-im\theta}\epsilon\right)\\
        &\epsilon_2 =
        e^{im\theta}\text{Im}\left(e^{-im\theta}\epsilon\right)
      \end{aligned}\right\}
    &&\Longrightarrow \epsilon = \epsilon_1 + i\epsilon_2,\\
    &\left.\begin{aligned}
        &\Sigma = e^{im\theta}\text{Re}\left(e^{-ib\frac{\abs{y}^2}{4}}\PmbT\right)\\
        & \Theta =
        e^{im\theta}\text{Im}\left(e^{-ib\frac{\abs{y}^2}{4}}\PmbT\right)
      \end{aligned}\right\}
    &&\Longrightarrow \QmbT = \Sigma + i\Theta.
  \end{aligned}
\end{equation}
Products between the components of $\epsilon$ and $\QmbT$ behave as if
they were real-valued, as does the modulus, for example,
$\abs{\epsilon}^2 = \abs{\epsilon_1}^2+\abs{\epsilon_2}^2$.  Moreover,
since $\abs{y}^2$ and the scaling operator, $\Lambda = 1 +
y\cdot\grad_y$, are radial operators, the algebraic relations for
$\abs{y}^2\QmbT$ and $\Lambda\QmbT$ are exactly the same as the case
$m=0$, \cite[Proposition 9 (iii)]{MR-UniversalityBlowupL2Critical-04}.
In particular, one may verify that, $\calLm_1(\epsilon_1) =
\frac{1}{2}\left[\Lm_+(\Lambda\epsilon_1) -
  \Lambda(\Lm_+\epsilon_1)\right]$, is true regardless of $m$. This is
the essential relationship for Lemma \ref{Lemma-SpectralProperty2},
below.  In the notation of (\ref{Eqn-Defn-Notation}), the
orthogonality conditions of equation
(\ref{Lemma-Mod-OrthogConditions}) can be written,
\[\begin{aligned}
  &\left\langle\epsilon_1,\abs{y}^2\Sigma\right\rangle +
  \left\langle\epsilon_2,\abs{y}^2\Theta\right\rangle = 0,\\
  &\left\langle\epsilon_2,\Lambda^2\Sigma\right\rangle -
  \left\langle\epsilon_1,\Lambda^2\Theta\right\rangle = 0,\\
  &\left\langle\epsilon_2,\Lambda\Sigma\right\rangle -
  \left\langle\epsilon_1,\Lambda\Theta\right\rangle = 0.
\end{aligned}\] These are exactly the same form as in the case
$m=0$. Indeed, the proof of Lemma \ref{Lemma-Modulation}, an implicit
function argument, is identical. See \cite[Lemma 2]{R-Zurich} for a
clear exposition.  For $m=0$, the following Lemma was proven by Merle
and Rapha\"el \cite[equation
(116)]{MR-UniversalityBlowupL2Critical-04}, and the same proof applies
here.
\begin{lemma}\label{Lemma-SpectralProperty2}
  Let $\epsilon\in \Hm$, and assume the Spectral Property is
  true. Then,
  \begin{equation}\label{Eqn-SpectralProperty2}\begin{aligned}
      \left\langle\Lm_1\epsilon_1,\epsilon_1\right\rangle& -
      \frac{\left\langle\epsilon_1,\Lm_+\Lambda^2\Qm\right\rangle
        \left\langle\epsilon_1,\Lambda\Qm\right\rangle}{\norm{\Lambda\Qm}_{L^2}^2}
      \geq\\
      &\delta_m\left(\int{\abs{\grad_y\epsilon}^2} +
        \int{\abs{\epsilon}^2e^{-\abs{y}}}\right) - \frac{1}{\delta_m}
      \left( \left\langle \epsilon_1,\Qm\right\rangle^2 +\left\langle
          \epsilon_1, \abs{y}^2 \Qm\right\rangle^2 \right),
    \end{aligned}\end{equation}
\end{lemma}

\begin{definition}[Description of Initial Data]\label{Defn-P}
  Define $\calPm$ to be those functions $u_0\in \Hm$ for which there
  are parameters $\lambda_0 > 0$, $b_0>0$ and $\gamma_0\in\R$ that
  satisfy the following conditions. Let $\epsilon_0$ denote the error
  in approximating $u_0$ with these particular parameters,
  \begin{equation}\label{Defn-P-GeoDecomp}\begin{aligned}
      u_0(x) &= \frac{1}{\lambda_0}\left(\QmT_{b_0}+\epsilon_0\right)
      \left(\frac{x}{\lambda_0}\right)e^{-i\gamma_0}.
    \end{aligned}\end{equation}
  We require that the orthogonality conditions
  (\ref{Lemma-Mod-OrthogConditions}) are satisfied, that there is,

  \begin{description}

  \item[\it proximity to $\Qm$,]
    \begin{equation}\label{Defn-P-Proximity}\begin{aligned}
	\text{in }L^2:
	&&& 0 < b_0^2 + \norm{\epsilon}_{L^2}^2 < (\alpha^*)^2,\\
	\text{in }\dot{H}^1: &&& \int{
          \abs{\grad_y\epsilon_0(y)}^2\,dy}
        +\int_{\abs{y}\leq\frac{10}{b_0}}{\abs{\epsilon_0(y)}^2e^{-\abs{y}}\,dy}
        <\Gamma_{b_0}^\frac{6}{7},
      \end{aligned}\end{equation}

  \item[\it parameters consistent with the log-log rate,]
    \begin{equation}\label{Defn-P-loglog}\begin{aligned}
	e^{-e^\frac{2\pi}{b_0}} < \lambda_0 <
        e^{-e^{\frac{\pi}{2}\frac{1}{b_0}}}, &&\text{ and,}
      \end{aligned}\end{equation}

  \item[\it normalized energy,]
    \begin{equation}\label{Defn-P-EnergyMomentum}
      \lambda_0^2\abs{E_0} 
      < \Gamma_{b_0}^{10}.
    \end{equation}
  \end{description}
\end{definition}

\begin{remark}[$\calPm$ is Non-Empty]\label{Remark-PmNonempty}
  Choose $b_0$ and $\lambda_0$ to satisfy (\ref{Defn-P-Proximity}) and
  (\ref{Defn-P-loglog}). Let $f\in \Hm$ satisfy orthogonality
  conditions (\ref{Lemma-Mod-OrthogConditions}) with $\norm{f}_{H^1} =
  1$, $\inner{f}{\Qm} = 1$. Such an $f$ may be computed explicitly
  from $\Qm$. Note that $\left.\partial_\nu E[\Qm+\nu f]\right|_{\nu =
    0} = - \inner{F}{\Qm} = -1$, and therefore we may choose
  $\epsilon_0 = \nu f$ with $\nu$ of the order of $E[\Qmb]$ to satisfy
  (\ref{Prop-Qb-EnerMomentum}).
\end{remark}

For the remainder of this section, we consider a fixed representative
$u_0 \in \calPm$. By the continuity of the flow of (\ref{Eqn-NLS}) in
$H^1$, and Lemma \ref{Lemma-Modulation}, there exists continuous
functions $\lambda(t) > 0$, $b(t)>0$ and $\gamma(t)\in\R$ and some
maximal $T_{\hyp}\in(0,T_{\max}]$ such that the following relaxations
of (\ref{Defn-P-Proximity}), (\ref{Defn-P-loglog}) and
(\ref{Defn-P-EnergyMomentum}) hold for all $t\in[0,T_{\hyp})$:
\begin{equation}\label{Hypo-Proximity-L2}
  0 < b^2(t) + \norm{\epsilon(t)}_{L^2}^2 < (\alpha^*)^\frac{1}{5},
\end{equation}
\begin{equation}\label{Hypo-Proximity-H1}
  \int{ \abs{\grad_y\epsilon(t,y)}^2\,dy}
  +\int_{\abs{y}\leq\frac{10}{b(t)}}{\abs{\epsilon(t,y)}^2e^{-\abs{y}}\,dy}
  <\Gamma_{b(t)}^\frac{3}{4},
\end{equation}
\begin{equation}\label{Hypo-loglog}\begin{aligned}
    e^{-e^\frac{10\pi}{b(t)}} < \lambda(t) <
    e^{-e^{\frac{\pi}{10}\frac{1}{b(t)}}}, &&\text{ and,}
  \end{aligned}\end{equation}
\begin{equation}\label{Hypo-EnergyMomentum}
  \lambda^2(t)\abs{E_0} 
  < \Gamma_{b(t)}^{2}.
\end{equation}
Note that as a consequence of these hypotheses, we may apply Lemma
\ref{Lemma-Modulation} at any $t\in[0,T_{\hyp})$. Therefore, one of
the following occurs:
\begin{description}
\item[Case 1:] $T_{\hyp} < T_{\max}$, and one of the hypotheses fails
  at $t = T_{\hyp}$, or
\item[Case 2:] $T_{\hyp} = T_{\max}$, $b\to 0$ as $t\to T_{\max}$, and
  due to (\ref{Hypo-loglog}) we have blowup.
\end{description}
In this section we will show that {\bf Case 1} cannot occur. Then,
assuming {\bf Case 2}, we will derive the conclusions of Theorem
\ref{Thm-MainResult}.
\begin{remark}[Parameters]
  The parameter $\eta>0$, already introduced, relates to the cutoff
  and shape of the singular profile $\QmbT$. Parameter $a>0$, to be
  introduced in Section \ref{SubSec-Lyapounov}, will be related to a
  cutoff point of the linear radiation associated with $\QmbT$. The
  value of $\eta$ is determined by the value of $a$ so that the
  argument of Subsection \ref{SubSubSec-HypoH1} is successful.  At all
  times, $\alpha^*>0$ is assumed sufficiently small for all the
  appropriate constants to cooperate.
\end{remark}

\subsection{Conservation Laws \& Basic Estimates}
By substitution of the time-dependent version of the geometric
decomposition (\ref{Defn-P-GeoDecomp}), the conservation laws of
(\ref{Eqn-NLS}) and the orthogonality conditions
(\ref{Lemma-Mod-OrthogConditions}) lead to some basic estimates.

\begin{lemma}
  \label{Lemma-PrelimEstimatesConservLaws}
  For all $t\in[0,T_{\hyp})$,
  \begin{description}

  \item[\it due to conservation of mass,]
    \begin{equation}\label{Eqn-prelimMassEst}
      b^2 + \int{\abs{\epsilon}^2} \lesssim \left(\alpha^*\right)^\frac{1}{2},
    \end{equation}
  \item[\it due to conservation of energy,]
    \begin{equation}\label{Eqn-prelimEnerEst}\begin{aligned}
        &2\Real\left\langle\epsilon,\QmbT-ib\Lambda\QmbT -
          \Psi_b\right\rangle \sim \int{\abs{\grad_y\epsilon}^2\,dy}
        -3\int_{\abs{y}\leq\frac{10}{b}}{\abs{\Qm\epsilon_1}^2}
        -\int_{\abs{y}\leq\frac{10}{b}}{\abs{\Qm\epsilon_2}^2},\\
        &\begin{aligned}\text{with error of the order, }&&
          \Gamma_b^{1-C\eta} +
          \delta(\alpha^*)\left(\epsNorm\right).\end{aligned}
      \end{aligned}\end{equation}
  \end{description}
\end{lemma}
\begin{proof}
  To prove (\ref{Eqn-prelimMassEst}), expand the conservation of mass,
  \[
  \int{\abs{\QmbT}^2\,dy} - \int{\abs{\Qm}^2\,dy}
  +2Re\left\langle\epsilon,\QmbT\right\rangle +\int{\abs{\epsilon}^2}
  =\int{\abs{u_0}^2} - \int{\abs{\Qm}^2\,dy}.
  \]
  Recognize $\partial_{(b^2)}\norm{\QmbT}_{L^2}^2$, and recall
  (\ref{Prop-Qb-Mass}). Use initial condition (\ref{Defn-P-Proximity})
  and hypothesis (\ref{Hypo-Proximity-H1}).  To prove
  (\ref{Eqn-prelimEnerEst}), expand the conservation of energy, as in
  \cite[eqn (188)]{MR-UniversalityBlowupL2Critical-04}.
  Use the normalized energy (\ref{Hypo-EnergyMomentum}) to estimate
  $\lambda^2E_0$. For the terms ${\mathcal O}(\epsilon^3)$, use the
  exponential decay of $\Qm$, the Hardy-type inequalities below, and
  hypothesis (\ref{Hypo-Proximity-H1}).
\end{proof}
\begin{lemma}[Hardy-type Inequalities]
  For any $\kappa > 0$ and for all $v \in H^1(\R^2)$,
  \begin{equation}\label{Eqn-ExpDecayByGrad}
    \int_{y\in\R^2}{\abs{v(y)}^2e^{-\kappa\abs{y}}}
    \leq C(\kappa)\left(\int{\abs{\grad v(y)}^2} + \int_{\abs{y} \leq 1}{\abs{v(y)}^2e^{-\abs{y}}}\right),
  \end{equation}
  \begin{equation}\label{Eqn-L2ByGrad}
    \int_{\abs{y} \leq \kappa}{\abs{v(y)}^2} \leq C\,\kappa^2\log \kappa
    \left(\int{\abs{\grad v(y)}^2} 
      + \int_{\abs{y}\leq 1}{\abs{v(y)}^2e^{-\abs{y}}}\right).
  \end{equation}
\end{lemma}
\begin{proof} Equation (\ref{Eqn-L2ByGrad}) is proven {\cite[equation
    (4.11)]{MR-SharpLowerL2Critical-06}}, and the same techniques
  prove (\ref{Eqn-ExpDecayByGrad}).
\end{proof}

Let us reiterate the notation $y = \frac{x}{\lambda(t)} \in \R^2$, and
introduce a rescaled time,
\begin{equation}\label{Defn-Eqn-s}
  \begin{aligned}
    s(t) = \int_0^t{\frac{d\,\tau}{\lambda^2(\tau)}} + s_0 && \text{
      where } && s_0 = e^\frac{3\pi}{4b_0} && \text{ and, } && s_1 =
    s(T_{\hyp}) \in (s_0,\infty].
  \end{aligned}\end{equation}
In these new variables, equation (\ref{Eqn-NLS}) now reads,
\begin{equation}\begin{aligned}
    \label{Eqn-NLSRescaled}
    ib_s\frac{\partial}{\partial b}\QmbT +i\epsilon_s -M(\epsilon)
    +ib\Lambda\epsilon = &
    i\left(\frac{\lambda_s}{\lambda}+b\right)\Lambda\QmbT
    +\tilde{\gamma}_s\QmbT\\
    & +i\left(\frac{\lambda_s}{\lambda}+b\right)\Lambda\epsilon
    +\tilde{\gamma}_s\epsilon\\
    & +\Psi_b -R[\epsilon],
  \end{aligned}\end{equation}
where we introduced the new variable, $\tilde{\gamma}(s) = -s -
\gamma(s),$ the term $R[\epsilon]$ corresponds to those terms of
$u\abs{u}^2$ that are formally ${\mathcal O}(\epsilon^2)$, and $M$ is
the linearized operator near $\QmbT$, analogous to $\Lm$,
(\ref{Eqn-Defn-L}).  Using our choice of notation
(\ref{Eqn-Defn-Notation}), equation (\ref{Eqn-NLSRescaled}) has
exactly the same form as that given by Merle \& Rapha\"el \cite[Lemma
7]{MR-UniversalityBlowupL2Critical-04} in the case of $m=0$.
Indeed, the algebraic structure in $\Hm$ is the same, and the
arguments of \cite[Appendix C]{MR-UniversalityBlowupL2Critical-04} (or
\cite[Appendix A]{R-StabilityOfLogLog-05}) prove the following Lemma,
without modification.

\begin{lemma}
  \label{Lemma-PrelimEstimatesOrthogConds}
  For all $s\in[s_0,s_1)$,
  \begin{description}
  \item[\it due to orthogonality with $\abs{y}^2\QmbT$,
    $\Lambda\QmbT$, and estimate (\ref{Eqn-prelimEnerEst}),]
    \begin{equation}\label{Eqn-prelimLambda+BEst}
      \abs{\frac{\lambda_s}{\lambda}+b} +\abs{b_s} 
      \lesssim
      \Gamma_b^{1-C\eta} + \left(\epsNorm\right),
    \end{equation}
  \item[\it due to orthogonality with $\Lambda^2\QmbT$,]
    \begin{equation}\label{Eqn-prelimGamma+REst}
      \abs{\tilde{\gamma}_s - 
	\frac{\left\langle\epsilon_1,\Lm_+\Lambda^2\Qm\right\rangle}{\norm{\Lambda \Qm}_{L^2}^2} }
      {
        \leq
        \Gamma_b^{1-C\eta} + \delta(\alpha^*)\left(\epsNorm\right)^\frac{1}{2}.
      }
    \end{equation}
  \end{description}
\end{lemma}

In order to show the coercive control (\ref{Hypo-Proximity-H1}) does
not fail, we prove the following Local Virial Identity.  This estimate
was originally shown by Merle \& Rapha\"el in
\cite{MR-BlowupDynamic-05} and was inspired by the work of Martel \&
Merle \cite{MartelMerle-LiouvilleThmGKdV-00} in a proof of soliton
stability for the generalized Korteweg-de Vries equation.

\begin{lemma}[Local Virial Identity]\label{Lemma-LocalVirial}
  For all $s\in[s_0,s_1)$,
  \begin{equation}\label{Eqn-LocalVirial}
    b_s \geq \delta_1\left(\epsNorm\right) - \Gamma_b^{1-C\eta},
  \end{equation}
  where $\delta_1 > 0$ is a universal constant.
\end{lemma}
\begin{proof}[Proof Outline.]
  We begin the same as the proof of (\ref{Eqn-prelimLambda+BEst}):
  take the real part of the inner product of (\ref{Eqn-NLSRescaled})
  by $\Lambda\QmbT$ and use (\ref{Eqn-prelimEnerEst}) to eliminate the
  terms ${\mathcal O}(\epsilon)$, as in \cite[Appendix
  C]{MR-UniversalityBlowupL2Critical-04}.
  The interim result is,
  \begin{equation}\label{Lemma-LocalVirial-Proof}\begin{aligned}
      b_s\,\frac{1}{4}\norm{y\Qm}_{L^2}^2 \gtrsim\;
      &\calHm(\epsilon,\epsilon) - \tilde{\gamma}_s\left(\epsilon_1,\Lambda \Qm\right)\\
      &-\Gamma_b^{1-C\eta} - \delta(\alpha^*)\left(\epsNorm\right),
    \end{aligned}\end{equation}
  where we have used the preliminary estimate
  (\ref{Eqn-prelimLambda+BEst}).
  We have also used the proximity of $\QmbT$ to $\Qm$,
  (\ref{Prop-Qb-closeToQ}), to isolate the $b$-dependence from
  interactions of the form $\left(\QmbT\right)^2\epsilon^2$ as a
  lower-order potential, the same form as \cite[equation
  (215)]{MR-UniversalityBlowupL2Critical-04}, here included as part of
  the final term.  To prove the Local Virial Identity, use the
  preliminary estimate on $\tilde{\gamma}_s$ and the Spectral
  Property, Proposition \ref{Prop-SpectralProperty}, as adapted by
  Lemma \ref{Lemma-SpectralProperty2}.
\end{proof}

\subsection{Lyapounov Functional}\label{SubSec-Lyapounov}

We cannot hope to prove $b^2$ is monotonically decreasing, since it is
a modulation parameter, and thus cannot hope to control $\epsilon$ by
the local virial identity at all times.  In this section we prove a
Lyapounov functional based on the mass ejection from the singular
region, to which $b^2$ is related, (\ref{Prop-Qb-Mass}), and which we
expect $b^2$ to track. To do this, we will further approximate the
central profile by including a linear radiative tail.

\begin{prop}[Linear Radiation]\label{Prop-Zb}
  For $\eta > 0$ sufficiently small, and all $\abs{b}>0$ sufficiently
  small depending on $\eta$, there exists a unique solution $\Zmb\in
  \dotHm(\R^2)$ to
  \begin{equation}\label{Prop-Zb-Eqn}
    \laplacian\Zmb - \Zmb + ib\Lambda\Zmb = \Psi_b,
  \end{equation}
  where $\Psi_b$ is the truncation error given by
  (\ref{Prop-Qb-Eqn}). Radiation $\Zmb\not\in L^2(\R^2)$, and,
  moreover, $\lim_{\abs{y}\to+\infty}\abs{y}\abs{\Zmb(y)}^2$
  exists. We denote this decay rate as, $\Gamma_b$.
  \begin{itemize}

  \item {\it Size in $\dot{H}^1$ and Derivative by $b$}:
    \begin{equation}\label{Prop-Zb-smallH1}\begin{aligned}
        \norm{\Zmb}_{\dot{H}^1}^2 \leq \Gamma_b^{1-C\eta}, &&\text{
          and, }&& \norm{\frac{\partial}{\partial b}\Zmb}_{C^1} \leq
        \Gamma_b^{\frac{1}{2}-C\eta}.
      \end{aligned}\end{equation}

  \item {\it Decay past the support of $\Psi_b$}:
    \begin{equation}\begin{aligned}\label{Prop-Zb-H1L2Near}
        \norm{\abs{y}\abs{\Zmb}+\abs{y}^2\abs{\grad\Zmb}}_{L^\infty(\abs{y}\geq
          R_b)} &\leq \Gamma_b^{\frac{1}{2}-C\eta} < +\infty.
      \end{aligned}\end{equation}

  \item {\it Stronger decay far past the support of $\Psi_b$}:
    \[\begin{aligned}
      e^{-(1+C\eta)\frac{\pi}{b}} \leq \frac{4}{5} \Gamma_b \leq
      \norm{\abs{y}^2\abs{\Zmb}^2}_{L^\infty(\abs{y}\geq R_b^2)} \leq
      e^{-(1-C\eta)\frac{\pi}{b}},
    \end{aligned}\] which we have already discussed, equation
    (\ref{Eqn-GammaBEstimate}), and,
    \begin{equation}
      \label{Prop-Zb-H1Far}
      \begin{aligned}
        \norm{\abs{y}^2\abs{\grad\Zmb}}_{L^\infty(\abs{y}\geq R_b^2)}
        \leq C\frac{\Gamma_b^\frac{1}{2}}{\abs{b}}.
      \end{aligned}\end{equation}
  \end{itemize}
\end{prop}

The proof of Proposition \ref{Prop-Zb} is given due to Merle and
Rapha\"el \cite[Appendix E]{MR-UniversalityBlowupL2Critical-04} and
\cite[Appendix A]{MR-SharpLowerL2Critical-06}. Brief discussion of the
necessary adapatations will be given at the end of Appendix
\ref{Appendix-ProofOfAlmostSelfSimilar}.

We denote,
\begin{equation}\begin{aligned}
    \label{DefnEqn-A}
    A(t) = e^{a\frac{\pi}{b(t)}}, && \text{ so that, } &&
    \Gamma_b^{-\frac{a}{2}} \leq A \leq \Gamma_b^{-\frac{3a}{2}},
  \end{aligned}\end{equation}
where $a>0$ is a universal parameter. Let $\phi_A$ denote a smooth
cutoff function of the region, ${\mathds 1}_{\{\abs{y} > 2A\}}$.
The truncated radiation, $\ZmbT = (1-\phi_A)\Zmb$, is algebraically
close to $\Zmb$ and satisfies,
\begin{equation}\label{DefnEqn-TildeZb}\begin{aligned}
    \laplacian\ZmbT - \ZmbT +ib\Lambda\ZmbT = \Psi_b + F, &&\text{
      where, }&& \abs{F}_{L^\infty}+\abs{y\cdot\grad F}_{L^\infty}
    \lesssim \frac{\Gamma^\frac{1}{2}_b}{A}.
  \end{aligned}\end{equation}
We will now repeat the calculation of the local virial identity, this
time including the linear radiation $\widetilde{\zeta}_b$ as part of
the central profile. That is we write,
\begin{equation}\label{DefnEqn-epsTilde}\begin{aligned}
    \epsTilde = \epsilon - \ZmbT && \Rightarrow && u(t,x) =
    \frac{1}{\lambda(t)}\left(\QmbT + \ZmbT + \epsTilde\right)
    \left(\frac{x}{\lambda}\right)e^{-i\gamma(t)},
  \end{aligned}\end{equation}
without affecting the parameters. This leads to a refined version of
equation (\ref{Eqn-NLSRescaled}) for $\epsTilde$. The proof of the
following three Lemmas is virtually identical\footnotemark to that of
Merle and Rapha\"el, \cite[Chapter 4]{MR-SharpLowerL2Critical-06}.
\footnotetext{Where Merle and Rapha\"el write, $\widetilde{\zeta}_b =
  \widetilde{\zeta}_{re}+i\widetilde{\zeta}_{im}$, one should instead
  read, $\ZmbT = \widetilde{\zeta}_1 + i\widetilde{\zeta}_2$, each
  component with spin $m$ following the convention of equation
  (\ref{Eqn-Defn-Notation}).}
\begin{lemma}[Radiative Virial Identity]\label{Lemma-RefinedVirial}
  For all $s\in[s_0,s_1)$,
  \begin{equation}\label{Eqn-RadiativeVirial}\begin{aligned}
      \partial_s f_1 \geq &\delta_2\left(\epsTildeNorm\right)
      +\Gamma_b - \frac{1}{\delta_2}\int_{A\leq\abs{y}\leq
        2A}{\abs{\epsilon}^2\,dy},
    \end{aligned}\end{equation}
  where $\delta_2 > 0$ is a universal constants and,
  \[
  f_1(s) = \frac{b}{4}\norm{y\QmbT}_{L^2}^2 +
  \frac{1}{2}\Imag\left(\int{y\cdot\grad\QmbT\ZmbTbar}\right)
  +\Imag\inner{\epsilon}{\Lambda\ZmbT}.
  \]
\end{lemma}
In the light of estimates such as (\ref{Eqn-L2ByGrad}) we cannot
expect the radiative virial identity to give a good control for
$\epsilon$.  Let $\phi_\infty$ denote a smooth cutoff function of the
region ${\mathds 1}_{\{\abs{y} > 3A\}}$ with steady derivative
$\phi'_\infty \approx \frac{1}{3A}$ on the region $A\leq\abs{y}\leq
2A$.
\begin{lemma}[Mass-Ejection] 
  \label{Lemma-MassDisperse}
  \begin{equation}\label{Eqn-MassDisperse}
    \partial_s\left(\int{\phi_{\infty}\left(\frac{y}{A}\right)\abs{\epsilon}^2\,dy}\right)
    \geq \frac{b}{400}\int_{A\leq \abs{y} \leq 2A}{\abs{\epsilon}^2\,dy}
    -\Gamma_b^\frac{a}{2}\int{\abs{\grad_y\epsilon}^2\,dy}
    -\Gamma_b^2.
  \end{equation}
\end{lemma}
\begin{remark}
  As a heuristic, assume that $\epsilon \approx \zeta_b$ on the
  region, $A \leq \abs{y} \leq 2A$. Use the definition of $\Gamma_b$
  to approximate the mass. Then with hypothesis
  (\ref{Hypo-Proximity-H1}), Lemma \ref{Lemma-MassDisperse} suggests
  continuous ejection of mass from the region $\abs{y} < \frac{A}{2}$,
  regardless of whether that region is growing or contracting.
\end{remark}
Together with the conservation of mass, Lemma
\ref{Lemma-RefinedVirial} and Lemma \ref{Lemma-MassDisperse} prove the
following Lemma. The argument relies on (\ref{Prop-Zb-smallH1}) and
the relation between parameters $a$ and $\eta$ stipulated by
Proposition \ref{Prop-Qb}.
\begin{lemma}[Lyapounov Functional]\label{Lemma-LyapounovFunc}
  For all $s\in[s_0,s_1)$,
  \begin{equation}\label{Eqn-Lyapounov}
    \partial_s{\mathcal J} \leq -Cb\left(
      \Gamma_b + \epsTildeNorm + \int_{A\leq\abs{y}\leq 2A}{\abs{\epsilon}^2}
    \right),
  \end{equation}
  where $C>0$ is a universal constant,
  \begin{equation}\label{DefnEqn-LyapounovFunctional}\begin{aligned}
      {\mathcal J}(s) =
      &\norm{\QmbT}_{L^2}^2 - \norm{\Qm}_{L^2}^2\\
      &+ 2\left\langle\epsilon_1,\Sigma\right\rangle +
      2\left\langle\epsilon_2,\Theta\right\rangle
      +\int{\left(1-\phi_{\infty}\right)\abs{\epsilon}^2\,dy}\\
      &-\frac{\delta_2}{800}\left( b\widetilde{f}_1(b) -
        \int_0^b{\widetilde{f}_1(v)\,dv}
        +b\,Im\left(\epsilon,\Lambda\ZmbT\right) \right),
    \end{aligned}\end{equation}
  and $\widetilde{f}_1$ is the principal part of $f_1$,
  \[
  \widetilde{f}_1(b) = \frac{b}{4}\norm{y\QmbT}_{L^2}^2 +
  \frac{1}{2}\Imag\left(\int{y\cdot\grad\ZmbT\ZmbTbar}\right).
  \]
\end{lemma}

\subsubsection{Estimates on \texorpdfstring{$\J$}{J}}
To first order, $\J$ quantifies the excess mass remaining in the
singular region.  After explicitly accounting for this mass, $\J$ is
comparable to $\norm{\epsilon}_{H^1}^2$, up to a power of $\Gamma_b$
that depends on our choice of truncation of the radiation.

\begin{lemma}\label{Lemma-CrudeLyapounov}
  For all $s\in[s_0,s_1)$ we have the crude estimate,
  \begin{equation}\label{Eqn-CrudeLyapounovEst}
    \abs{ {\mathcal J} - d_mb^2 } < \delta_3b^2,
  \end{equation}
  where $0 < \delta_3 \ll 1$ is a universal constant, and $d_mb^2$ is
  the approximate excess mass of profile $\QmbT$.
\end{lemma}

\begin{lemma}\label{Lemma-RefinedLyapounov}
  Let $f_2$ denote those terms of ${\mathcal J}$ that are formally
  ${\mathcal O}(b^2)$,
  \[
  f_2(b) = \norm{\QmbT}_{L^2}^2 - \norm{\Qm}_{L^2}^2
  -\frac{\delta_2}{800}\left( b\tilde{f}_1(b) -
    \int_0^b{\tilde{f}_1(v)\,dv} \right).
  \]
  These are the terms concerned with the excess mass.  For all
  $s\in[s_0,s_1)$ we have the refined estimate,
  \begin{equation}\label{Eqn-RefinedLyapounovEst} {\mathcal J}(s) -
    f_2(b(s)) \left\{\begin{aligned}
	&\leq
        \Gamma_b^{1-Ca} &+& CA^2\log A\left(\epsNorm\right)\\
	&\geq
        -\Gamma_b^{1-Ca} &+& \frac{1}{C}\left(\epsNorm\right).
      \end{aligned}\right.
  \end{equation}
\end{lemma}
\begin{proof}
  The crude estimate (\ref{Eqn-CrudeLyapounovEst}) can be either
  proven directly or seen as a special case of
  (\ref{Eqn-RefinedLyapounovEst}) and hypothesis
  (\ref{Hypo-Proximity-H1}). The estimate
  (\ref{Eqn-RefinedLyapounovEst}) and its proof is exactly as given by
  Merle \& Rapha\"el, \cite[equation
  (5.6)]{MR-SharpLowerL2Critical-06}. The essential point is that the
  most difficult term of ${\mathcal J}(s) - f_2(b(s))$ can be handled
  with the conservation of energy (\ref{ConserveEnergy}), here written
  in rescaled variables,
  \[\begin{aligned}
    2\left\langle\epsilon_1,\Sigma\right\rangle +
    2\left\langle\epsilon_2,\Theta\right\rangle +
    \int{(1-\phi_A)\abs{\epsilon}^2} = \left\langle
      \Lm_+\epsilon_1,\epsilon_1\right\rangle +&\left\langle
      \Lm_-\epsilon_2,\epsilon_2\right\rangle
    -\int{\phi_A\abs{\epsilon}^2}\\
    &+ E[\QmbT] -\lambda^2E_0 + {\mathcal O}(\epsilon^3).
  \end{aligned}\]
  To establish the lower bound of (\ref{Eqn-RefinedLyapounovEst}) we
  need $\Lm$ to be coercive. We claim that,
  \begin{lemma}\label{Lemma-Maris}
    For $v = v_1 + iv_2\in \Hm$,
    \begin{equation}\label{Eqn-Coercivity-L}
      \left\langle \Lm_+v_1,v_1\right\rangle
      +\left\langle \Lm_-v_2,v_2\right\rangle
      \geq \delta_3\norm{v}_{H^1}^2 - \frac{1}{\delta_3}\left(
        \left\langle v_1,\phi_+\right\rangle^2
        +\left\langle v_1,\grad\Qm\right\rangle^2
        +\left\langle v_2,\Qm\right\rangle^2
      \right),
    \end{equation}
    where $\phi_+$ is the normalized eigenvector corresponding to the
    smallest eigenvalue of $\Lm_+$.
  \end{lemma}
  Merle and Rapha\"el \cite[Appendix D]{MR-SharpLowerL2Critical-06}
  remark that $\phi_+$ lies in the span of $\Qm$ and $\abs{y}^2\Qm$,
  and that (\ref{Eqn-Coercivity-L}) may be localized to,
  \[
  \begin{aligned}
    \left\langle \Lm_+v_1,v_1\right\rangle +&\left\langle
      \Lm_-v_2,v_2\right\rangle
    -\int{\phi_A\abs{v}^2}\\
    &\geq \delta_2\norm{v}_{H^1}^2 - \frac{1}{\delta_2}\left(
      \left\langle v_1,\Qm\right\rangle^2 +\left\langle
        v_1,\abs{y}^2\Qm\right\rangle^2
      +\left\langle v_2,\Qm\right\rangle^2 \right),
  \end{aligned}\] assuming $A$ is sufficiently large relative to the
  exponential decay of $\Qm$.
\end{proof}
\begin{proof}[Proof of Lemma \ref{Lemma-Maris}.]
  Following the variational characterization of $\Qm$, Weinstein
  \cite[Prop 2.7]{Weinstein85} argues (in the case $m=0$) that for all
  $f\in\Hm$, $\left.\partial^2_\epsilon J[\Qm+\epsilon
    f]\right|_{\epsilon=0} \geq 0$. By an explicit calculation we
  conclude,
  \[
  \inf_{f\in\Hm, \inner{f}{\Qm}=0}\inner{\Lm_+f}{f} \geq 0.
  \]
  Let $\mu_+<0$ be the lowest eigenvalue of $\Lm_+$, and $\phi_+\in
  L^2$ the corresponding normalized eigenvector.
  If there were two linearizely independent negative directions, then
  there would be one perpendicular to $\Rm$. Therefore,
  \[
  \inf_{f\in\Hm, \inner{f}{\phi_+}=0}\inner{\Lm_+f}{f} \geq 0.
  \]
  The following argument due to \cite{MartelMerle01} is an improvement
  on the proof of \cite[Prop 2.9]{Weinstein85}. Consider,
  \[
  \delta_+ = \inf\left\{\begin{aligned}\inner{\Lm_+f}{f} &&|&&
      \norm{f}_{\Hm} = 1 && \text{ and } && \inner{f}{\phi_+} =
      0\end{aligned}\right\} \geq 0.
  \]
  Assume $\delta_+ = 0$. Then by weak convergence of a minimizing
  sequence there exists a minimizer $f_+\in\Hm$, and there are
  lagrange multipliers so that,
  \[
  \left(\Lm_+-\ell_1\right)f_+ = \ell_2\phi_+.
  \]
  An inner product with $f_+$ implies $\ell_1 = 0$, and then an inner
  product by $\phi_+$ implies $\ell_2=0$. As we remarked in Section
  \ref{SubSec-VortexSolitons}, Pego \& Warchall
  \cite{PegoWarchall-VorticesNLS-02} found that the nullspace of
  $\Lm_+$ restricted to $\Hm$ is empty, and we have a contradiction.
\end{proof}

\subsection{Description of the Blowup}

Let us consider the hypotheses of Section \ref{SubSec-DecompModulate}
in turn. In each case, we will show that if the solution exists for $t
= T_{hyp}$, then the hypothesis holds for some interval
$[0,T_{\hyp}+\delta)$. This will prove that {\bf Case 1}, introduced
in Section \ref{SubSec-DecompModulate}, cannot occur, and that,
therefore, $T_{\hyp} = T_{\max}$.  This means that the dynamics
described by (\ref{Hypo-Proximity-L2}), (\ref{Hypo-Proximity-H1}),
(\ref{Hypo-loglog}) and (\ref{Hypo-EnergyMomentum}) persist for the
remaining lifetime of the solution. Indeed, we will show that that
lifetime is finite, equation (\ref{Corollary-THypSmall}), and use
these dynamics to prove the behaviour claimed for Theorem
\ref{Thm-MainResult}.

\subsubsection{Hypothesis (\ref{Hypo-Proximity-L2}).}

\noindent Preliminary estimate (\ref{Eqn-prelimMassEst}) shows that
hypothesis (\ref{Hypo-Proximity-L2}) cannot fail.

\subsubsection{Hypothesis
  (\ref{Hypo-Proximity-H1}).}\label{SubSubSec-HypoH1}

\begin{lemma} For all $s\in[s_0,s_1)$,
  \begin{equation}\label{Eqn-epsImproved}\begin{aligned}
      \epsNorm \leq \Gamma_b^\frac{4}{5},
    \end{aligned}\end{equation}
  which shows that hypothesis (\ref{Hypo-Proximity-H1}) cannot fail.
\end{lemma}
\begin{proof}
  Consider arbitrary fixed $s \in[s_0,s_1)$.
  \begin{enumerate}
  \item[\bf (a)] If $\partial_sb(s) \leq 0$, then
    (\ref{Eqn-epsImproved}) follows from the local virial identity,
    (\ref{Eqn-LocalVirial}).
  \item[\bf (b)] If $\partial_sb(s) > 0$, then there exists a largest
    interval $(s_+,s)$, with $s_0\leq s_+$, on which $\partial_sb >0$.
    \[
    \begin{aligned}
      \text{This implies, } && b(s_+) < b(s) && \text{ and either, }
      && \left.\begin{aligned}
          \left({\bf c}\right) &&&s_+ = s_0, \\
          \;\;\text{ or,}\\
          \left({\bf d}\right) &&&\partial_sb(s_+) = 0.
        \end{aligned}\right.
    \end{aligned}\] In case {\bf (c)} or {\bf (d)}, either by the
    initial condition (\ref{Defn-P-Proximity}) or the local virial
    identity, respectively,
    \[
    \int{\abs{\grad_y\epsilon(s_+,y)}^2\,dy}
    +\int_{\abs{y}\leq\frac{10}{b(s_+)}}{\abs{\epsilon(s_+,y)}^2e^{-\abs{y}}\,dy}
    \leq \Gamma_{b(s_+)}^\frac{6}{7}.
    \]
    From the upper bound of refined estimate
    (\ref{Eqn-RefinedLyapounovEst}), and assuming $a>0$ is
    sufficiently small,
    \begin{equation}\label{Proof-LowerBound-eqn5} {\mathcal J}(s_+) -
      f_2(b(s_+)) \leq \Gamma_{b(s_+)}^\frac{5}{6} <
      \Gamma_{b(s)}^\frac{5}{6}.
    \end{equation}
    Since ${\mathcal J}$ is non-increasing, 
    and from the lower bound of refined estimate
    (\ref{Eqn-RefinedLyapounovEst}),
    \begin{equation}\label{Proof-LowerBound-eqn6}\begin{aligned}
        \Gamma_{b(s)}^\frac{5}{6} \geq &
        {\mathcal J}(s) - f_2(b(s_+))\\
        \gtrsim & \left(\int{\abs{\grad_y\epsilon(s,y)}^2\,dy}
          +\int_{\abs{y}\leq\frac{10}{b(s)}}{\abs{\epsilon(s,y)}^2e^{-\abs{y}}\,dy}
	\right)\\
        &\qquad-\Gamma_{b(s)}^{1-Ca} + \left(f_2(b(s)) -
          f_2(b(s_+))\right).
      \end{aligned}\end{equation}
    As noted in the proof of crude estimate
    (\ref{Eqn-CrudeLyapounovEst}), we may assume the constant
    $\delta_2$ of equation (\ref{Eqn-RadiativeVirial}) is sufficiently
    small relative to $d_0$, such that $0 < \left.\frac{\partial
        f_2}{\partial_{b^2}}\right|_{b^2=0} < \infty$, and proving
    that $\left(f_2(b(s)) - f_2(b(s_+))\right) > 0$. Assuming $a>0$ is
    sufficiently small, this proves (\ref{Eqn-epsImproved}).
  \end{enumerate}
\end{proof}

\subsubsection{Hypothesis (\ref{Hypo-loglog}).}

\begin{lemma}
  \label{Lemma-Upperbound}
  For all $s\in[s_0,s_1)$,
  \begin{equation}\label{Eqn-bLowerBound+lambdaUpperBound}\begin{aligned}
      b(s) \geq \frac{3\pi}{4\log s}, &&\text{ and, }&& \lambda(s)
      \leq \sqrt{\lambda_0}e^{-\frac{\pi}{3}\frac{s}{\log s}}.
    \end{aligned}\end{equation}
\end{lemma}
\begin{proof}
  Recall the bounds on $\Gamma_b$, (\ref{Eqn-GammaBEstimate}),
  hypothesis (\ref{Hypo-Proximity-H1}), inject into the local virial
  identity (\ref{Eqn-LocalVirial}), carefully integrate, and recall
  the clever choice of $s_0$, (\ref{Defn-Eqn-s}),
  \[\begin{aligned}
    \partial_se^{+\frac{3\pi}{4b}} =
    -\frac{b_s}{b^2}\frac{3\pi}{4}e^{+\frac{3\pi}{4b}} \leq 1 &&
    \Longrightarrow && e^{+\frac{3\pi}{4b}} \leq s -s_0 +
    e^{+\frac{3\pi}{4b_0}} = s.
  \end{aligned}\] Next, we take hypothesis (\ref{Hypo-Proximity-H1})
  and preliminary estimate (\ref{Eqn-prelimLambda+BEst}) to
  approximate the dynamics of $\lambda$,
  \begin{equation}\label{Eqn-lambda-prelimDynamics}\begin{aligned}
      \abs{\frac{\lambda_s}{\lambda}+b}+\abs{b_s} <
      \Gamma_b^\frac{1}{2} && \Longrightarrow &&
      -\frac{\lambda_s}{\lambda} \geq \frac{2b}{3} \geq
      \frac{\pi}{2\log s}.
    \end{aligned}\end{equation}
  By the initial condition on $b_0$, (\ref{Defn-P-Proximity}), we may
  assume $s_0$ is sufficiently large so that,
  $\int_{s_0}^s{\frac{\pi}{2\log\sigma}d\sigma} \geq
  \frac{\pi}{3}\left(\frac{s}{\log s}-\frac{s_0}{\log s_0}\right).$ By
  the initial choice of a log-log relationship, (\ref{Defn-P-loglog}),
  we may assume, $-\log\lambda_0 \geq e^{\frac{\pi}{2b_0}} =
  s_0^\frac{3}{2}$. That is, by integrating
  (\ref{Eqn-lambda-prelimDynamics}) we have,
  \[
  -\log\lambda \geq -\frac{1}{2}\log\lambda_0 +
  \frac{\pi}{3}\frac{s}{\log s}.
  \]
\end{proof}
\begin{corollary}\label{Corollary-THypSmall}
  \[
  T_{\hyp} = \int_{s_0}^{s_1}{\lambda^2(\sigma)\,d\sigma} \leq
  \lambda_0\int_{2}^{+\infty}{e^{-\frac{2\pi}{3}\frac{\sigma}{\log
        \sigma}}\,d\sigma} < \alpha^*.
  \]
\end{corollary}
\begin{corollary}
  \[\begin{aligned}
    \lambda \leq e^{-e^{\frac{\pi}{5b}}}, && \text{ which shows that
      half of hypothesis (\ref{Hypo-loglog}) cannot fail.}
  \end{aligned}\]
\end{corollary}
\begin{proof}
  Due to equation (\ref{Eqn-bLowerBound+lambdaUpperBound}), and again
  assuming $s_0 > 0$ is sufficiently large,
  \[
  -\log\left(s\lambda(s)\right) \geq \frac{\pi}{3}\frac{s}{\log s} -
  \log s \geq \frac{s}{\log s}.
  \]
  Take the logarithm and apply equation
  (\ref{Eqn-bLowerBound+lambdaUpperBound}) again,
  \begin{equation}\label{Proof-Improv1-loglog-eqn2}\begin{aligned}
      \log\abs{-\log\left(s\lambda(s)\right)} \geq
      \log\left(\frac{s}{\log s}\right) \geq \frac{4}{15}\log s \geq
      \frac{\pi}{5 b(s)} && \Longrightarrow && s\lambda(s) \leq
      e^{-e^{\frac{\pi}{5b}}}.
    \end{aligned}\end{equation}
\end{proof}

\begin{lemma}
  \label{Lemma-Lowerbound}
  For all $s\in[s_0,s_1)$,
  \begin{equation}\label{Eqn-bUpperBound}
    b(s) \leq \frac{4\pi}{3\log s}.
  \end{equation}
\end{lemma}

\begin{proof}[Proof of Lemma \ref{Lemma-Lowerbound}]
  Due to the crude estimate (\ref{Eqn-CrudeLyapounovEst}) and the
  Lyapounov inequality (\ref{Eqn-Lyapounov}),
  \[\begin{aligned}
    \partial_se^{+\frac{5\pi}{4}\sqrt{\frac{d_0}{{\mathcal J}}}}
    \gtrsim \frac{b}{{\mathcal
        J}}\,\Gamma_be^{\frac{5\pi}{4}\sqrt{\frac{d_0}{{\mathcal J}}}}
    \geq 1.
  \end{aligned}\] where the final inequality is due to $\frac{5}{4} >
  1+C\eta$, the bound for $\Gamma_b$, (\ref{Eqn-GammaBEstimate}), and
  assumes $\alpha^*$ is sufficiently small. By integrating the
  inequality,
  \begin{equation}\label{Lemma-Lowerbound-ProofEqn1}
    e^{+\frac{5\pi}{4}\sqrt{\frac{d_0}{{\mathcal J}(s)}}} \geq 
    e^{+\frac{5\pi}{4}\sqrt{\frac{d_0}{{\mathcal J}(s_0)}}}
    +s -s_0.\end{equation}
  Finally, by the crude estimate (\ref{Eqn-CrudeLyapounovEst}) and the definition of $s_0$ (\ref{Defn-Eqn-s}),
  \[
  e^{+\frac{5\pi}{4}\sqrt{\frac{d_0}{{\mathcal J}(s_0)}}} >
  e^\frac{\pi}{b_0} > s_0,
  \]
  which, again with the crude estimate (\ref{Eqn-CrudeLyapounovEst}),
  proves (\ref{Eqn-bUpperBound}) from
  (\ref{Lemma-Lowerbound-ProofEqn1}).  Finally, we note here a related
  estimate that will be used in Subsection
  \ref{SubSubSec-FinalDynamic}. Divide the Lyapounov inequality
  (\ref{Eqn-Lyapounov}) by $\sqrt{\mathcal J}$, integrate in time, and
  use the crude estimate once again to get,
  \begin{equation}\label{Eqn-epsIntegral}
    \int_{s_0}^s{\left(
        \Gamma_{b(\sigma)} + \epsNorm
      \right)\,d\sigma}
    \lesssim \sqrt{{\mathcal J}(s_0)} - \sqrt{{\mathcal J}(s)}
    \lesssim b_0.
  \end{equation}

\end{proof}

\begin{corollary}
  \[\begin{aligned}
    e^{-e^{\frac{5\pi}{b}}} \leq \lambda, && \text{ which shows the
      other half of hypothesis (\ref{Hypo-loglog}) cannot fail.}
  \end{aligned}\]
\end{corollary}
\begin{proof}
  From the approximate dynamics of $\lambda$, equation
  (\ref{Eqn-lambda-prelimDynamics}),
  \[\begin{aligned}
    -\frac{\lambda_s}{\lambda} \leq 3b \leq \frac{4\pi}{\log s} &&
    \Longrightarrow && -\log \lambda(s) \leq -\log\lambda_0 +
    4\pi\int_{s_0}^s{\frac{1}{\log\sigma}\,d\sigma}
  \end{aligned}\] Bound the integral with $4\pi (s-s_0)$, and by using
  (\ref{Eqn-bUpperBound}) again, $e^{4\pi(s-s_0)} \leq
  e^{4\pi\left(e^{\frac{4\pi}{3b(s)}}-s_0\right)}$.  With the
  definition of $s_0$ (\ref{Defn-Eqn-s}) and initial condition
  (\ref{Defn-P-loglog}),
  \[
  \lambda(s) \geq \lambda_0e^{4\pi s_0}\,e^{-4\pi
    e^\frac{4\pi}{3b(s)}} > e^{-e^\frac{5\pi}{b(s)}}.
  \]
\end{proof}

\subsubsection{Hypothesis (\ref{Hypo-EnergyMomentum}).}

\noindent As another consequence of the approximate dynamics of
$\lambda$, equation (\ref{Eqn-lambda-prelimDynamics}),
\[\begin{aligned}
  \frac{d}{ds}\left(\lambda^2e^\frac{5\pi}{b}\right) =
  2\lambda^2e^\frac{5\pi}{b}\left(\frac{\lambda_s}{\lambda}+b-b-\frac{5\pi
      b_s}{2b^2}\right)
  \leq& -\lambda^2be^{5\pi}{b} < 0,\\
  &\Longrightarrow \lambda^2(t)e^{\frac{5\pi}{b(t)}} \leq
  \lambda_0^2e^\frac{5\pi}{b_0}.
\end{aligned}\] Then, due to the bounds on $\Gamma_b$,
(\ref{Eqn-GammaBEstimate}), and the initial condition
(\ref{Defn-P-EnergyMomentum}),
\[
\lambda^2(t)\abs{E_0} <
\Gamma_{b(t)}^4\,e^\frac{5\pi}{b_0}\lambda_0^2\abs{E_0}
<\Gamma_{b(t)}^4\,e^\frac{5\pi}{b_0}\Gamma_{b_0}^{10} \ll
\Gamma_{b(t)}^4.
\]
This shows that hypothesis (\ref{Hypo-EnergyMomentum}) cannot fail.

\subsubsection{Dynamics of Theorem
  \ref{Thm-MainResult}}
\label{SubSubSec-FinalDynamic}
\begin{proof}[Proof of Log-log Rate]
  By proving $T_{\hyp} = T_{\max}$, we have already shown blowup in
  finite time, due to Corollary \ref{Corollary-THypSmall}. Here we
  establish the rate.
  By direct calculation and a change of variable,
  \begin{equation}\label{Proof-MainResult-loglogEqn1}\begin{aligned}
      -\partial_t\left(\lambda^2\log\abs{\log\lambda}\right) =
      &-\frac{\lambda_s}{\lambda}\log\abs{\log\lambda} \left(2 +
        \frac{1}{\abs{\log\lambda}\log\abs{\log\lambda}}\right).
    \end{aligned}\end{equation}
  Recall the approximate dynamic $ -\frac{\lambda_s}{\lambda} \approx
  b, $ and with hypothesis (\ref{Hypo-loglog}), equation
  (\ref{Proof-MainResult-loglogEqn1}) reads,
  \[
  \frac{1}{C} \leq
  -\partial_t\left(\lambda^2\log\abs{\log\lambda}\right) \leq C.
  \]
  Integrate over $[t,T_{\max})$.
  Since $\lambda$ is very small we may estimate,
  \begin{equation}\label{Proof-MainResult-loglogEqn2}\begin{aligned}
      \frac{1}{C}\left(\frac{T_{\max}-t}{\log\abs{\log(T_{\max}-t)}}\right)^\frac{1}{2}
      \leq \lambda(t) \leq C
      \left(\frac{T_{\max}-t}{\log\abs{\log(T_{\max}-t)}}\right)^\frac{1}{2}.
    \end{aligned}\end{equation}
  Moreover, the relationship between $\lambda(t)$ and the log-log rate
  has a universal asymptotic value as $t \to T_{\max}$, see \cite[Prop
  6]{MR-SharpLowerL2Critical-06}.
\end{proof}

\begin{proof}[Proof of Singularity Description in $L^2$]
  The proof here is heavily inspired by that given by Merle and
  Rapha\"el, \cite[Section 4]{MR-ProfilesQuantization-05}.
  First, we show for any $R>0$ there exists $u^*$ such that,
  \begin{equation}\label{Proof-MainResult-L2-SpatialConvergence}\begin{aligned}
      \tilde{u}(t) \to u^* && \text{ in } && L^2_x\left(\abs{x} \geq
        R\right) && \text{ as } && t \rightarrow T_{\max}.
    \end{aligned}\end{equation}

  Second, to establish equation (\ref{Thm-MainResult-L2}), we will
  prove that both,
  \begin{equation}\label{Proof-MainResult-L2-NormConvergence}\begin{aligned}
      u^*\in L^2(\R^{2}), && \text{ and, } &&\int{\abs{u^*}^2} =
      \lim_{t\rightarrow T_{\max}}\int{\abs{\tilde{u}(t)}^2}.
    \end{aligned}\end{equation}

  Let $\epsilon_0 > 0$ be arbitrary. Choose some $T_{\max} -
  \epsilon_0 <t(R) < T_{\max}$. By hypothesis (\ref{Hypo-loglog}) we
  may assume that, $u(t) = \tilde{u}$ on $\left\{\abs{x} >
    \frac{R}{4}\right\}$ for $t\in[t(R),T_{\max})$ and by equation
  (\ref{Eqn-epsIntegral}) we may assume that, $
  \int_{t(R)}^{T_{\max}}{\int{\abs{\grad\widetilde{u}}^2\,dx}\,dt} <
  \epsilon_0.  $ For a parameter $\tau > 0$, to be fixed later, we
  denote,
  \begin{equation}\label{Proof-MainResult-L2-DefnvTau}
    v^\tau(t,x) = u(t+\tau,x) - u(t,x).
  \end{equation}
  Since $t(R) < T_{\max}$, $u(t)$ is strongly continuous in $L^2$ at
  time $t(R)$. Thus, there exists $\tau_0$ such that,
  \begin{equation}\label{Proof-MainResult-L2-vTauSmall}\begin{aligned}
      \int{\abs{v^\tau(t(R))}^2\,dx} < \epsilon_0 && \text{ for all }
      \tau \in [0,\tau_0].
    \end{aligned}\end{equation}
  Denote $\phi_R$ a smooth cutoff of the region ${\mathds
    1}_{\left\{\abs{x} > R\right\}}$.  By direct calculation,
  \begin{equation}\label{Proof-MainResult-L2-FluxEqn1}\begin{aligned}
      \frac{1}{2}\partial_t\left(\int{\phi_R\abs{v^\tau}^2}\right)
      =& Im\left(\int{\grad\phi_R\cdot\grad v^\tau \overline{v^\tau}\,dx}\right)\\
      &+ Im\left(\int{\phi_R v^\tau\left(\overline{u\abs{u}^2(t+\tau)
              - u\abs{u}^2(t)}\right)\,dx}\right).
    \end{aligned}\end{equation}
  Regarding the first RH term of (\ref{Proof-MainResult-L2-FluxEqn1}),
  use H\"older, (\ref{Proof-MainResult-L2-vTauSmall}), and the choice
  of $t(R)$,
  \[
  \int_{t(R)}^{T_{\max}}\abs{ Im\left(\int{\grad\phi_R\cdot\grad
        v^\tau \overline{v^\tau}\,dx}\right) \,dt} \leq C
  \left(\int_{t(R)}^{T_{\max}}{1^2\,dt}\right)^\frac{1}{2}\epsilon_0^\frac{1}{2}
  < C\epsilon_0.
  \]
  Regarding the second RH term of
  (\ref{Proof-MainResult-L2-FluxEqn1}), control with
  $\norm{\tilde{u}}_{L^4}^4 \leq
  \norm{\tilde{u}}_{L^2}^2\norm{\tilde{u}}_{\dot{H}^1}^2 \ll
  \norm{\tilde{u}}_{H^1}^2$, and integrate in time to get control by
  $\epsilon_0$.  We have proven that $\tilde{u}$ is Cauchy on
  $\abs{x}>R$,
  \[\begin{aligned}
    \int{\phi_R\abs{v^\tau(t)}^2\,dx} < C \epsilon_0 && \text{ for all
    } \tau \in [0,\tau_0] \text{ and } t\in[t(R),T_{\max}-\tau).
  \end{aligned}\]

  We now turn our attention to
  (\ref{Proof-MainResult-L2-NormConvergence}). The profile and
  radiation have support of radius $R(t) = A(t)\lambda(t)$, which, due
  to hypothesis (\ref{Hypo-loglog}), is going to zero with a bound,
  $R(t) \leq (T_{\max}-t)^{\frac{1}{2}-\delta}$.  From the definition
  of $A(t)$ and equation (\ref{Eqn-L2ByGrad}) we may bound
  $\int(1-\phi_{R(t)})\abs{\widetilde{u}(t)}^2$ to prove,
  \[
  \lim_{t\to T_{\max}}\int{\phi_{R(t)}\abs{u(t)}^2} =\lim_{t\to
    T_{\max}}\int{\phi_{R(t)}\abs{\widetilde{u}(t)}^2} =\lim_{t\to
    T_{\max}}\int{\abs{\widetilde{u}(t)}^2},
  \]
  which proves that the following limit exists,
  \[
  \int{\abs{u^*}^2} = \lim_{t\to
    T_{\max}}\int{\phi_{R(t)}\abs{u(t)}^2}.
  \]
  This completes the proof of
  (\ref{Proof-MainResult-L2-NormConvergence}) and
  (\ref{Thm-MainResult-L2}).
\end{proof}

\section{The Spectral Property}
\label{Section-SpectralProperty}

We now provide a numerically assisted proof of the spectral property
for the case $m=1$.  We also present some computations on higher order
vortices and discuss why they do not work.  Before proving the
Spectral Property of Section \ref{s:intro_specprop}, we will establish
the following variant:
\begin{prop}
  \label{p:specprop}
  Let $\eps \in H^1_\m$ satisfy the orthogonality conditions,
  \begin{equation}
    \label{e:general_ortho_conds}
    \inner{\eps_1}{Q^\m} = \inner{\eps_1}{\Lambda Q^\m} =
    \inner{\eps_2}{\Lambda Q^\m} = \inner{\eps_2}{\Lambda^2 Q^\m} = 0.
  \end{equation}
  Then, for the case $m=1$, there is a universal constant $C_m > 0$,
  so that,
  \begin{equation}
    \calHm(\eps,\eps) \geq  C_m \int \paren{\abs{\nabla_y \eps}^2 + e^{-\abs{y}}\abs{\eps}^2}dy.
  \end{equation}
\end{prop}
Proposition \ref{Prop-SpectralProperty} is an immediate
corollary\footnotemark.  \footnotetext{See the end of Subsection
  \ref{Subsection-FinalProofSpectral} for details.}  Following
\cite{FMR-ProofOfSpectralProperty-06, Marzuola2010}, we proceed in two
steps.  First we count the number of negative eigenvalues of the
operators $\calL^\m_1$ and $\calL^\m_2$.  We then show that the
assumed $L^2$ orthogonality conditions are sufficient to project away
from the negative directions of the bilinear forms, $\calH^\m_1$ and
$\calH^\m_2$, associated with $\calL^\m_1$ and $\calL^\m_2$.

We now restrict ourselves to $\eps \in H^1_\m$, $\eps = e^{i m \theta}
\eps_{\rad}$, where $\eps_{\rad}$ is a purely radial function,
\begin{equation}
  \eps_{\rad} \in H^1_{\rad+} \equiv H_\rad^1(\R^2) \cap \set{u  \mid
    \abs{x}^{-1} u \in L^2(\R^2)}.
\end{equation}
Given $\eps \in H^1_\m$, we calculate
\begin{equation}
  \calL_1^\m \eps = e^{i m \theta} \paren{-\frac{d^2}{dr^2 } -
    \frac{1}{r}\frac{d}{dr} + \frac{m^2}{r^2} + 3 R^\m y \cdot \nabla R^\pm}\eps_\rad
\end{equation}
This motivates defining the two operators and inner products on
$H^1_{\rad+}$
\begin{subequations}
  \begin{gather}
    \calL_{1,\rad}^{(m)} \equiv-\Delta_\rad + \frac{m^2}{r^2} + 3 R^\m
    y\cdot \nabla R^\m =
    -\Delta_\rad + \frac{m^2}{r^2} + \calV_{1,\rad}\\
    \calL_{2,\rad}^{(m)} \equiv-\Delta_\rad + \frac{m^2}{r^2} + R^\m
    y\cdot \nabla R^\m =
    -\Delta_\rad + \frac{m^2}{r^2} + \calV_{2,\rad}\\
    \calH_{1,\rad}^{(m)}(\cdot, \cdot) \equiv
    \inner{\calL_{1,\rad}^{(m)}\cdot }{\cdot},\quad
    \calH_{2,\rad}^{(m)}(\cdot, \cdot) \equiv
    \inner{\calL_{2,\rad}^{(m)}\cdot }{\cdot},
  \end{gather}
\end{subequations}
where $\Delta_\rad$ is the radial Laplacian, $\Delta_\rad \equiv
\frac{d^2}{dr^2} + \frac{1}{r}\frac{d}{dr}$.

The orthogonality conditions, \eqref{e:general_ortho_conds}, are now
formulated as
\begin{equation}
  \label{e:radial_ortho_conds}
  \inner{\eps^\rad_1}{R^\m} = \inner{\eps^\rad_1}{\Lambda R^\m} =
  \inner{\eps^\rad_2}{\Lambda R^\m} = \inner{\eps^\rad_2}{\Lambda^2 R^\m} = 0,
\end{equation}
where
\[
\eps = e^{i m \theta} \paren{\eps_1^\rad + i \eps_2^\rad}
\]
and $\eps \in H^1_\m$, $\eps_j^\rad \in H^1_{\rad+}$.

All that follows relies on the reduction to a series of one
dimensional radial problems.

\subsection{The Index of Bilinear Forms}

\begin{defn}
  The {\it index} of a bilinear form $\calB$ with respect to vector
  space $V$, denoted $\ind_V \calB$, is the minimal co-dimension over
  all subspaces of $V$ on which $\calB$ is a positive.
\end{defn}
For bilinear forms induced by self-adjoint operators ({\it i.e.}
$\calB = \inner{\calL \cdot }{\cdot}$), the index corresponds to the
number of negative eigenvalues of the operator.  To calculate the
index, we extend Theorem XIII.8 of Reed \& Simon \cite{RSv4} to:
\begin{theorem}
  \label{thm:rs_idx}
  Let $U$ solve,
  \begin{equation*}
    \calL\, U = - \frac{d^2}{dr^2} U - \frac{1}{r} \ddr U  + \calV(r) U  + \frac{m^2}{r^2} U = 0,
  \end{equation*}
  with initial conditions given by the limits,
  \begin{equation*}
    \lim_{r\to 0} {r^{-m}}{U (r) }  = 1, \quad\lim_{r\to 0} \ddr\paren{
      {r^{-m}}{U (r) }}  = 0,
  \end{equation*}
  and where the potential $\calV$ is sufficiently smooth and decaying
  at $\infty$.  Then, the number of roots of $U$ not at the origin,
  $N(U)$, is finite and equal to the index of the bilinear form
  associated to $\calL $ over the vector space $H^1_{\rad+}$.
\end{theorem}
\begin{proof}
  The proof, which we omit, is quite similar to the proof of the
  indicated Theorem of Reed \& Simon.  In turn, that proof is a
  generalization of the Sturm Oscillation theorem for two point
  boundary value problems.
\end{proof}

\begin{prop}[Numerically Verified]
  \label{p:idx_computations}
  For the cases $m=1,2,3$,
  \begin{equation}\begin{aligned}
      \ind_{H^1_{\rad+}}\, \calH_{1,\rad}^\m = 2,&& \text{and}, &&
      \ind_{H^1_{\rad+}}\, \calH_{2,\rad}^\m = 1.
    \end{aligned}\end{equation}
\end{prop}
\begin{proof}
  Using the methods described in Appendix \ref{s:numerics}, we solve,
  \begin{subequations}
    \label{e:idx_m0}
    \begin{gather}
      \calL_{1,\rad}^\m U^\m = 0, \quad \lim_{r\to 0} r^{-m}U^\m(0)=0,
      \quad \lim_{r\to 0} \ddr
      {r^{-m}}{U^\m (r) }  = 0,\\
      \calL_{2,\rad}^\m Z^\m = 0, \quad \lim_{r\to 0} r^{-m}Z^\m(0)=0,
      \quad \lim_{r\to 0} \ddr {r^{-m}}{Z^\m (r) } = 0.
    \end{gather}
  \end{subequations}
  Plotting the solutions in Figure \ref{f:idx}, we can see that $U^\m$
  has two zero crossings and $Z^\m$ has one zero crossing.  Subject to
  the acceptance of these computations, Theorem \ref{thm:rs_idx}
  yields the result.
\end{proof}

\ifpdf
\begin{figure}
  \centering
  \subfigure{\includegraphics[width=2.35in]{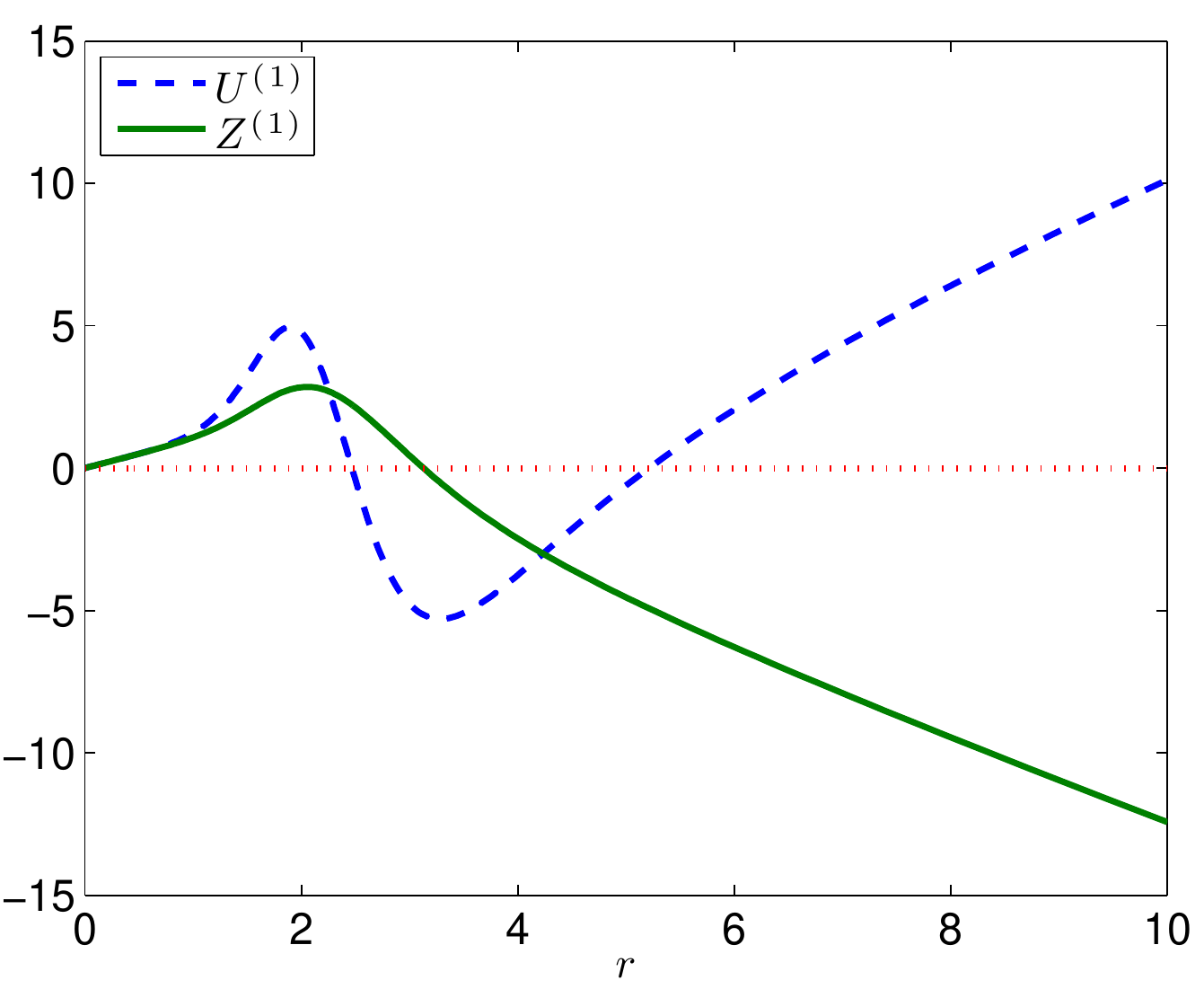}}
  \subfigure{\includegraphics[width=2.35in]{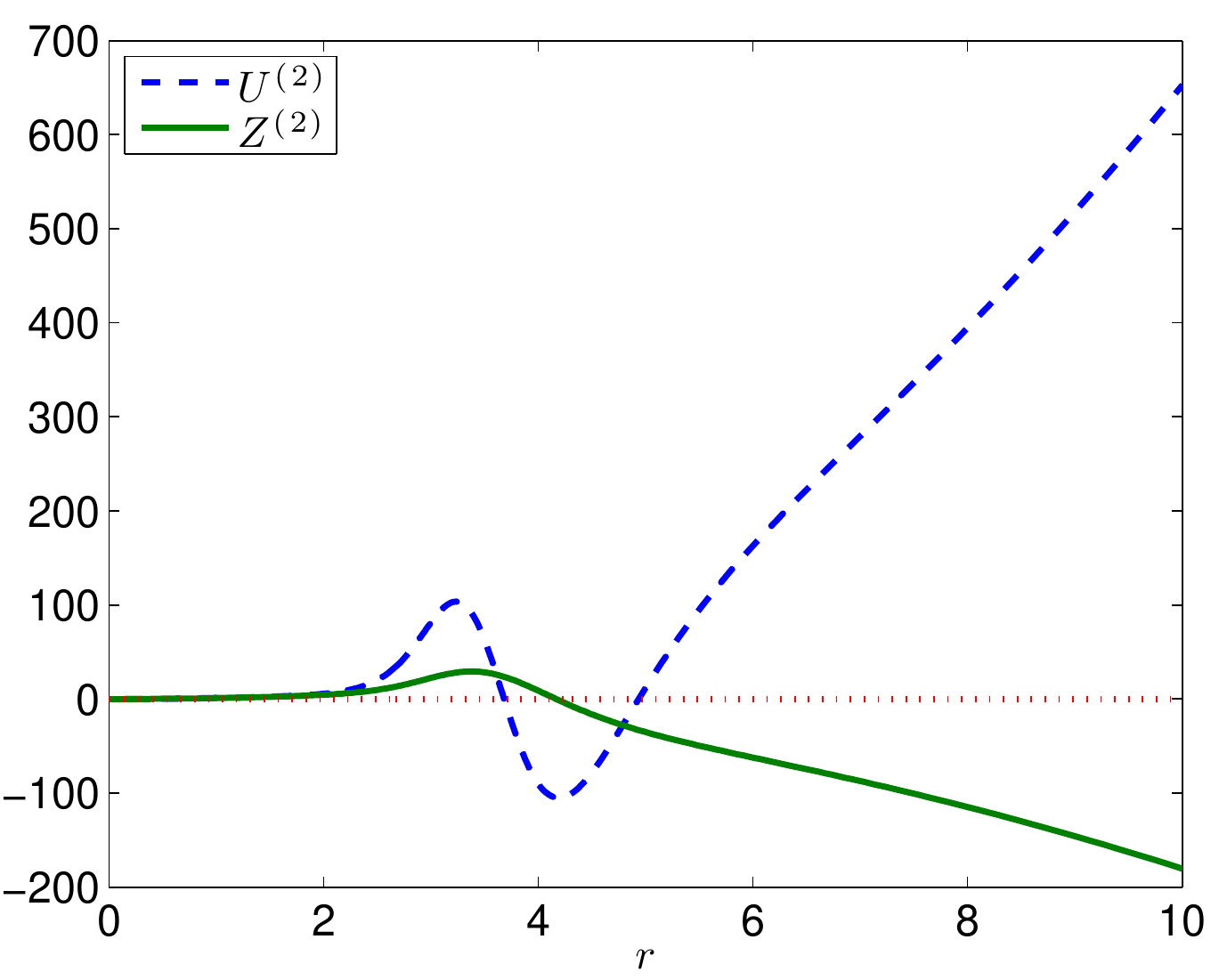}}
	
  \subfigure{\includegraphics[width=2.35in]{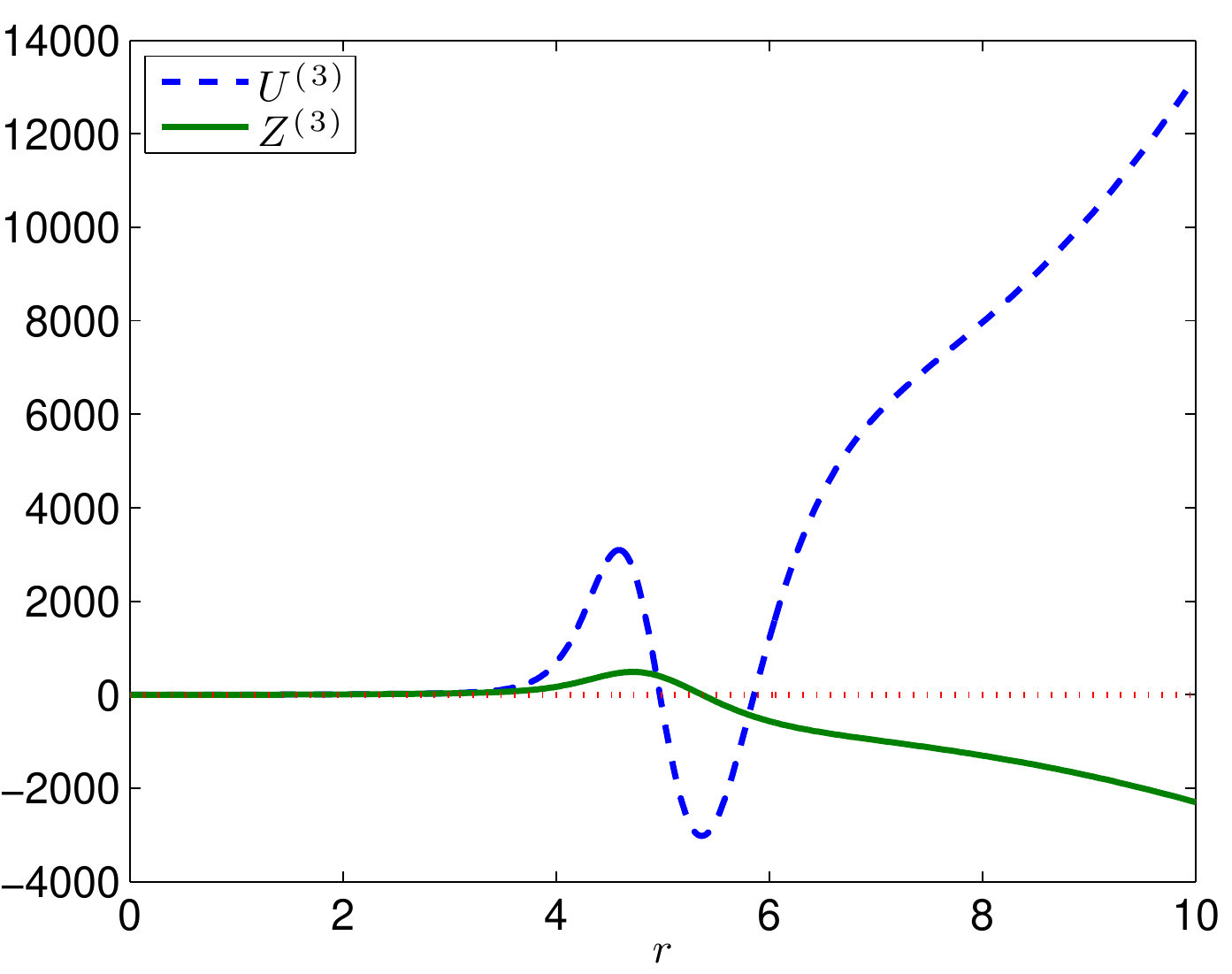}}
  \caption{Plots of the functions $U^{(m)}$ and $Z^{(m)}$ solving
    \eqref{e:idx_m0} for $m=1,2,3$.  The number of non-zero zero
    crossings determines the indexes of $H_1$ and $H_2$.}
  \label{f:idx}
\end{figure}
\else (No figures in DVI) \fi

\begin{prop}
  \label{prop:idx_perturbation}
  There exists a constant, $\delta_0>0$, depending on $m$, such that
  for $\delta\in (0, \delta_0)$, the bilinear
  forms 
  associated with the perturbed operators,
  \[
  \overline{\calL}_{j,\rad}^{(m)} \equiv {\calL}_{j,\rad}^{(m)}-
  \delta e^{-\abs{y}},
  \]
  have the same index, {\it i.e.}
  \[
  \ind \, \calH_{j,\rad}^\m =\ind \, \overline{\calH}_{j,\rad}^\m.
  \]
\end{prop}
\begin{proof}
  We briefly sketch the proof, which follows from three observations.
  First, the of the solutions of the perturbed form of
  \eqref{e:idx_m0} are continuous with respect to $\delta$. In
  particular, there is $C^1_{\loc}$ convergence. Second, the roots of
  the index functions, in the perturbed and unperturbed cases, must be
  simple.  For a sufficiently small $\delta_0$, we can ensure that on
  any compact interval the perturbed and unperturbed solutions have
  the same number of zeros.  Finally, for a sufficiently large compact
  interval, outside the interval the equation is approximately
  ``free'' (the localized potentials are negligible), and we can
  ensure there are no additional zeros; this may require further
  shrinking $\delta_0$.
\end{proof}

\subsection{Orthogonality Conditions and Inner Products}

To verify that orthogonality conditions \eqref{e:radial_ortho_conds}
project away from the negative subspaces, we need to compute a number
of inner products of the form $\inner{\calL_{j,\rad}^{(m)} u}{u}$,
where $u$ solves $\calL_{j,\rad}^{(m)} u = f$.  Although these
products are computed numerically, we justify their existence:

\begin{prop}[Numerically Verified]
  \label{prop:eu_bvp_2d}
  Let $f$ be a continuous, radially symmetric, localized function
  satisfying the bound $\abs{f(r)} \leq C e^{-\kappa r}$ for some
  positive constants $C$ and $\kappa$.  There exists a unique radially
  symmetric solution,
  \[\begin{aligned}
    \overline{\calL}_{j,\rad}^{(m)} u = f, &&j=1,2.
  \end{aligned}\] that belongs to the class, $u \in
  L^\infty([0,\infty))\cap C^2([0,\infty))$.

\end{prop}
\begin{proof}
  This is Proposition 2 and 4 of
  \cite{FMR-ProofOfSpectralProperty-06}, along with our computations
  of the indexes in Lemma \ref{e:idx_m0}.  See \cite{Marzuola2010} for
  some additional details and a full proof in dimension $d=1$.
\end{proof}
\begin{remark}
  The solutions in Proposition \ref{prop:eu_bvp_2d} may not vanish as
  $r\to \infty$.  Indeed, they can only be expected to be bounded.
\end{remark}

\begin{prop}[Numerically Verified]
  \label{prop:ip_computations}
  Let $U_1$, $U_2$, $Z_1$, and $Z_2$ be $L^\infty$ radially symmetric
  functions solving,
  \begin{subequations}
    \label{e:ip_bvps}
    \begin{align}
      \calL_{1,\rad}^{(m)} U_1 &= R^\m,\\
      \calL_{1,\rad}^{(m)} U_2 &= \Lambda R^\m,\\
      \calL_{2,\rad}^{(m)} Z_1 &= \Lambda R^\m,\\
      \calL_{2,\rad}^{(m)} Z_2 &= \Lambda^2 R^\m.
    \end{align}
  \end{subequations}
  Then the inner products,
  \[\begin{aligned}
    &K_1^{(m)}\equiv \inner{\calL_{1,\rad}^{(m)} U_1}{U_1}, &&
    &K_2^{(m)}\equiv\inner{\calL_{1,\rad}^{(m)} U_2}{U_2},
    &&K_3^{(m)}\equiv\inner{\calL_{1,\rad}^{(m)} U_1}{U_2},\\
    &J_1^{(m)}\equiv\inner{\calL_{2,\rad}^{(m)} Z_1}{Z_1}, &&
    &J_2^{(m)}\equiv\inner{\calL_{2,\rad}^{(m)} Z_2}{Z_2},
    &&J_3^{(m)}\equiv\inner{\calL_{2,\rad}^{(m)} Z_1}{Z_2},
  \end{aligned}\] take the values given in Tables \ref{t:k_ips} and
  \ref{t:j_ips}.
\end{prop}
\begin{proof}
  Using the methods described in Appendix \ref{s:numerics}, these are
  computed numerically.
\end{proof}

\begin{table}
  \centering
  \caption{Inner products associated with $\calL^\m_{1,\rad}$ for
    different winding numbers.}
  \label{t:k_ips}
  \begin{tabular}{l l l l l }
    \hline
    $m$ & $K_1^{(m)}$ & $K_2^{(m)}$ & $K_3^{(m)}$ & $K_1^{(m)} K_2^{(m)} - \left(K_3^{(m)}\right)^2$\\
    \hline
    1 &-0.48237 &-25.798 & 1.28129&10.8025\\
    2 & 0.520152&-13.1545 &1.7983 &  -10.0762\\
    3 &2.59249 & 5.1232&-1.54694 & 10.8888\\
    \hline
  \end{tabular}
\end{table}

\begin{table}
  \centering
  \caption{Inner products associated with $\calL^\m_{2,\rad}$ for
    different winding numbers.}
  \label{t:j_ips}
  \begin{tabular}{l l l l l }
    \hline
    $m$ & $J_1^{(m)}$ & $J_2^{(m)}$ & $J_3^{(m)}$ & $J_1^{(m)} J_2^{(m)} - \left(J_3^{(m)}\right)^2$\\
    \hline
    1 &6.6985 &163.548 & -47.7764&-1.1871e+03\\
    2 & 25.1685&1319.28 &-235.186 &-2.2108e+04\\
    3 &82.6396& 8426.22&-936.752 & -1.8116e+05 \\
    \hline
  \end{tabular}
\end{table}

As with the indices, we have stability of the inner products with
respect to perturbation by a small portential:
\begin{prop}
  \label{prop:ip_perturbation}
  Let $\overline{U}_l$ and $\overline{Z}_l$ denote the solutions and
  $\overline{K}_l^\m$ and $\overline{J}_l^\m$ the inner products,
  analogous to those of Proposition \ref{prop:ip_computations}, for
  the boundary value problems with the perturbed operators,
  $\overline{\calL}_{j,\rad}^\m$.  For $\delta_0>0$ sufficiently
  small, the solutions and inner products are continuous with respect
  to $\delta$.
\end{prop}
\begin{proof}
  This follows from the invertiblity and continuity with respect to
  $\delta$ of the operators.
\end{proof}

\subsection{Proof of the Spectral Property}
\label{Subsection-FinalProofSpectral}

We are now able to prove Proposition \ref{p:specprop}.  The arguement
closely follows the proofs found in
\cite{FMR-ProofOfSpectralProperty-06,Marzuola2010}. The two bilinear
forms, $\overline{\calH}_1^\m$ and $\overline{\calH}_2^\m$, are
treated seperately. First, we will show that $L^2$ orthogonality to
$Q^\m$ and $\Lambda Q^\m$ suffices to project away from the negative
subspace of $\calH_1^\m$. This will only be successful for
$m=1$. Later, we will show that $L^2$ orthogonality to $\Lambda Q^\m$
and $\Lambda^2 Q^\m$ projects away from the negative subspace of
$\calH_2^\m$.

\begin{proof}[Spectral Property for $\calH_1^\m$.]
  Given an element $u\in H^1_\m$, $u = e^{i m \theta} u_\rad$,
  satisfying orthogonality conditions \eqref{e:general_ortho_conds},
  showing positivity of $\calH_1^\m$ on such a $u$ is equivalent to
  showing posiviity of $\calH_{1,\rad}^\m$ on $u_\rad \in H^1_{\rad+}$
  satisfying orthogonality conditions \eqref{e:radial_ortho_conds}.

  By Propositions \ref{p:idx_computations} and
  \ref{prop:idx_perturbation}, $\overline{\calH}_{1,\rad}^\m$ has a
  two-dimensional subspace of negative directions.  Recall the
  notation of equation \eqref{e:ip_bvps}. Let $\mathrm{V} = \spn
  \set{\overline{U}_1, \overline{U}_2}$.  We will prove that, for
  $m=1$, $\overline{\calH}^\m_{1,\rad}$ is negative on all of
  $\mathrm{V}$.  Indeed, consider an arbitrary element of this space,
  \[
  \hat{U} = c_1 \overline{U}_1 + c_2 \overline{U}_2,
  \]
  and compute,
  \begin{equation}\label{e:matrixH_on_V}
    \begin{split}
      \overline{\calH}^\m_{1,\rad}(\hat{U},\hat{U})
      & = c_1^2 \overline{K}_1^{(m)} + 2 c_1 c_2 \overline{K}_3^{(m)}
      + c_2^2
      \overline{K}_2^{(m)} \\
      & = \begin{pmatrix} c_1 &
        c_2\end{pmatrix} \begin{pmatrix}\overline{K}_1^{(m)}
        &\overline{K}_3^{(m)}
        \\ \overline{K}_3^{(m)}  &   \overline{K}_2^{(m)} \end{pmatrix}\begin{pmatrix} c_1 \\
        c_2\end{pmatrix}.
    \end{split}
  \end{equation}
  If the above matrix is negative definite, then the bilinear form is
  negative on the two dimensional space $\mathrm{V}$.  We examine the
  matrix using the computations in Table \ref{t:k_ips} and elementary
  properties of matrices.
  For $m=1$,
  \[
  \mathrm{tr} = \overline{K}_1^{(1)} + \overline{K}_2^{(1)} =-26.2804
  + \littleo(1),
  \]
  where $\littleo(1)$ corresponds to taking the perturbation
  parameter, $\delta$, sufficiently small.  Therefore the sum of the
  two eigenvalues is negative; at least one is negative.  Next,
  \[
  \mathrm{det}= \overline{K}_1^{(1)}\overline{K}_2^{(1)} -
  (\overline{K}_3^{(1)})^2 = 10.8025 + \littleo(1),
  \]
  so the two eigenvalues have the same sign.  Therefore
  $\overline{\calH}^{(1)}_{1,\rad}$ is negative on $\mathrm{V}$.
  Table \ref{t:k_ips} shows that this is false for $m=2,3$.  We
  restrict our attention to $m=1$.


  Pretending that $\mathrm{V}\subset H^1_{\rad+}(\R^2)$, we could
  decompose the space as
  \begin{equation}\label{e:pretendDecomp}
    H^1_{\rad+}(\R^2) = \mathrm{V} \oplus_{\overline{\calH}^{(1)}_{1,\rad}} \mathrm{V}^\perp
  \end{equation}
  where our notation indicates that we have formed the orthogonal
  complement with respect to the $\overline{\calH}^{(1)}_{1,\rad}$
  bilinear form.  The non-degeneracy of the matrix
  \eqref{e:matrixH_on_V} justifies this decomposition.

  It follows that $\overline{\calH}^{(1)}_{1,\rad}$ is positive on
  $\mathrm{V}^\perp$.  Otherwise, there would be $W \in
  \mathrm{V}^\perp$ such that
  $\overline{\calH}^{(1)}_{1,\rad}(W,W)<0$, which implies by
  construction that, $\spn\set{W, \overline{U}_1, \overline{U}_2}$, is
  a negative definite space of $\overline{\calH}^{(1)}_{1,\rad}$ with
  dimension three.  But then, given any subspace $\mathrm{U} \subset
  H^1_{\rad+}$ of codimension two, $\mathrm{U} \,\cap\, \spn\set{W,
    \overline{U}_1, \overline{U}_2} \neq \emptyset$, which contradicts
  the index calculation.

  Finally, given any function $u \in H^1_{\rad+}$ and $L^2$ orthogonal
  to $R^{(1)}$ and $\Lambda R^{(1)}$, we decompose $u$ as
  \[
  u = c_1 \overline{U}_1 + c_2 \overline{U}_2 + u^\perp
  \]
  where, $u^\perp \in \mathrm{V}^\perp$, again in the sense of
  \eqref{e:pretendDecomp}.  Then,
  \[\begin{aligned}
    0 = \inner{u}{R^{(1)}}_{L^2} &= c_1
    \inner{\overline{U}_1}{R^{(1)}}_{L^2} +
    c_2 \inner{\overline{U}_2}{R^{(1)}}_{L^2} + \inner{u^\perp}{R^{(1)}}_{L^2}\\
    &= c_1 \overline{\calH}^{(1)}_{1,\rad}(\overline{U}_1,
    \overline{U}_1) +
    c_2 \overline{\calH}^{(1)}_{1,\rad}(\overline{U}_2, \overline{U}_1) + \overline{\calH}^{(1)}_{1,\rad} (u^\perp, \overline{U}_1)\\
    &= c_1 \overline{K}^{(1)}_1 + c_2 \overline{K}^{(1)}_3,\\
    0 = \inner{u}{\Lambda R^{(1)}}_{L^2} &= c_1
    \inner{\overline{U}_1}{\Lambda R^\m}_{L^2} +
    c_2 \inner{\overline{U}_2}{\Lambda R^{(1)}}_{L^2} + \inner{u^\perp}{\Lambda R^{(1)}}_{L^2}\\
    &= c_1 \overline{\calH}^{(1)}_{1,\rad}(\overline{U}_1,
    \overline{U}_2) +
    c_2 \overline{\calH}^{(1)}_{1,\rad}(\overline{U}_2, \overline{U}_2) + \overline{\calH}^{(1)}_{1,\rad} (u^\perp, \overline{U}_2)\\
    &= c_1 \overline{K}^{(1)}_3 + c_2 \overline{K}^{(1)}_2.
  \end{aligned}\] Due to the non-degeneracy of \eqref{e:matrixH_on_V},
  the only solution is $c_1=c_2=0$.  Therefore, for all such $u$,
  \begin{equation*}
    \overline{\calH}^{(1)}_{1,\rad}(u,u) \geq 0.
  \end{equation*}
  This yields the positivity of $\overline{\calH}^{(1)}_{1}$ on
  $H^1_{(1)}$.

  Of course, $\overline{U}_1$ and $\overline{U}_2$ are {\it not} in
  $H^1_{\rad+}$.  The above argument is made rigorous by introducing
  an appropriate cutoff function and then taking limits.  We refer the
  reader to \cite{FMR-ProofOfSpectralProperty-06,Marzuola2010}; we
  will not reproduce this here.

\end{proof}
\begin{proof}[Spectral Property for $\calH_2^\m$.]
  
  As in the case of $\overline{\calH}_{1}^\m$, we will prove
  positivity of $\overline{\calH}_{2}^\m$ subject to the orthogonality
  conditions, by working with the associated radial form,
  $\overline{\calH}_{2,\rad}^\m$.  By Propositions
  \ref{p:idx_computations} and \ref{prop:idx_perturbation},
  $\overline{\calH}^\m_{2,\rad}$ has one negative direction.  Examing
  Table \ref{t:j_ips}, neither $\overline{Z}_1^{(m)}$ nor
  $\overline{Z}_2^{(m)}$ appears to point in the negative direction.
  Define,
  \begin{subequations}
    \begin{align}
      \hat{R}^{(m)} &\equiv \Lambda R^\m - \frac{\overline{J}_3^{(m)}}{\overline{J}_2^{(m)}} \Lambda^2 R^\m,\\
      \hat{Z} &\equiv
      \overline{Z}_1-\frac{\overline{J}_3^{(m)}}{\overline{J}_2^{(m)}}\overline{Z}_2.
    \end{align}
  \end{subequations}
  Then $\calL_2^\m \hat{Z} = \hat{R}^{(m)} $ and,
  \begin{equation}\label{e:H2_ZhatZhat}
    \overline{\calH}^\m_2(\hat{Z} , \hat{Z} ) = \frac{1}{\overline{J}_2^{(m)}}\paren{\overline{J}_1^{(m)}\overline{J}_2^{(m)}-\left(\overline{J}_3^{(m)}\right)^2}<0.
  \end{equation}
  Now that we have constructed a negative direction, we apply a
  similar argument as in the case of $\overline{\calH}^\m_{1,\rad}$;
  however, this will hold not just for $m=1$, but also for $m=2,3$.
  We decompose $H_{\rad+}^1(\R^2)$ as
  \begin{equation}
    H_{\rad+}^1(\R^2) = \spn\set{\hat{Z} }\oplus_{\overline{\calH}^\m_{2,\rad}} \spn\set{\hat{Z} }^\perp
  \end{equation}
  Since the index of $\overline{\calH}^\m_{2,\rad}$ is one, we are
  assured that it is positive on $\spn\set{\hat{Z} }^\perp$.  Finally,
  given $v\in H_{\rad+}^1$ orthogonal to $\Lambda R^\m$ and $\Lambda^2
  R^\m$, it may be decomposed as $v = c_1 \hat{Z} + v^\perp$, and,
  \[\begin{aligned}
    0 = \inner{v}{\hat{R}^{(m)} }_{L^2}
    &= c_1 \overline{\calH}^\m_2(\hat{Z} ,\hat{Z} ) + \overline{\calH}^\m_2(v^\perp, \hat{Z} )\\
    &= c_1 \overline{\calH}^\m_2(\hat{Z} ,\hat{Z} ).
  \end{aligned}\] Invoking \eqref{e:H2_ZhatZhat}, this implies that,
  $v = v^\perp \in \spn\set{\hat{Z} }^\perp$.  Therefore, for such
  $v$,
  \[
  \overline{\calH}^\m_{2,\rad}(v,v) \geq 0
  \]
  for $m=1,2,3$.  Posivitiy of $\overline{\calH}^\m_2$ on $H^1_\m$,
  subject to orthgonality to $\Lambda Q^\m$ and $\Lambda^2 Q^\m$,
  follows.
\end{proof}

\begin{proof}[Proof of Proposition \ref{p:specprop}.]
  Given $\eps = \eps_1 + i \eps_2$ satisfying the orthogonality
  conditions of Proposition \ref{p:specprop} we have proven that,
  \begin{equation*}
    \overline{\calH}^{(1)}(\eps, \eps)  = \overline{\calH}^{(1)}_1 (\eps_1, \eps_1) + \overline{\calH}^{(1)}_2 (\eps_2,
    \eps_2) \geq 0,
  \end{equation*}
  from which we infer,
  \begin{equation*}
    \calH^{(1)} (\eps, \eps)\geq \delta \int e^{-\abs{y}} \abs{\eps}^2 dy.
  \end{equation*}
  Let $\theta >0$.  Then,
  \begin{equation*}
    (1+\theta) \calH^{(1)} (\eps, \eps) \geq \theta \int \abs{\nabla \eps}^2 dy+ \theta
    \int \calV_1 \abs{\eps_1}^2 + \calV_2 \abs{\eps_2}^2 dy + \delta \int e^{-\abs{y}} \abs{\eps}^2 dy.
  \end{equation*}
  Although the potentials are sign indefinite, for $\theta$
  sufficiently small,
  \begin{equation}
    \theta
    \int \calV_1 \abs{\eps_1}^2 + \calV_2 \abs{\eps_2}^2 dy + \delta \int
    e^{-\abs{y}} \abs{\eps}^2 dy \geq \frac{\delta}{2}\int e^{-\abs{y}} \abs{\eps}^2 dy.
  \end{equation}
  We now have the result,
  \begin{equation*}\begin{aligned}
      \calH^{(1)} (\eps, \eps) &\geq \frac{\theta}{1+\theta} \int
      \abs{\nabla \eps}^2dy +
      \frac{\delta}{2(1+\theta)}\int e^{-\abs{y}} \abs{\eps}^2 dy\\
      &\geq \delta_0 \int \abs{\nabla \eps}^2 + e^{-\abs{y}}
      \abs{\eps}^2 dy.
    \end{aligned}\end{equation*}
\end{proof}

\begin{proof}[Proof of Proposition \ref{Prop-SpectralProperty}.]
  Let $\eps\in H^1_{(1)}(\R^2)$ with $\eps = \eps_1 + i \eps_2$, and
  further decompose this as:
  \begin{subequations}
    \begin{align}
      \eps_1 &= e^{i\theta} \paren{u + c_1 R^{(1)} + c_2\Lambda R^{(1)}},\\
      \eps_2 &= e^{i\theta} \paren{v + d_1 \Lambda R^{(1)} +
        d_2\Lambda^2 R^{(1)}},
    \end{align}
  \end{subequations}
  where $u\perp_{L^2} R^{(1)}, \Lambda R^{(1)}$ and $v\perp_{L^2}
  \Lambda R^{(1)}, \Lambda^2 R^{(1)}$.  Expaning,
  \[\begin{aligned}
    &\calH^{(m)}(\eps, \eps)
    &=&&& \calH^{(m)}_1(\eps_1, \eps_1) + \calH^{(m)}_2(\eps_2, \eps_2),\\
    &\calH^{(m)}_1(\eps_1, \eps_1) & =&&& \calH^{(m)}_{1,\rad}(u,u) +
    2 c_1 \inner{\calL^{(m)}_{1,\rad} u}{R^{(1)}}
    + 2 c_2 \inner{\calL^{(m)}_{1,\rad} u}{\Lambda R^{(m)}}\\
    &&&&& + c_1^2 M^{(m)}_1 + c_2^2 M^{(m)}_2 + 2 c_1 c_2 M^{(m)}_3,\\
    &\calH^{(m)}_2(\eps_2, \eps_2) & =&&& \calH^{(m)}_{2,\rad}(v,v) +
    2 d_1 \inner{\calL^{(m)}_{2,\rad} v}{\Lambda R^{(m)}}
    + 2 d_2 \inner{\calL^{(m)}_2 v}{\Lambda^2 R^{(m)}}\\
    &&&&& + d_1^2 N^{(m)}_1 + d_2^2 N^{(m)}_2 + 2 d_1 d_2 N^{(m)}_3,
  \end{aligned}\] where $M^{(m)}_j, N^{(m)}_j$ are fixed terms arising
  from applications of the $\mathcal{H}_{j,\rad}^{(m)}$ bilinear forms
  to combinations of $R^{(m)}$, $\Lambda R^{(m)}$, and $\Lambda^2
  R^{(m)}$.

  We now construct a lower bound.  Let $\theta>0$.  Then
  \begin{equation}
    \begin{split}
      c_1 \inner{\calL^{(m)}_1
        u}{R^{(m)}}&\leq\frac{1}{2} \paren{\theta^{-2}c_1^2 +
        \theta^2\inner{\calL^{(m)}_{1,\rad} u}{R^{(m)}}^2}\\
      &\leq \frac{1}{2}\bracket{\theta^{-2} c_1^2 +\theta^2\paren{\int \abs{u} \abs{\calL^{(m)}_{1,\rad} R^{(m)}} }^2}\\
      & \leq \frac{1}{2}\bracket{\theta^{-2} c_1^2
        +\theta^2\paren{\int \abs{u} \abs{\calL^{(m)}_{1,\rad}
            R^{(m)}}^{1/2}\abs{\calL^{(1)}
            R^{(m)}}^{1/2}  }^2}\\
      &\leq C\bracket{\theta^{-2} c_1^2 +  \theta^2\int \abs{\calL^{(m)}_{1,\rad} R^{(m)}}\abs{u}^2 }\\
      &\leq C \paren{\theta^{-2} c_1^2 + \theta^2\int
        e^{-\abs{y}}\abs{u}^2 }.
    \end{split}
  \end{equation}
  The other terms in which $u$ or $v$ appears once are similarly
  controlled.  Therefore,
  \[
  \calH^{(m)}(\eps, \eps) \geq \calH_{1,\rad}^{(m)}(u,u) +
  \calH_{2,\rad}^{(m)}(v,v) - C\theta^{-2}(c_1^2 + c_2^2 + d_1^2
  +d_2^2) -C \theta^2 \int e^{-\abs{y}}\abs{u+iv}^2,
  \]

  For the case $m=1$, we apply Proposition \ref{p:specprop} to get
  \[
  \begin{split}
    \calH^{(1)}(\eps, \eps) \geq &
    C_{(1)} \int{\abs{\nabla(u+iv)}^2}\\
    &+ \left(C_{(1)} - C \theta^2\right) \int{ e^{-\abs{y}}\abs{u+iv}^2 } - C\theta^{-2}(c_1^2 + c_2^2 + d_1^2 +d_2^2)\\
    \geq & \frac{C_{(1)}}{2} \int \abs{\nabla(u+iv)}^2 +
    e^{-\abs{y}}\abs{u+iv}^2 - C\theta^{-2}(c_1^2 + c_2^2 + d_1^2
    +d_2^2),
  \end{split}
  \]
  where we take $\theta >0$ sufficienty small. Finally,
  \[\begin{aligned}
    \int \abs{\nabla e^{i\theta}(u+iv)}^2 + e^{-\abs{y}}\abs{
      &e^{i\theta}(u+iv)}^2\\
    &\geq C\paren{ \int \abs{\nabla\eps}^2 + e^{-\abs{y}}\abs{\eps}^2
    } - \bigo\left(c_1^2 + c_2^2 + d_1^2 + d_2^2\right).
  \end{aligned}\]
\end{proof}

\appendix

\section{Almost-Self Similar Profiles}\label{Appendix-ProofOfAlmostSelfSimilar}
In this Appendix, we outline the proof of Proposition \ref{Prop-Qb},
showing modifications of the proof given in the case $m=0$,
\cite{MR-SharpUpperL2Critical-03, MR-UniversalityBlowupL2Critical-04,
  MR-SharpLowerL2Critical-06}. We then briefly discuss the proof of
Proposition \ref{Prop-Zb}. Recall that for $e^{ib\frac{r^2}{4}}\Qmb =
e^{im\theta}\Pmb(r)$ we have equation (\ref{Eqn-Pmb}),
\[
\laplacian \Pmb - \left(1+\frac{m^2}{r^2}-\frac{b^2}{4}r^2\right)\Pmb
+ \Pmb\abs{\Pmb}^2 = 0.
\]
This is not a scale-invariant equation, and there is no clear
representative solution.  Fibich and Gavish
\cite{FG-TheorySingularVortex-08} chose to consider the solution where
the boundary condition $\lim_{r\to 0}r^{-m}\Pmb(r) \neq 0$ is chosen
to minimize the amplitude of the asymptotic oscillation.
Since we intend to truncate anyways, it is more convenient to choose
boundary conditions,
\begin{equation}\label{Eqn-Pmb-Conditions}
  \Pmb(r) \left\lbrace\begin{aligned}
      \neq 0 &&\text{ for }&& 0 < r &< (1-\eta)R_b,\\
      = 0	&&\text{ for }&&  r &= (1-\eta)R_b.
    \end{aligned}\right.
\end{equation}
Recall that $R_b$ was chosen, (\ref{Eqn-Defn-Rb}), so that the strong
maximum principle applies to, $\laplacian -
\left(1+\frac{m^2}{r^2}-\frac{b^2}{4}r^2\right)$, on a region larger
than, $r \leq (1-\eta)R_b$.

\noindent {\bf Step 1:} Existence of $\Pmb$.


Following the argument of \cite[p605-606]{MR-SharpUpperL2Critical-03},
let ${\mathcal F}_\m$ denote the space of radial profiles of functions
in $H^1_\m$. That is, radial $H^1$ functions $f(x)$ for which
$x^{-1}f(x)\in L^2$.  Perform a constrained minimization of,
\[
2\,J_b[w] = \int{\abs{\grad w}^2} + \int{\abs{w}^2} +
m^2\int{\abs{\frac{w}{r}}^2} -\frac{b^2}{4}\int{\abs{rw}^2},
\]
over the subspace of finite-variance functions in $\calFm$ with
$w((1-\eta)R_b) = 0$ and $\int{\abs{w}^4} = 1$, where all integrals
are taken over a larger compact set, for example $r \leq
(1-\eta^2)R_b$. Note that $J_b$ is coercive on $\Hm(\real^2)$,
\begin{equation}\label{Eqn-AppA-FuncLowerBound}
  J_b[w] \geq C(\eta)\left(\int{\abs{\grad w}^2} + \int{\abs{w}^2} +
    m^2\int{\abs{\frac{w}{r}}^2}\right).
\end{equation}
This minimizing sequence can be assumed to converge weakly in
$H^1_\m$, which is simply a subspace of $H^1(\R^2)$, and thus strongly
in $L^4$ due to Sobolev embedding on a compact domain. Here we use
that equation (\ref{Eqn-NLS}) is energy subcritical. The Lagrange
multiplier of the Frechet derivative shows that (\ref{Eqn-Pmb}) is
satisfied. Interior regularity estimates show that the weak limit is
$C^3$ on $r < (1-\eta)R_b$. The weak limit is also strictly positive
due to $w((1-\eta)R_b) = 0$ and the maximum principle.

\noindent {\bf Step 2:} $L^\infty$ Estimates, Uniform in $b$.


There exists a fixed constant $C>0$ for all $\abs{b}>0$ sufficiently
small so that,
\begin{equation}\label{Eqn-AppA-P_LInfty}
  \abs{\Pmb}_{L^\infty} \leq C.
\end{equation}
\noindent Moreover, there is uniform decay of the tail of the
solutions. For the same $b$,
\begin{equation}\label{Eqn-AppA-P_DecayTail}
  \sup_{\abs{b}\sim 0}\abs{\Pmb}_{L^\infty(r>R)} \longrightarrow 0 \text{ as }{R\to+\infty}.
\end{equation}
Both bounds are proven in \cite[p606]{MR-SharpUpperL2Critical-03}.
Equation (\ref{Eqn-AppA-P_LInfty}) is a simple consequence of the
variational characterization of Step 1, whereas to prove equation
(\ref{Eqn-AppA-P_DecayTail}), truncate to $r > R$, treat
$r^\frac{N-1}{2}\abs{\Pmb}$ as a one-dimensional function, and control
by the standard Sobolev embedding $H^\frac{1}{2}(\real)
\hookrightarrow L^\infty(\real)$.

\noindent {\bf Step 3:} Local Convergence to $\Rm$ (in $C^3$).


As $b\to 0$, $\Pmb$ converges weakly to some positive radial function
$P$, with decay to $0$ as $r\to+\infty$, and which satisfies,
$\laplacian P - \left(1+\frac{m^2}{r^2}\right)P + P\abs{P^2}$. This
characterizes $P$ as the unique groundstate $\Rm$,
\cite{Mizumachi-VortexSolitons-05}.  Moreover, due to interior
regularity estimates, on any compact set the convergence of $\Pmb$ is
strong in $C^3$, up to a subsequence in $b$.

\noindent {\bf Step 4:} Uniform Convergence to $\Rm$ (in $C^3$ with
exponential weight)


Here we adapt the argument of
\cite[p658-659]{MR-UniversalityBlowupL2Critical-04}. Consider the
operator $\calK = \laplacian -
\left(1+\frac{m^2}{r^2}-\frac{b^2}{4}r^2 - \frac{\eta^2}{2}\right)$,
which satisfies the maximum principle on $1 < r < (1-\eta)R_b$, for
any $\eta>0$ sufficiently small. Restate (\ref{Eqn-Pmb}) as,
\begin{equation}\label{Eqn-AppA-Pmb}
  \calK \Pmb = \frac{\eta^2}{2}\Pmb-\left(\Pmb\right)^3.
\end{equation}
Consider the new function $f_{b}(r) = e^{-(1-\eta)R_b
  \Theta\left(\frac{r}{R_b}\right)}$, with,
\[
\Theta(\xi) = {\mathds 1}_{0 < \xi < 1}\int{\sqrt{1 - z^2}\,dz} +
{\mathds 1}_{1 \leq \xi}\,\Theta(1)\,\xi.
\]
Note the dependence on $m$. By direct calculation,
\[\begin{aligned}
  f_b^{-1}\calK f_b = (1-\eta)\frac{ \frac{r}{R_b^2} }{\sqrt{ 1 -
      \left(\frac{r}{R_b}\right)^2 }}
  &+(1-\eta)^2\left(1-\left(\frac{r}{R_b}\right)^2\right)\\
  &- \frac{1}{r}\sqrt{ 1 - \left(\frac{r}{R_b}\right)^2 } -
  \left(1+\frac{m^2}{r^2}-\frac{b^2}{4}r^2 - \frac{\eta^2}{2}\right).
\end{aligned}\] We now approximate each term on the region
$\frac{1}{\eta} < r \leq (1-\eta)R_b$,

\[\begin{aligned}
  f_b^{-1}\calK f_b \leq \frac{(1-\eta)^2}{R_b} &+\left((1-\eta)^2 -
    1\right)\left(1-\left(\frac{r}{R_b}\right)^2\right)
  +\left(\frac{b^2}{4}-\frac{1}{R_b^2}\right)r^2 \\
  &- \eta^\frac{3}{2}\sqrt{2-\eta} - m^2\eta^2 + \frac{\eta^2}{2}.
\end{aligned}\] Recall that, $R_b =
\frac{\sqrt{2+2\sqrt{1+b^2m^2}}}{b}$. By assuming $b>0$ is
sufficiently small with respect to $\eta$, we conclude $f_b^{-1}\calK
f_b$ is strictly negative for the given range of $r$.

From Step 2, and the exponential decay of $\Rm$, there exists a fixed
value $r_0 > \frac{1}{\eta}$ such that for all $b>0$ sufficiently
small,
\[\begin{aligned}
  \frac{\eta^2}{2}\Pmb - \left(\Pmb\right)^3 > 0 &&\text{ for } && r
  \in \Omega = r_0 < r < (1-\eta)R_b.
\end{aligned}\] We have shown that $\calK\left(c\, f_b -
  \Pmb\right)<0$ for $r\in\Omega$ and any arbitrary constant
$c>0$. Now we note that,
\[
\lim_{b\to 0} f_b(r_0) = e^{-(1-\eta)r_0} > 0,
\]
so that we may choose our constant $c = 2\Rm(r_0)e^{+(1-\eta)r_0}$
and, with our boundary condition (\ref{Eqn-Pmb-Conditions}), conclude
that,
\[
\left. c\, f_b(r) - \Pmb(r) \right|_{\partial\Omega} > 0.
\]
The maximum principle may now be applied. The same argument can be
applied to $\Rm$, $b=0$, and the weight $f(r) = e^{-(1-\eta)r}$.
With Step 3, this proves the first precursor of
(\ref{Prop-Qb-closeToQ}),
\begin{equation}\label{Eqn-AppA-GoalStep4}\begin{aligned}
    \left.
      \norm{e^{(1-C\eta)R_b\Theta\left(\frac{r}{R_b}\right)}\left(\Pmb-\Rm\right)}_{C^3}
    \right.  \longrightarrow 0 &&\text{ as } && b \rightarrow 0.
  \end{aligned}\end{equation}

To prove the bound for the energy, (\ref{Prop-Qb-EnerMomentum}), note
that without loss of generality $(1+C\eta)(1-a) =
(1-\delta)<1$. Introduce a new operator $\calK$ and function $f_b$ in
terms of $\delta$ in place of $\eta$ and argue Step 4 again. In
particular, we may assume that $r_0 < \delta R_b \ll (1-\eta)^2R_b < r
< (1-\eta)R_b \ll (1-\delta)R_b$.

\noindent {\bf Step 5:} Uniqueness of $\Pmb$; Continuity in $b$


For fixed $b_0>0$ sufficiently small, and $b\approx b_0$, consider,
\begin{equation}\label{Eqn-AppA-Tbb}
  T_{b,b_0} = \left(\frac{R_b}{R_{b_0}}\right)\Pmb\left(\frac{R_b}{R_{b_0}}r\right)
\end{equation}
Then $T_{b,b_0} \in \calFm$ and vanishes for $r=(1-\eta)R_{b_0}$, and
we consider the differential, $T_\Delta = T_{b,b_0} - \Pm_{b_0}$, with
the same domain. The goal is to prove,
\begin{equation}\label{Eqn-AppA-GoalStep5}
  \norm{e^{im\theta}T_\Delta(r)}_{H^1(\real^2)} \leq C\frac{\abs{b-b_0}}{b_0},
\end{equation}
for some fixed constant $C$. To do so, consider the equation for
$T_\Delta$ written as,
\begin{equation}\label{Eqn-AppA-TDelta}\begin{aligned}
    \left(\Lm_+ - \frac{b_0^2}{4}r^2\right) T_\Delta =
    &-\left(\left(1-\frac{R_b^2}{R_{b_0}^2}\right)\left(1 -
        \frac{b_0^2}{4}r^2\right)
      +\frac{R_b^2}{R_{b_0}^2}\frac{b_0^2 - b^2\frac{R_b^2}{R_{b_0}^2}}{4}r^2 \right) T_\Delta\\
    &-3R_\m^2 T_\Delta +\left(T_\Delta + \Pm_{b_0}\right)^3 - \left(\Pm_{b_0}\right)^3\\
    &+\left( \left(1-\frac{R_b^2}{R_{b_0}^2}\right)\left(1 -
        \frac{b_0^2}{4}r^2\right) +\frac{R_b^2}{R_{b_0}^2}\frac{b_0^2
        - b^2\frac{R_b^2}{R_{b_0}^2}}{4}r^2 \right) \Pm_{b_0},
  \end{aligned}\end{equation}
where $\Lm_+$ is the linerized operator from equation
(\ref{Eqn-Defn-L}). We will use ${\mathcal F}_{b,b_0}$ to denote the
final right hand term of (\ref{Eqn-AppA-TDelta}).  Note that in the
case $m=0$, and thus $R_b = \frac{2}{b}$, the final multiples of
$T_\Delta$ and $\Pm_{b_0}$ collapse. All three right hand terms of
(\ref{Eqn-AppA-TDelta}) are bounded in the same way as in
\cite[p609]{MR-SharpUpperL2Critical-03}, with only minor
adaptations\footnotemark.
To conclude the argument from \cite{MR-SharpUpperL2Critical-03} and
establish (\ref{Eqn-AppA-GoalStep5}) there only remains to show the
following Lemma:

\footnotetext{The terms due to $R_b \neq \frac{2}{b}$ have no
  effect. Part of the error term $G_1(R)$ that appears in
  \cite{MR-SharpUpperL2Critical-03} has been moved to the left hand
  side of (\ref{Eqn-AppA-TDelta}), so that the constant $A_0$ that
  appears in \cite{MR-SharpUpperL2Critical-03} can be ignored.  }
\begin{lemma}\label{Lemma-AppA-Univ212}
  Let $\mu_+<0$ be the lowest eigenvalue of $\Lm_+$, and $\phi_+\in
  L^2$ the corresponding normalized eigenvector. For $b>0$
  sufficiently small with respect to $\eta$, and assuming $\eta>0$ is
  itself sufficiently small,
  \[
  \inner{\left(\Lm_+-\frac{b^2}{4}r^2\right)w}{w} \geq
  \delta_+\norm{w}_{H^1}^2 - \frac{1}{\delta_+}\inner{w}{\phi_+}^2,
  \]
  for $\delta_+>0$ constant and any $w\in H^1_\m$ vanishing at
  $r=(1-\eta)R_b$.
\end{lemma}
Lemma \ref{Lemma-AppA-Univ212} is analogous to \cite[equation
(212)]{MR-UniversalityBlowupL2Critical-04}, and is adapted from Lemma
\ref{Lemma-Maris} by using a cutoff and the exponential decay of
$\phi_+$. Details can be found,
\cite[p660]{MR-UniversalityBlowupL2Critical-04}.

\noindent {\bf Step 6:} Frechet Derivative on Fixed Domain


The aim is to prove that there exists,
\begin{equation}\label{Eqn-AppA-GoalStep6}\begin{aligned}
    \left.\frac{\partial}{\partial_b}T_{b,b_0}\right|_{b=b_0} \in
    H^1_\m.
  \end{aligned}\end{equation}
We will follow the argument of
\cite[p610]{MR-SharpUpperL2Critical-03}, and revisit equation
(\ref{Eqn-AppA-TDelta}). In the limit $b \to b_0$ we have, with
respect to $L^2$-norm,
\begin{equation}\label{Eqn-AppA-TDeltaTilde}\begin{aligned}
    \left(\Lm_+ - \frac{b_0^2}{4}r^2\right)\frac{T_\Delta}{b-b_0} = 0
    -3\left(\left(\Rm\right)^2 -
      \left(\Pm_{b_0}\right)^2\right)\frac{T_\Delta}{b-b_0}
    +\left.\frac{\partial}{\partial_b}{\mathcal
        F}_{b,b_0}\right|_{b=b_0}.
  \end{aligned}\end{equation}
Note that by direct calculation,
\[
\left.\frac{\partial}{\partial_b}{\mathcal F}_{b,b_0}\right|_{b=b_0} =
\frac{2}{b_0}\left( 1-\frac{b_0^2}{4}r^2
  -\frac{1}{2}\frac{\sqrt{1+b_0^2m^2}-1}{\sqrt{1+b_0^2m^2}}
\right)\Pm_{b_0},
\]
and clearly exists.  To show equation (\ref{Eqn-AppA-GoalStep6}), we
recall from Step 5 that, for $b_0>0$ sufficiently small,
$\Lm_+-\frac{b_0^2}{4}r^2$ is invertible over the subspace of $L^2_\m$
functions that vanish at $r=(1-\eta)R_b$.

\noindent {\bf Step 7:} Uniform Bound for
$\left.\partial_bT_{b,b_0}\right|_{b=b_0}$ (in $C^2$ with exponential
weight)


Revisit equation (\ref{Eqn-AppA-TDelta}), again in the limit $b\to
b_0$ with respect to $L^2$ norm,
\[\begin{aligned}
  \left(\Lm_+ - \frac{b_0^2}{4}r^2 + 3\left(\left(\Rm\right)^2 -
      \left(\Pm_{b_0}\right)^2\right)\right)
  \left.\frac{\partial}{\partial b}T_{b,b_0}\right|_{b=b_0} &=
  +\left.\frac{\partial}{\partial_b}{\mathcal
      F}_{b,b_0}\right|_{b=b_0}.
\end{aligned}\] Similar to Step 4, we apply a maximum principle
argument on the region $\frac{1}{\eta} < r \leq (1-\eta)R_b$ to prove,
\[
\norm{ e^{(1-C\eta)R_b\Theta\left(\frac{r}{R_b}\right)}
  \left.\frac{\partial}{\partial b}T_{b,b_0}\right|_{b=b_0}
}_{C^2(r<(1-\eta)R_b)} \lesssim \frac{1}{b_0}.
\]
The full argument is the same as
\cite[p610-611]{MR-SharpUpperL2Critical-03} with only minor
adaptations.

\noindent {\bf Step 8:} Uniform Bound for
$\left.\partial_b\PmbT\right|_{b=b_0}$ (in $C^2$ with exponential
weight)


Let $\PmbT = \phi_b\Pmb$ where $\phi_b$ are the smooth cutoff
functions,
\begin{equation}\label{DefnEqn-phiB}
  \phi_b(r) = \left\{\begin{aligned}
      1 & \text{ for } r < (1-\eta)^2R_b\\
      0 & \text{ for } r > (1-\eta)R_b,
    \end{aligned}\right.
\end{equation}
with the good behaviour,
$\norm{\grad\phi_b}_{L^\infty}+\norm{\laplacian\phi_b}_{L^\infty} \to
0$, as $b\to 0$.  Alternately,
\begin{equation}\label{Eqn-AppA-AlterPTilde}
  \PmbT = \left(\phi_b - \phi_{b_0}\right)\Pmb + \phi_{b_0}\left(\Pmb-\Pm_{b_0}\right) + \PmT_{b_0}.
\end{equation}
The goal is to prove that,
\begin{equation}\label{Eqn-AppA-GoalStep8}\begin{aligned}
    \norm{e^{(1-C\eta)R_b\Theta\left(\frac{r}{R_b}\right)}\frac{\partial}{\partial_b}\PmbT}_{C^2(\real^2)}
    \longrightarrow 0 &&\text{ as }&& b\to 0.
  \end{aligned}\end{equation}
which is the second precursor to (\ref{Prop-Qb-closeToQ}).
Regarding the first right hand term of (\ref{Eqn-AppA-AlterPTilde}),
we may re-express $\Pmb$ in terms of $T_{b,b_0}$. Then by Step 7 and
the support of $\phi_b-\phi_{b_0}$, the contribution from that term is
neglible.  The remaining term,
$\phi_{b_0}\left(\Pmb-\Pm_{b_0}\right)$, is treated with calculations
similar to those applied to $T_\Delta$ in Steps 5, 6 and 7. Details
can be found, \cite[p611-612]{MR-SharpUpperL2Critical-03}.

\noindent {\bf Step 9:} Supercritical Mass

The proof of (\ref{Prop-Qb-Mass}) is due to \cite[Lemma
1]{MR-SharpLowerL2Critical-06}. 
Here, we give a summary for the reader's convenience.  To begin, note
from equation (\ref{Eqn-Pmb}) that $\PmbT$ is formally a function of
$b^2$. Then from Step 8 and the chain rule we conclude that, with an
exponential weight, $\partial_{(b^2)}\PmbT$ is bounded in $C^2$. From
equation (\ref{Eqn-Pmb}) it can be shown in the limit $b\to 0$ that,
\begin{equation}\label{Eqn-AppA-Pb2}
  L_+\frac{\partial}{\partial(b^2)}\PmbT = \frac{r^2}{4}\PmbT.
\end{equation}
Consider then a product of (\ref{Eqn-AppA-Pb2}) by $\Lambda \Rm$,
\[\begin{aligned}
  \frac{1}{4}\int{\abs{x}^2\abs{\Rm}^2\,dx}
  &&&= -\frac{1}{4}\inner{r^2\,\Rm}{\Lambda \Rm}\\
  &&&=-\lim_{b\to 0}\inner{L_+\partial_{(b^2)}\Pmb}{\Lambda \Rm}\\
  &&&=-\lim_{b\to 0}\inner{\partial_{(b^2)}\Pmb}{-2\Rm} &&= \lim_{b\to
    0}\partial_{b^2}\norm{\Pmb}_{L^2(\real^2)}.
\end{aligned}\]

This concludes our summary of the proof of Proposition \ref{Prop-Qb}.
\begin{proof}[Proof of Proposition \ref{Prop-Zb}.]

  Apply the point transformation, $e^{ib\frac{r^2}{4}}\Zmb =
  e^{-im\theta}r^mZ(r)$. Then equation (\ref{Prop-Zb-Eqn}) reads,
  \[
  \partial_r^2Z + \frac{2m+1}{r}\partial_rZ - Z + \frac{b^2r^2}{4}Z =
  \widetilde{\Psi}_b,
  \]
  where $r^m\widetilde{\Psi}_b = \laplacian\phi_b\Pmb
  +\grad\phi_b\cdot\grad\Pmb + \left(\phi_b^3-\phi_b\right)\Pmb$. The
  arguments of \cite[Appendix E]{MR-UniversalityBlowupL2Critical-04}
  and \cite[Appendix A]{MR-SharpLowerL2Critical-06}, then prove a
  version of Proposition \ref{Prop-Zb} for $e^{ib\frac{r^2}{4}}Z(r)$,
  as a radial function on $\real^{2m+2}$. By accounting for the
  equivalences of norms, this proves Proposition \ref{Prop-Zb}.

\end{proof}

\section{Details of Numerical Methods}
\label{s:numerics}
Our numerical methods closely follow those detailed in
\cite{Marzuola2010}, employing the Fortran 90/95 boundary value
problem software described in \cite{shampine2006user}.  We briefly
review it here.

The software is designed to solve two point boundary value problems of
the form
\begin{equation}
  \frac{d}{dr}\mathbf{y} = \frac{1}{r}S\mathbf{y} + \mathbf{f}(r,
  \mathbf{y}),
\end{equation}
by nonlinear collocation. Note that the algorithm handles $r^{-1}$
singuralities.  All of our computations were performed on the domain
$[0, 50]$ with tolerance $10^{-10}$.

Codes that can be used to reproduce the computations presented here
are available at
\url{http://www.math.toronto.edu/simpson/files/vortex_dist.tgz}.

\subsection{Point Transformations}
Unfortunately, the equation for the vortex state,
\eqref{Eqn-VortexSoliton-R}, and the operators $\calLm_{1,\rad}$ and
$\calLm_{2,\rad}$, include $r^{-2}$ singularities.  We address this
with the point transformation
\begin{equation}
  R^\m(r) = r^m \tilde{R}^\m(r).
\end{equation}
Similarly, $U = e^{im\theta}r^m \tilde{U}^{(m)}(r)$ for any of the dependent
variables.  With this transformation, the vortex equation becomes,
\begin{subequations}
  \begin{gather}
    (\tilde{R}^\m)'' + \frac{2m+1}{r}(\tilde{R}^\m)' - \tilde{R}^\m + r^{2m} (\tilde{R}^\m)^3 = 0,\\
    (\tilde{R}^\m)'(0) = 0, \quad \lim_{r\to \infty}\tilde{R}^\m(r) =
    0,
  \end{gather}
\end{subequations}
and the operators $\calLm_{1,\rad}$,$\calLm_{2,\rad}$ become,
\begin{subequations}\begin{align}
    \begin{split}
      \calLm_{1,\rad} U &= r^{m} \set{- \tilde{U}'' -
        \frac{2m+1}{r}\tilde{U}' + 3 r^{2m} \tilde{R}^\m (m
        \tilde{R}^\m + r (\tilde{R}^\m)' )
        \tilde{U}}\\
      &= r^m \tilde{\calL}_1 \tilde{U}
    \end{split}\\
    \begin{split}
      \calLm_{2,\rad} Z &= r^{m} \set{- \tilde{U}'' -
        \frac{2m+1}{r}\tilde{U}' + r^{2m} \tilde{R}^\m (m
        \tilde{R}^\m+ r (\tilde{R}^\m)' )
        \tilde{U}}\\
      &=r^m \tilde{\calL}_2 \tilde{Z}
    \end{split}
  \end{align}
\end{subequations}

The right hand sides of \eqref{e:ip_bvps} conveniently become,
\begin{align}
  R^{(m)} &= r^m \tilde{R}^{(m)}\\
  \Lambda R^{(m)} & = r^m\set{(m+1) \tilde{R}^{(m)}+ r
    (\tilde{R}^{(m)})'}\\
  \begin{split}
    \Lambda^2 R^{(m)} & = r^m\set{(m+1)^2 \tilde{R}^{(m)}+(3+2m) r
      (\tilde{R}^{(m)})' +
      r^2 (\tilde{R}^{(m)})''}\\
    &= r^m\set{\bracket{(m+1)^2+r^2 }\tilde{R}^{(m)}+2 r
      (\tilde{R}^{(m)})' - r^{2(m+1)} (\tilde{R}^{(m)})^3 }
  \end{split}
\end{align}

\subsection{Artificial Boundary Conditions}
As the algorithm is designed to compute on finite intervals of
$[a,b]$, we must compute on $[0, \rmax]$, where $\rmax$ is
sufficiently large.  This neccessitates the introduction of an
artificial boundary condition on $\tilde{R}^{(m)}$, the vortex state,
and $U_j$ and $Z_j$ solving the boundary value problems
\eqref{e:ip_bvps}.  The analogous question in the index function
computations is verifying that there are no zeros beyond $\rmax$ which
might have been missed.

To develop the artifiical boundary conditions, we examine the
asymptotic behaviour of the solutions, using that potential terms are
exponentially decaying.  For the vortex state,
\begin{equation}
  \label{e:vortex_abc}
  \tilde{R}^{(m)}(r) \propto r^{-m -\frac{1}{2}}e^{-r}
\end{equation}
This gives us the boundary condition at $\rmax$
\begin{equation}
  (\tilde{R}^{(m)})'(\rmax) + \paren{1 + \frac{2m+1}{2\rmax}}\tilde{R}^{(m)}(\rmax) = 0
\end{equation}
which is accurate to $\bigo(\rmax^{-2})$.

By similar analysis the solutions to the linear boundary value
problems, generically denoted by $W$, are
\begin{equation}
  W(r) \propto r^{-2m}
\end{equation}
as $r\to \infty$.  Thus
\begin{equation}
  \label{e:bvp_abc}
  \tilde{W}'(\rmax) + \frac{2m}{\rmax}\tilde{W}(\rmax)=0
\end{equation}
This too is accurate to $\bigo(\rmax^{-2})$.

\subsection{Verification of Results}
With these approximations, we solve the following sets of equations,
as single first order systems:
\begin{itemize}
\item The vortex $\tilde{R}^{(m)}$, and the index functions $U$ and
  $Z$,
\item The vortex $\tilde{R}^{(m)}$, the boundary value problem
  solutions $U_1$ and $U_2$, and the $K_j$ inner products.
\item The vortex $\tilde{R}^{(m)}$, the boundary value problem
  solutions $Z_1$ and $Z_2$, and the $J_j$ inner products.
\end{itemize}
In computing the index functions, or alternatively the inner products,
we are actually solving mixed initial value/boundary value problems.

We now present several {\it a postiori} checks on the accuracy of our
results.  All are based on checking that the behaviour of the
solutions for large $r$ is consistent with the anticipated asymptotic
behavior.
\subsubsection{Verification of the Vortex States}
Two related ways of checking that we have adequately computed the
vortex states are to examine its decay as $r$ becomes large and to see
that \eqref{e:vortex_abc} becomes small as $r\to \infty$.  For the
vortices appearing in Figure \ref{f:vortex_radial}, we plot these two
metrics in Figures \ref{f:vortex_decay} and \ref{f:vortex_abc}.  With
this artificial boundary condition, the exponential decay is well
captured.

\ifpdf
\begin{figure}
  \centering
  \includegraphics[width=3in]{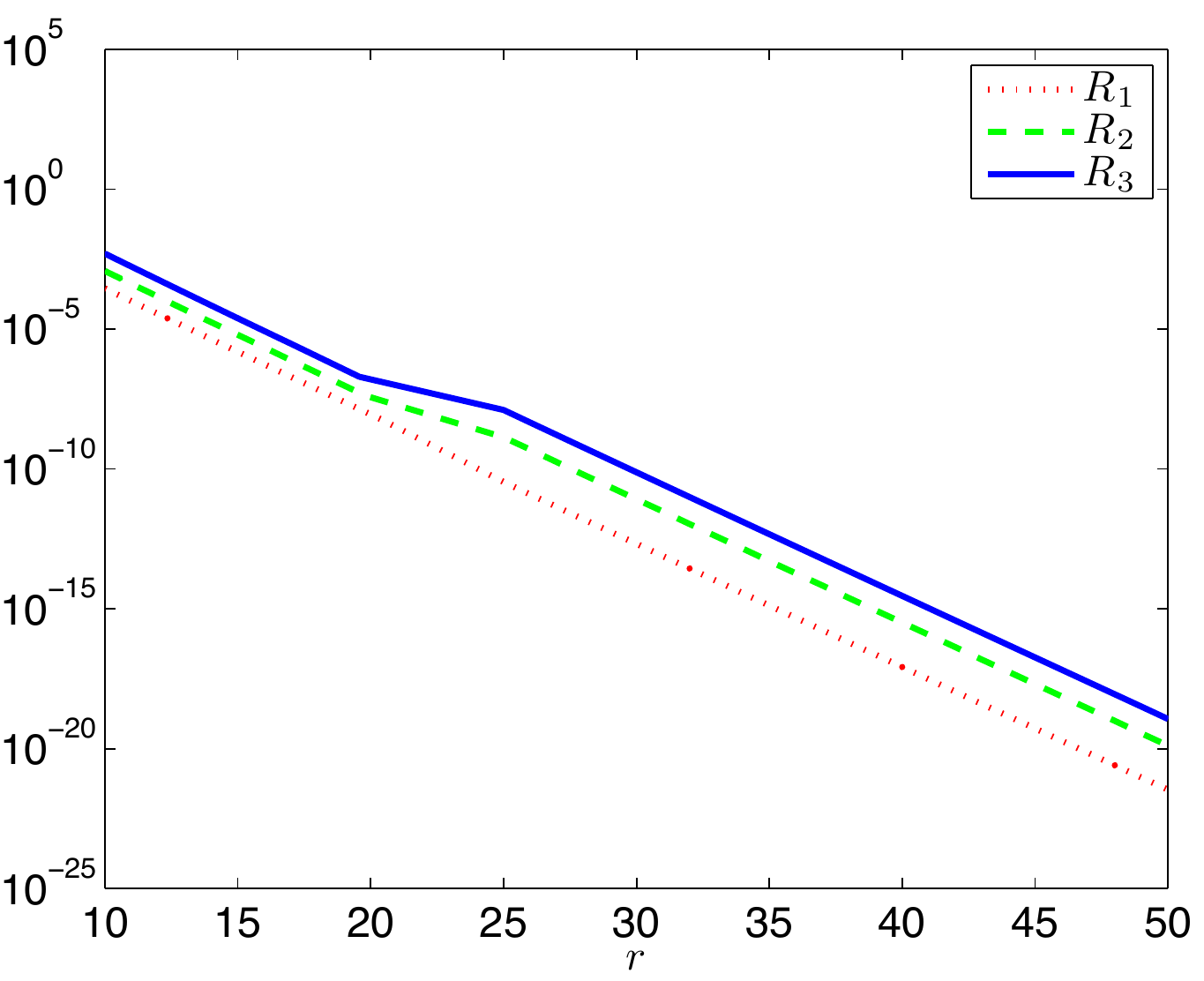}
  \caption{Behavior of the computed vortices as $r$ becomes large.  We
    recover the exponential decay.}
  \label{f:vortex_decay}
\end{figure}

	\begin{figure}
          \centering
          \includegraphics[width=3in]{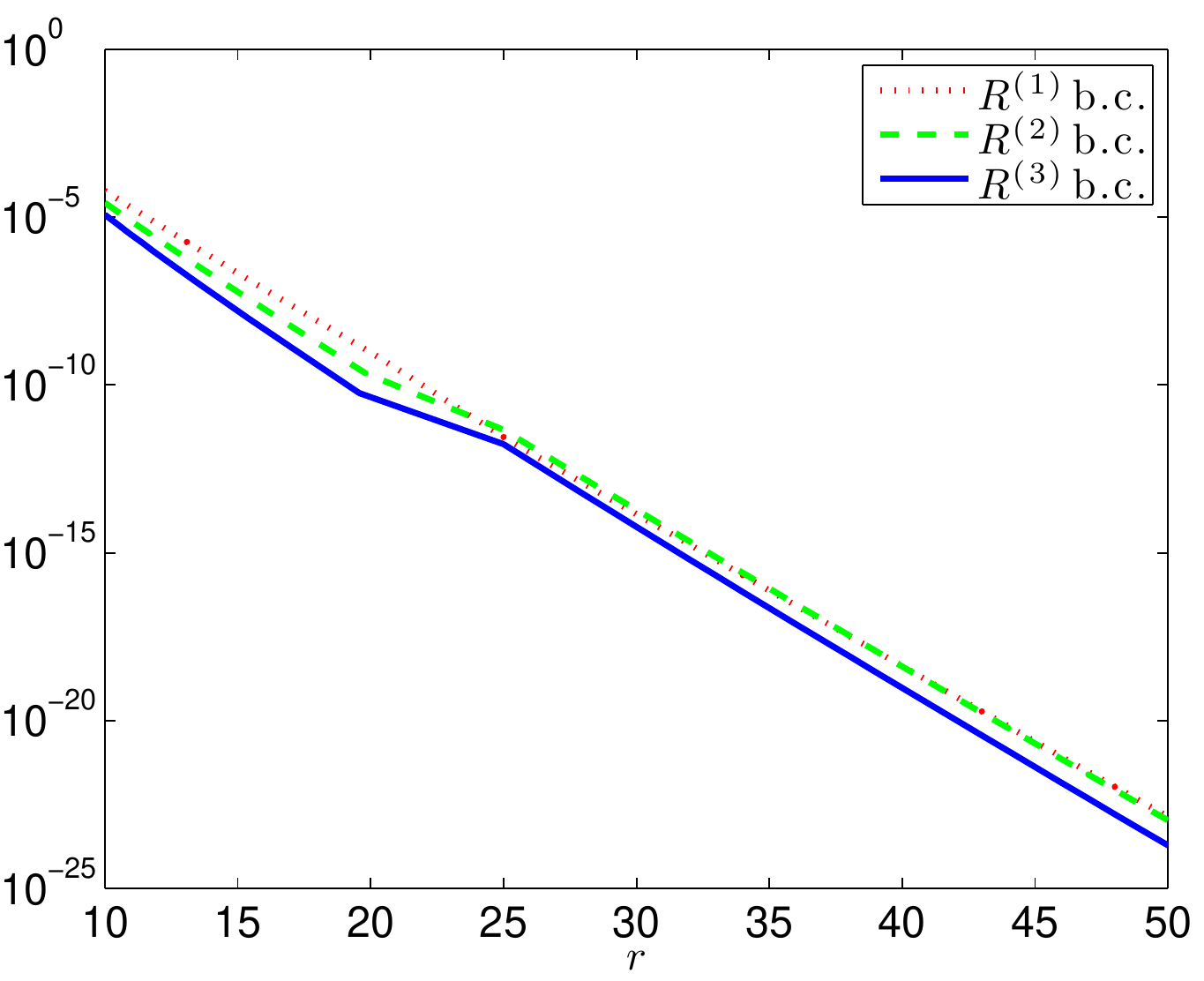}
          \caption{Asymptotically, the artificial boundary condition
            on the vortex is well satisfied.}
          \label{f:vortex_abc}
	\end{figure}
        \else (No figures in DVI) \fi

\subsubsection{Verification of the Index Count}
In counting the zeros of the index functions from Figure \ref{f:idx},
there is the concern that there may be another root located beyond
$\rmax$.  To assess this, we note that the asympotically free behavior
of $\tilde{U}$ and $\tilde{Z}$ is
\begin{subequations}
  \begin{align}
    \tilde{U}^{(m)} &\sim C_0^{(m)} + C_1^{(m)} r^{-2m}\\
    \tilde{Z}^{(m)} &\sim D_0^{(m)} + D_1^{(m)} r^{-2m}
  \end{align}
\end{subequations}
We can estimate these constants by noting
\begin{subequations}
  \begin{align}
    \frac{\tilde{U}' r^{1+2m}}{-2m} &\sim C_1^{(m)}\\
    \tilde{U} + \frac{\tilde{U}' r}{2m} &\sim C_0^{(m)}\\
    \frac{\tilde{Z}' r^{1+2m}}{-2m} &\sim D_1^{(m)}\\
    \tilde{Z} + \frac{\tilde{Z}' r}{2m} &\sim D_0^{(m)}
  \end{align}
\end{subequations}
These constants are plotted in Figure \ref{f:idx_consts}.  As they
show, we have certainly computed into the ``free'' equation regime.
More importantly, since $C^{(m)}_0 >0$ and $D^{(m)}_0<0$ in all cases,
we should not expect any additional zeros in the $U^{(m)}$ or
$Z^{(m)}$ functions appearing in Figure \ref{f:idx}.

\ifpdf
\begin{figure}
  \centering
  \subfigure[$m=1$]{\includegraphics[width=2.35in]{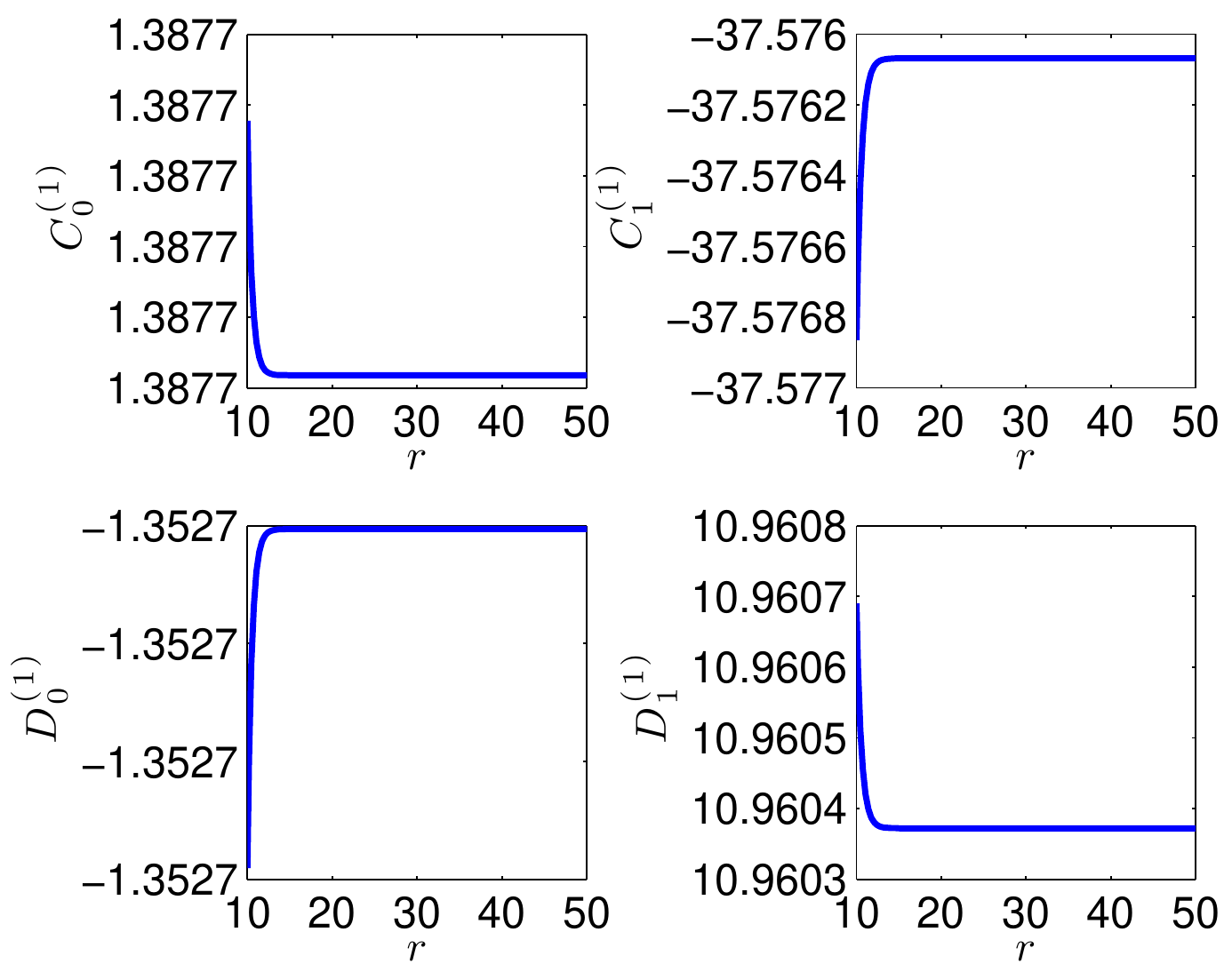}}
  \subfigure[$m=2$]{\includegraphics[width=2.35in]{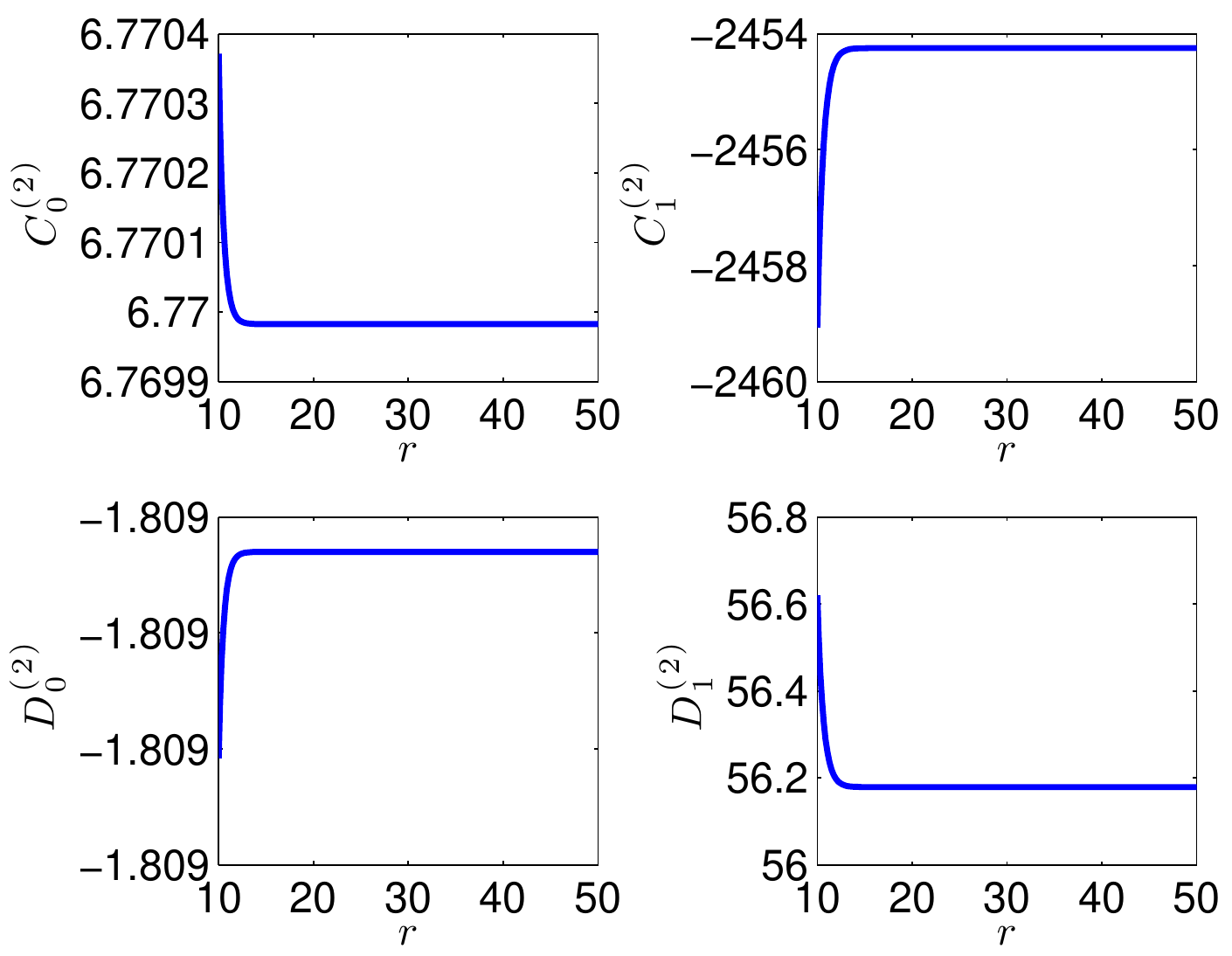}}
  \subfigure[$m=3$]{\includegraphics[width=2.35in]{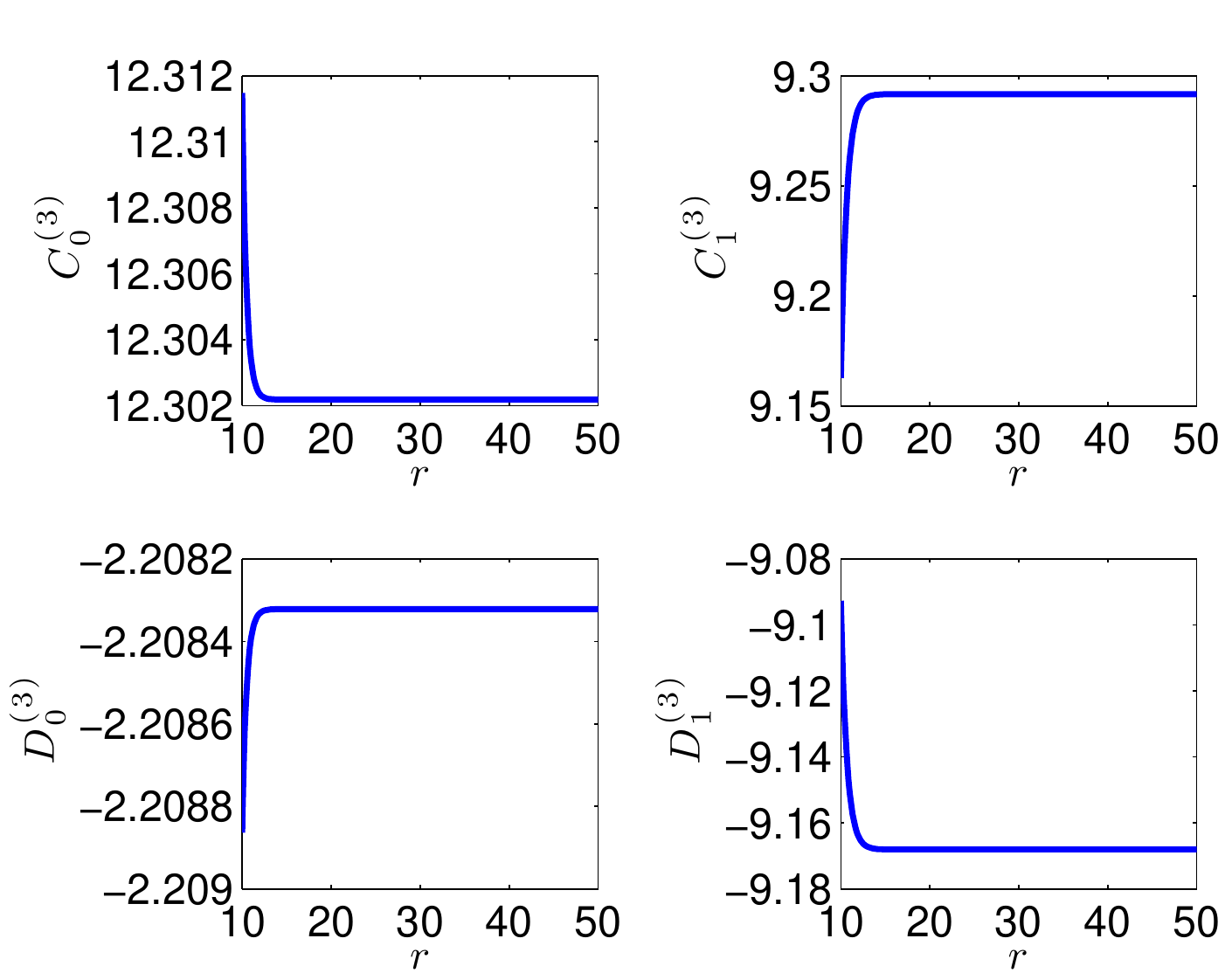}}
  \caption{The asymptotical index function constants for winding
    numbers $m=1,2,3$.}
  \label{f:idx_consts}
\end{figure}
\else (No figures in DVI) \fi

\subsubsection{Verification of the Inner Products}
For the inner product computations, we verify that in solving the
boundary value problems, $U_l^{(m)}$, $Z_l^{(m)}$ adequately satisfy
the artificial boundary conditions \eqref{e:bvp_abc}, and that the
$K_l$, $J_l$ values are ``constant''.  The check on the boundary
conditions is given in Figures \ref{f:ipp_abcs} and \ref{f:ipm_abcs}.
As these figures show, \eqref{e:bvp_abc} is well approximated.

\ifpdf
\begin{figure}
  \centering
  \subfigure[$m=1$]{\includegraphics[width=2.35in]{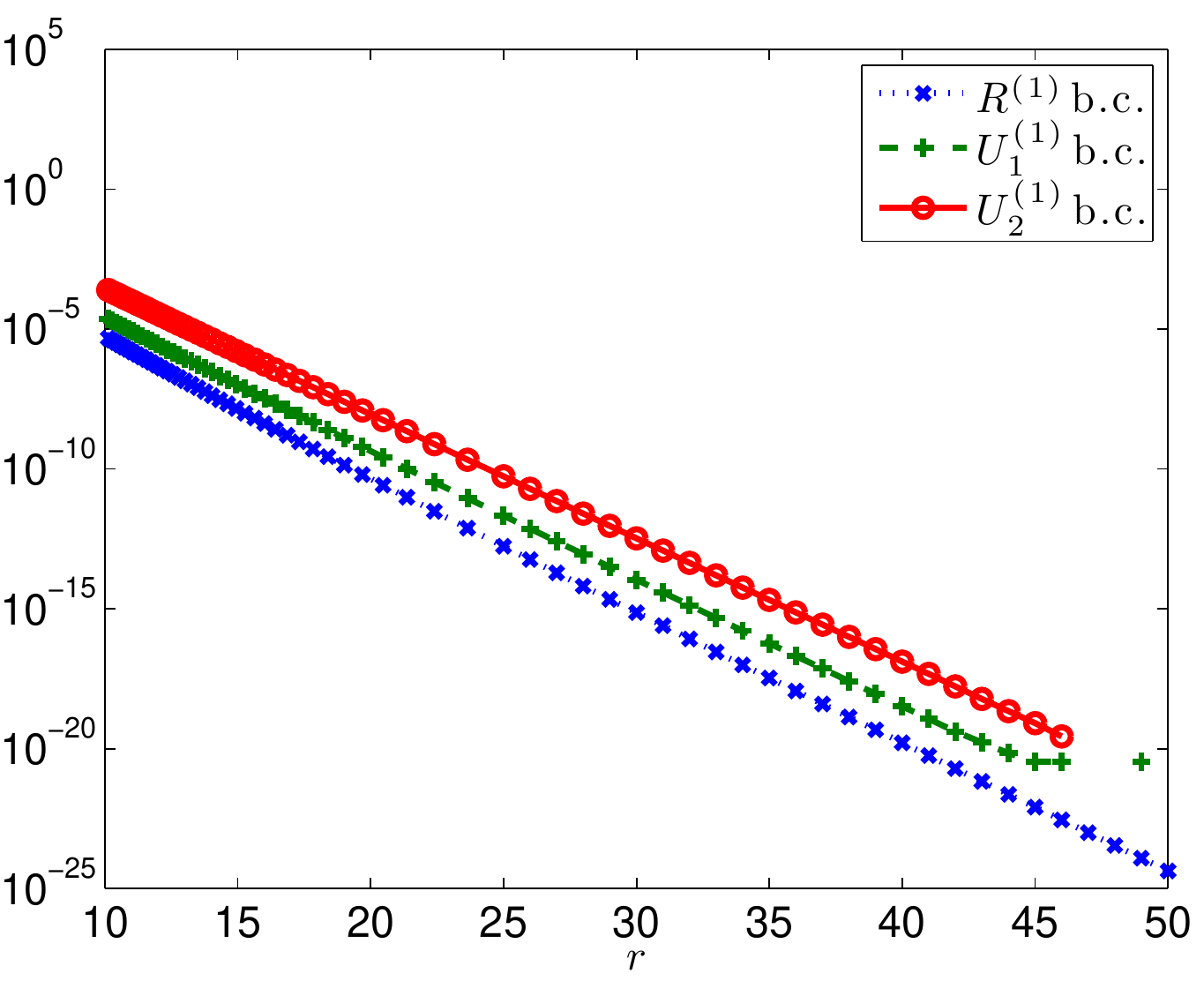}}
  \subfigure[$m=2$]{\includegraphics[width=2.35in]{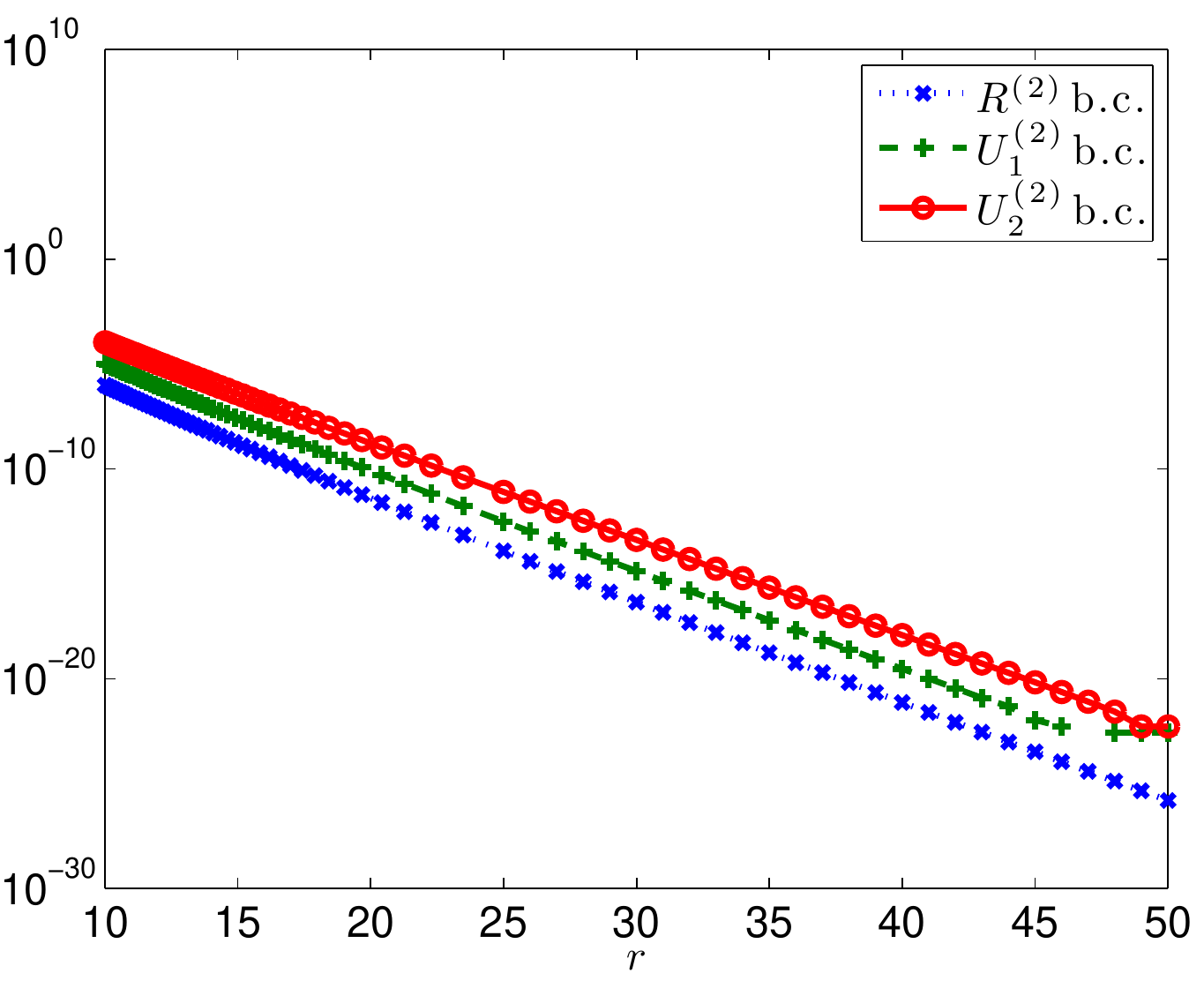}}
  \subfigure[$m=3$]{\includegraphics[width=2.35in]{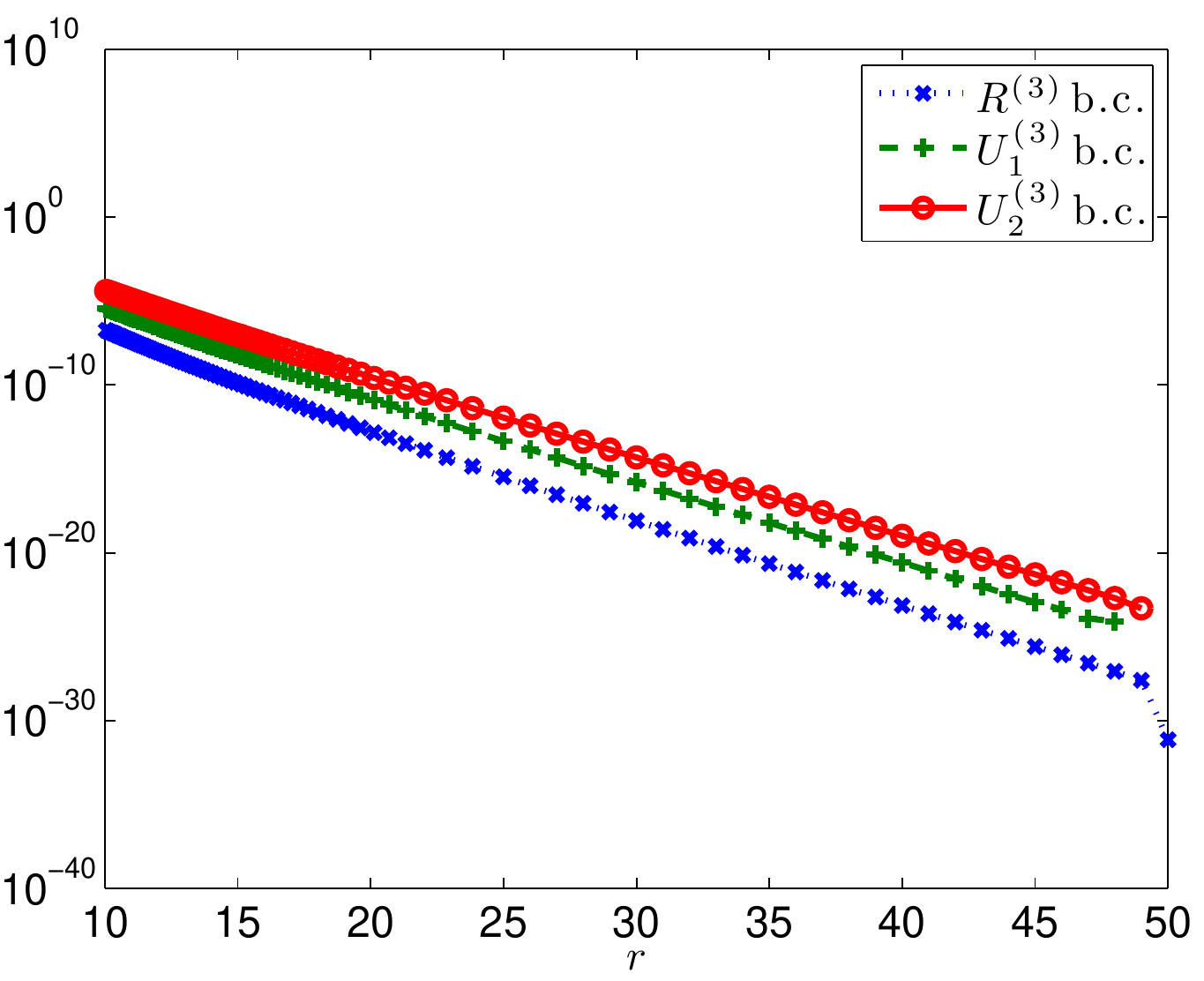}}
  \caption{Check of the artificial boundary conditions on the
    $U_l^{(m)}$ functions.}
  \label{f:ipp_abcs}
\end{figure}
	
	\begin{figure}
          \centering
          \subfigure[]{\includegraphics[width=2.35in]{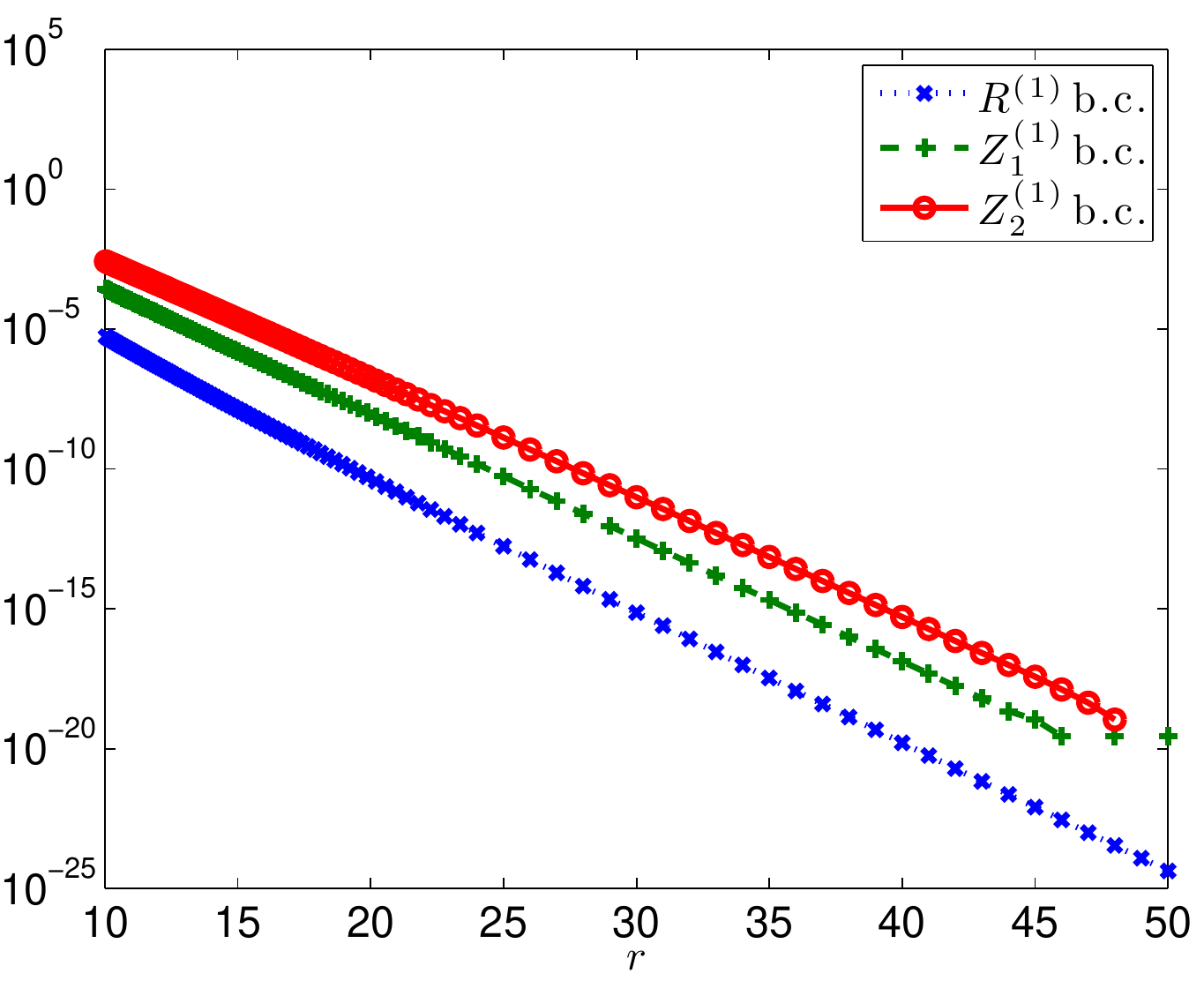}}
          \subfigure[]{\includegraphics[width=2.35in]{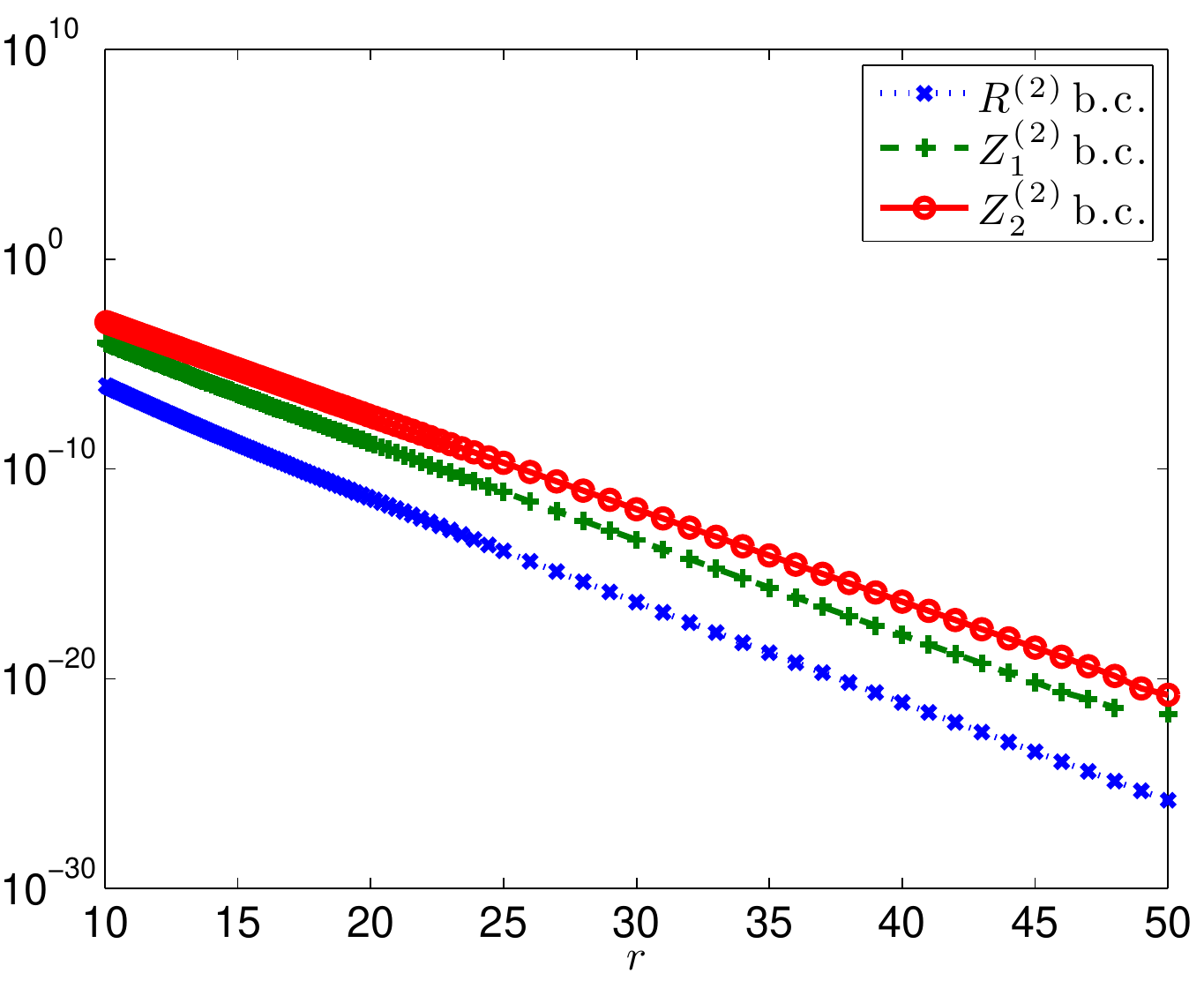}}
          \subfigure[]{\includegraphics[width=2.35in]{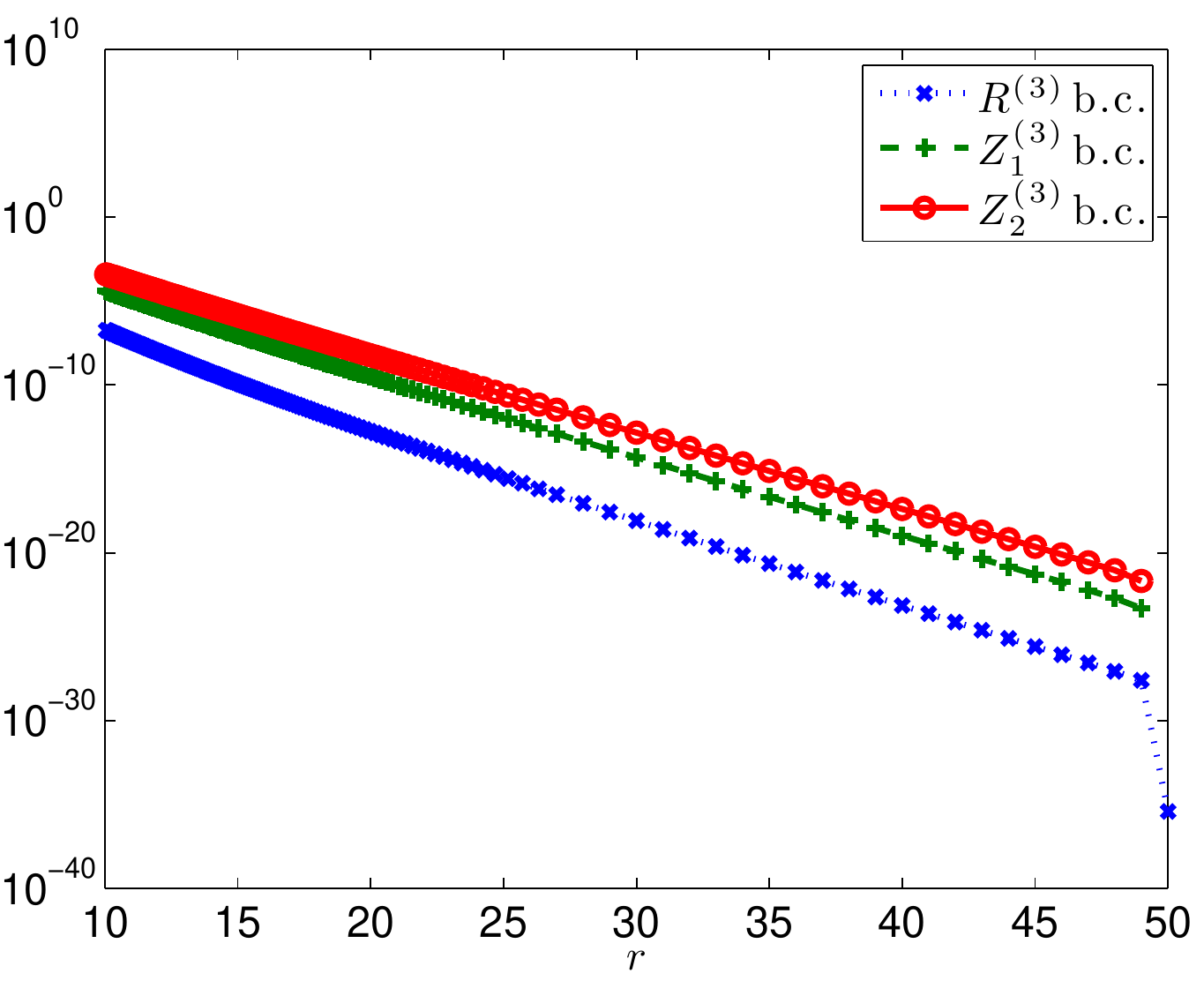}}
          \caption{Check of the artificial boundary conditions on the
            $Z_l^{(m)}$ functions.}
          \label{f:ipm_abcs}
	\end{figure}
        \else (No figures in DVI) \fi

        In computing the inner products, we define
        \begin{equation}
          k_1^{(m)}(r) \equiv \int_0^{r} U_1^{(m)} R_1^{(m)} r dr.
        \end{equation}
        $k_2^{(m)}$, $k_3^{(m)}$, $j_1^{(m)}$, $j_2^{(m)}$, and
        $j_3^{(m)}$ are defined analogously.  Clearly,
        \begin{equation}
          \lim_{r\to\infty}k_1^{(m)}(r) = K_1^{(m)}
        \end{equation}
        and analogously for the other inner product values.  We
        approximate,
        \begin{equation}
          K_1^{(m)} \approx k_1^{(m)}(\rmax),
        \end{equation}
        for $\rmax$ sufficiently large that these converge to their
        limiting values.  As Figures \ref{f:kfuncs} and \ref{f:jfuncs}
        show, this is indeed the case.

        \ifpdf
	\begin{figure}
          \centering
          \subfigure[$m=1$]{\includegraphics[width=2.35in]{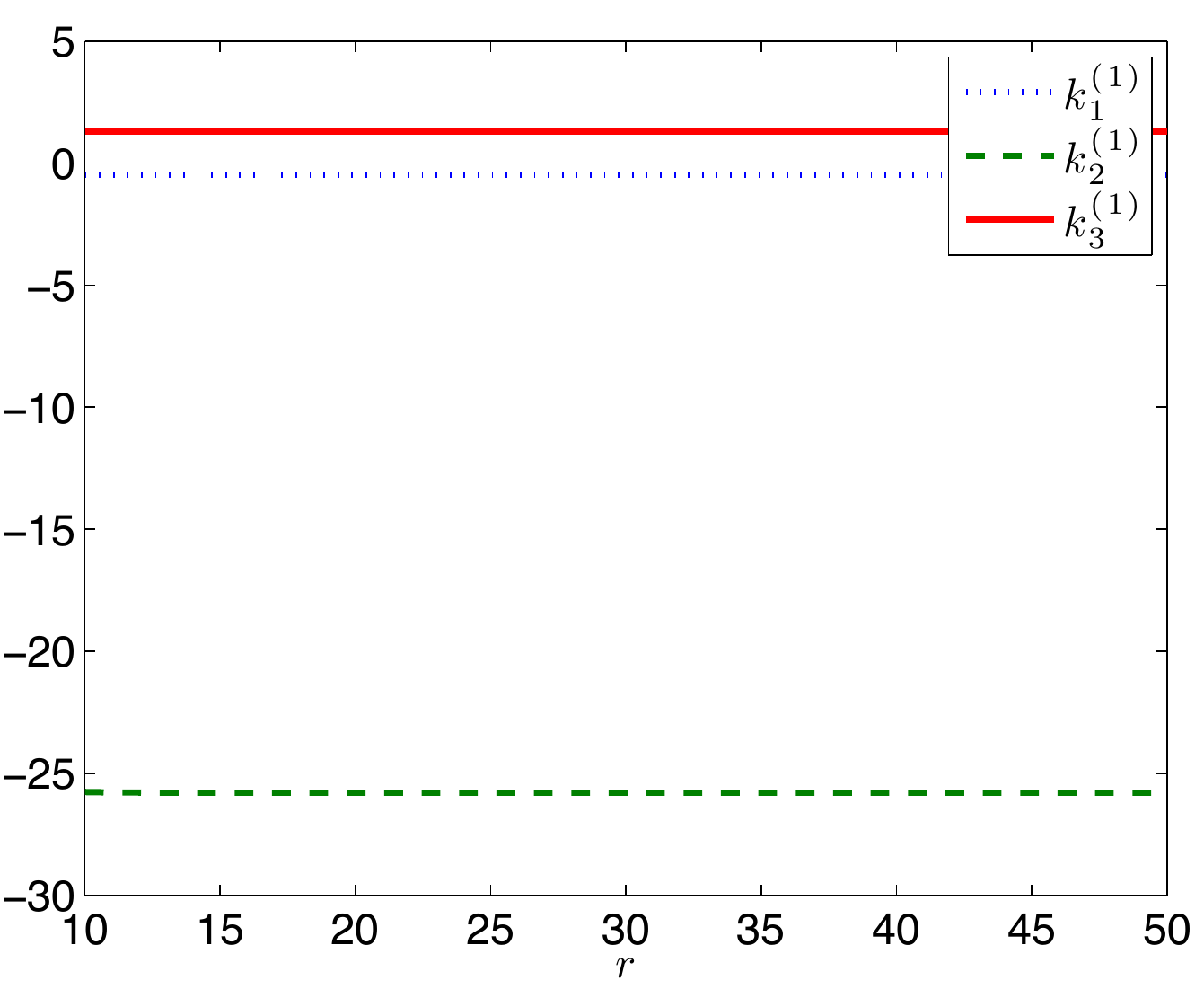}}
          \subfigure[$m=2$]{\includegraphics[width=2.35in]{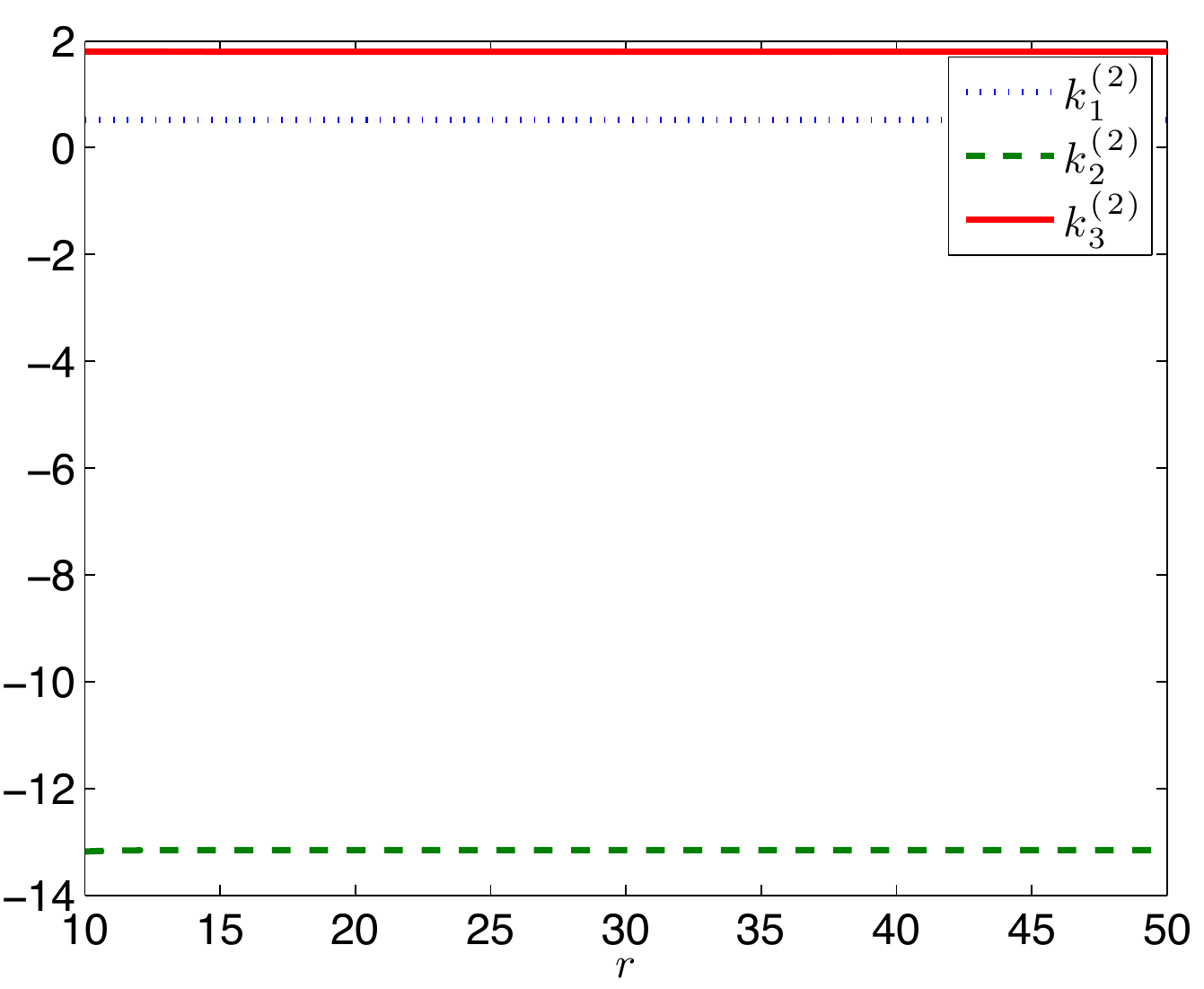}}
          \subfigure[$m=3$]{\includegraphics[width=2.35in]{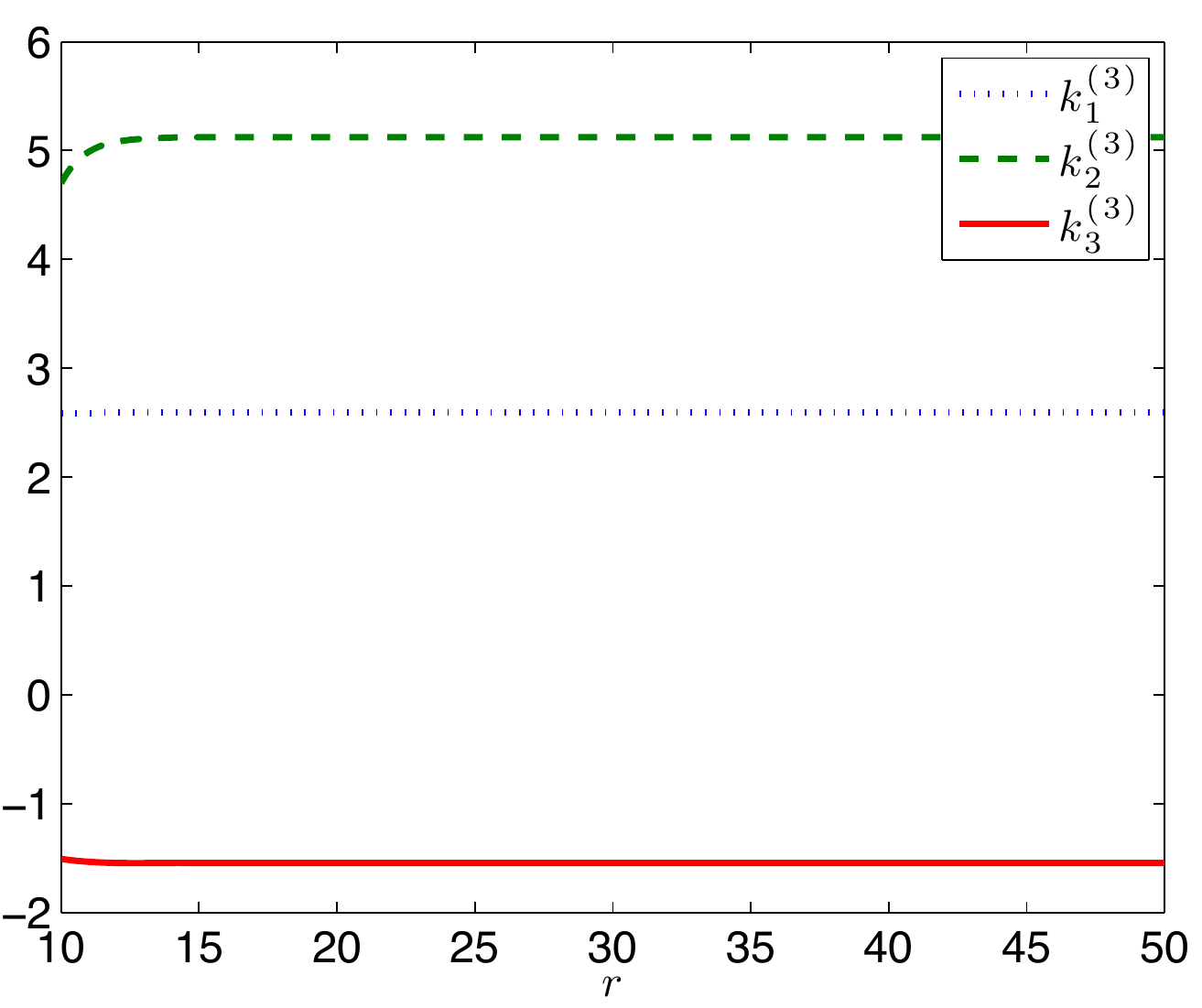}}
          \caption{Convergence of the approximate inner product values
            to their limiting states.}
          \label{f:kfuncs}
	\end{figure}
	
	\begin{figure}
          \centering
          \subfigure[$m=1$]{\includegraphics[width=2.35in]{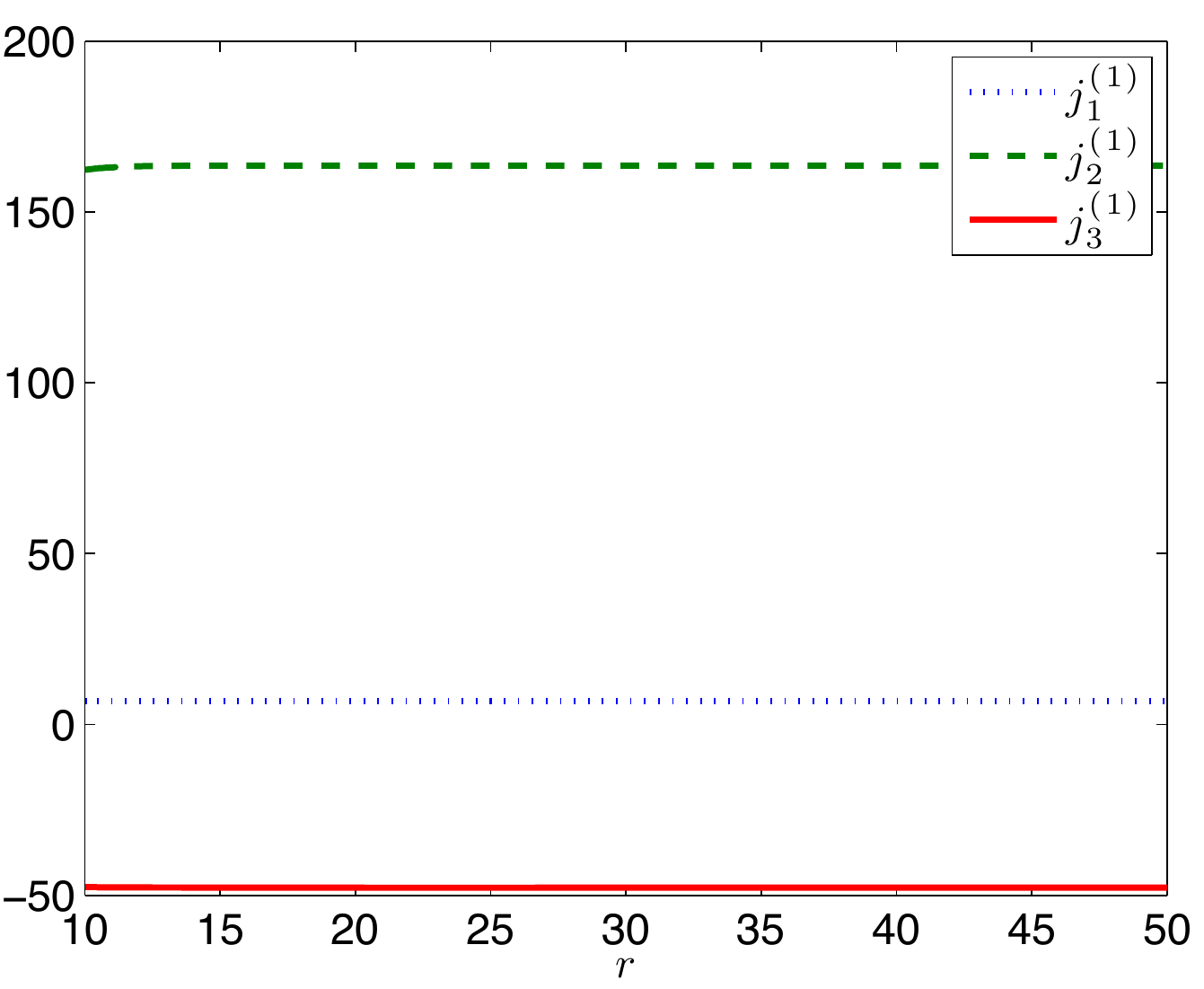}}
          \subfigure[$m=2$]{\includegraphics[width=2.35in]{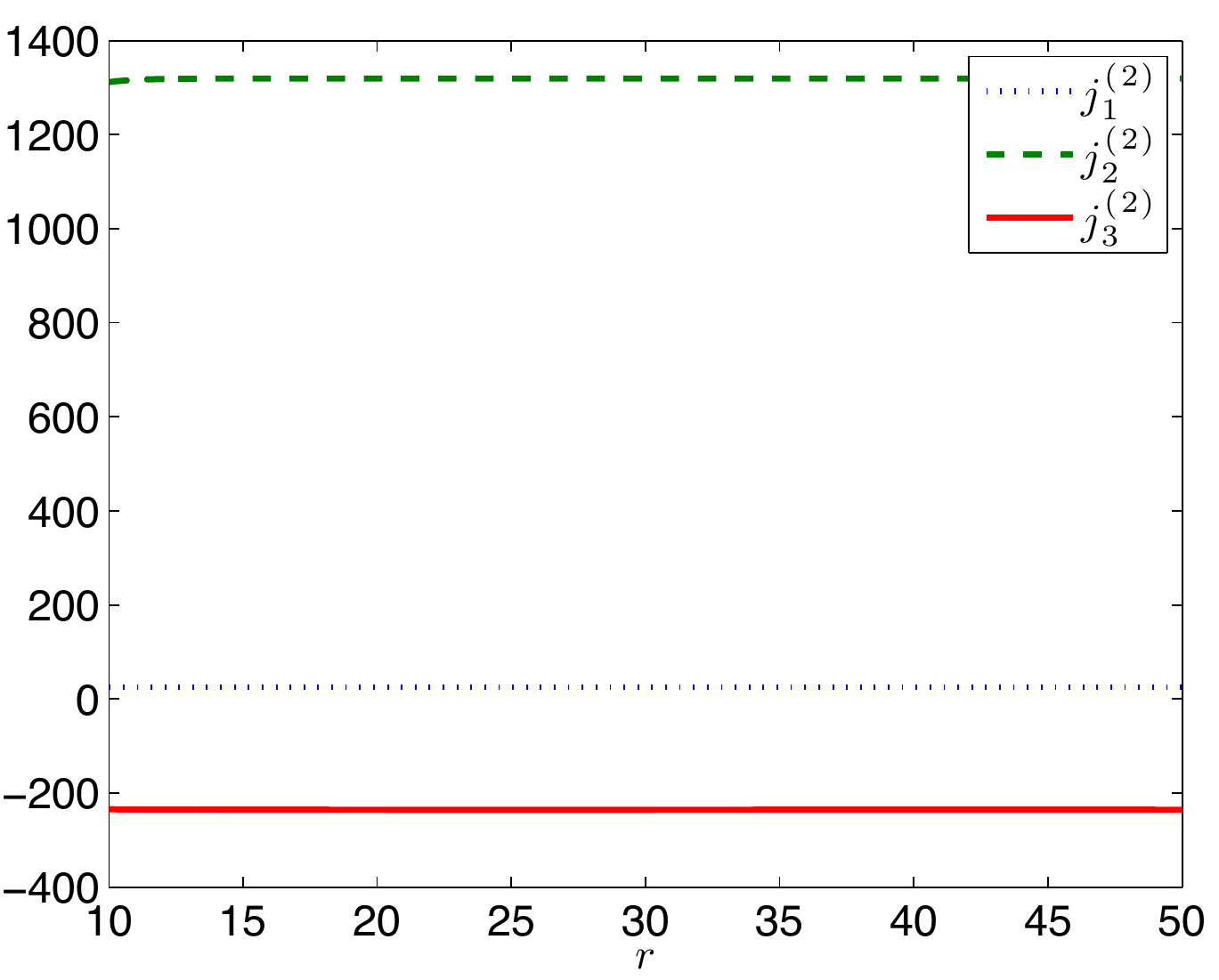}}
          \subfigure[$m=3$]{\includegraphics[width=2.35in]{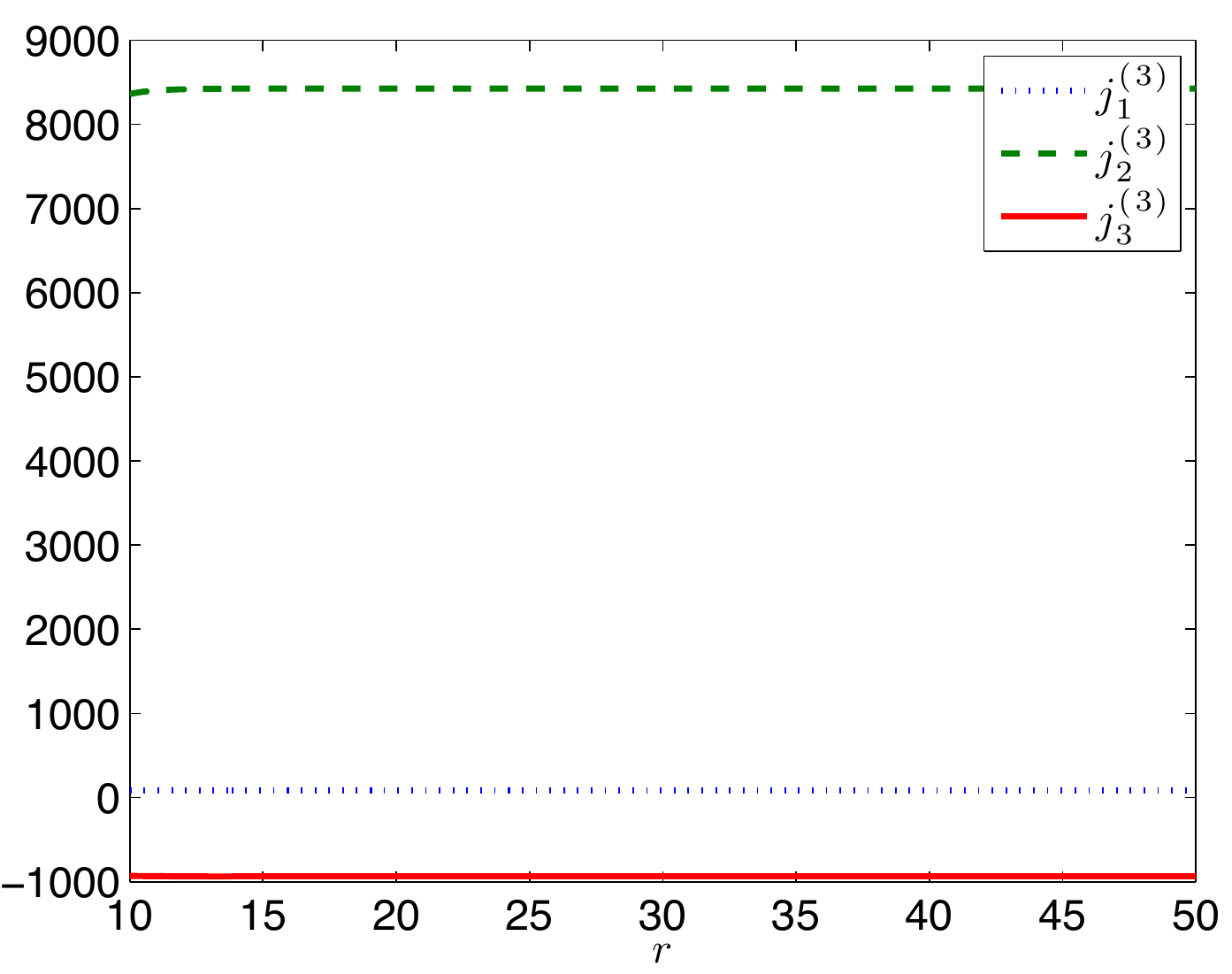}}
          \caption{Convergence of the approximate inner product values
            to their limiting states.}
          \label{f:jfuncs}
	\end{figure}
        \else (No figures in DVI) \fi

        \clearpage \bibliographystyle{plain}
        \bibliography{VortexBlowup}

      \end{document}